\documentclass[mnsc,nonblindrev]{informs3_hide} %

\OneAndAHalfSpacedXI %

\usepackage{amsmath,amsfonts,amssymb}
\usepackage{mathtools}
    \allowdisplaybreaks
\usepackage{bm}
\usepackage[mathscr]{euscript}
\usepackage[inline]{enumitem}
\usepackage{xcolor}
\usepackage{graphicx}
\usepackage{multirow}
\usepackage{booktabs}
\usepackage{microtype}
\usepackage{fix-cm}
\usepackage{physics}
\usepackage{xfrac}
\usepackage{tabularx}

\usepackage{algorithm}  
\usepackage{algorithmicx}  
\usepackage{algpseudocode}
\usepackage{psfrag}
\usepackage{makecell}
\usepackage{tikz}
\usetikzlibrary{arrows.meta}

\DeclareMathOperator{\pr}{\mathbb P}
\DeclareMathOperator{\E}{\mathbb E}
\DeclareMathOperator{\CD}{\mathsf{CoD}}
\DeclareMathOperator{\Var}{\mathsf{Var}}
\DeclareMathOperator{\Cov}{\mathsf{Cov}}

\usepackage{comment}

\usepackage{natbib}
 \bibpunct[, ]{(}{)}{,}{a}{}{,}%
 \def\bibfont{\small}%
 \def\BIBand{and}%

\TheoremsNumberedThrough     %

\EquationsNumberedThrough    %

\MANUSCRIPTNO{MS-SMS-2023-03763.R2}

\begin{document}

\RUNTITLE{Staffing under Taylor's Law}

\TITLE{Staffing under Taylor's Law: A Unifying Framework for Bridging Square-root and Linear Safety Rules}

\ARTICLEAUTHORS{%
\AUTHOR{L. Jeff Hong}
\AFF{Department of Industrial and Systems Engineering, University of Minnesota, Minneapolis, MN 55455, U.S., \EMAIL{lhong@umn.edu}}
\AUTHOR{Weihuan Huang}
\AFF{School of Management \& Engineering, Nanjing University, Nanjing 210093, China, 
\EMAIL{hwh@nju.edu.cn}}
\AUTHOR{Jiheng Zhang, Xiaowei Zhang}
\AFF{Department of Industrial Engineering and Decision Analytics, Hong Kong University of Science and Technology, Clear Water Bay, Hong Kong SAR, \EMAIL{jiheng@ust.hk}, \EMAIL{xiaoweiz@ust.hk}}
}

\ABSTRACT{%
Staffing rules are an essential management tool in service industries for meeting target service levels.
The square-root safety rule, based on the Poisson arrival assumption, has been commonly used.
However, empirical findings suggest that arrivals often exhibit ``over-dispersion'', meaning that the variance exceeds the mean.
In this paper, we develop a new doubly stochastic Poisson process model that captures two key features of over-dispersed arrivals: (i) Taylor's law, which links the variance to the mean through a power-law relationship, and (ii) temporal correlation decay, where the correlation between arrival counts in disjoint time intervals decreases as the time gap grows.
Using this model, we study how over-dispersion affects staffing and derive a closed-form staffing formula to ensure a desired service level.
Our formula shows that the safety level grows as a power of the nominal load. The exponent lies between $1/2$ (the square-root safety rule) and $1$ (the linear safety rule).
It depends on the degree of over-dispersion, and it implies that Taylor's law is the dominant factor in determining staffing levels in heavy traffic.
Extensive numerical experiments with both simulated and real arrival data show that our model and staffing rules significantly outperform various alternatives.
}


\maketitle

\vspace{-18pt}

\section{Introduction}\label{Intro}

According to the World Bank's Development Data Group  (\url{https://data.worldbank.org}), as of 2021, service industries contribute nearly 80\% of the United States' gross domestic product and over 60\% at a global level.
These service industries are predominantly labor-intensive, with staffing costs constituting a major portion of operating expenses. As highlighted by \cite{Gans}, call centers exemplify this pattern, where personnel costs account for 60\%--70\% of total operating expenses. Consequently, determining optimal staffing levels becomes a strategic imperative, requiring careful balance between service quality and operational efficiency. 

The study of staffing under uncertainty originates from \cite{Erlang}'s Erlang C formula and the heavy-traffic analysis of \cite{HalfinWhitt81}.
Their work shows that achieving a target delay probability requires staffing beyond the nominal load $\lambda/\mu$ (the arrival-to-service rate ratio) by a safety level of $\beta \lambda^{1/2}$, where $\beta$ depends on the target delay probability. This is the well-known \emph{square-root safety rule}.

The efficacy of the square-root safety rule hinges critically on the Poisson arrival assumption, where the variance of arrival counts equals their mean over any fixed time period. 
However, empirical evidence from large-scale service systems  increasingly demonstrates that actual arrivals exhibit \emph{over-dispersion}---variation exceeding that of Poisson processes. 
For detailed analyses of this phenomenon, see \cite{Avramidis04}, \cite{ibrahim2016modeling}, and  \cite{OreshkinReegnardLEcuyer16}.

\begin{figure}[t]
  \FIGURE{
    \includegraphics[width=0.5\textwidth]{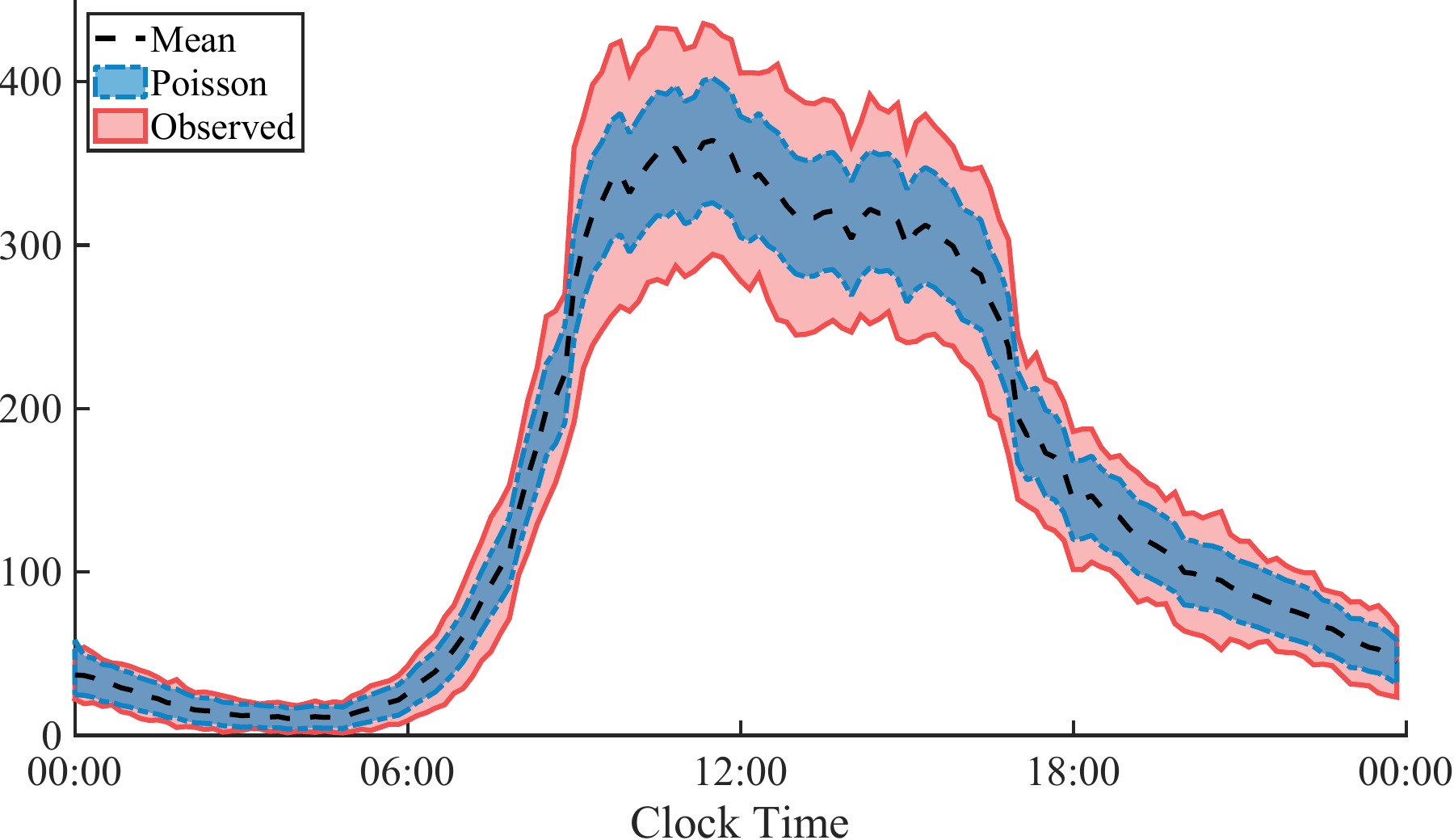} 
  }
  {Daily Arrival Pattern of the NYC 311 Call Center. \label{fig:time-of-day}}
  {The black dashed line shows the mean arrival count per 10-minute interval, calculated from weekday data (July--November 2017). The half-width of the shaded areas represents two standard deviations. 
The red solid boundary represents observed variability in the data, while the blue dash-dotted boundary shows the theoretical variability under a Poisson process assumption.}	
\end{figure}

Figure~\ref{fig:time-of-day} illustrates over-dispersion in arrival data from the New York City (NYC) 311 Call Center. The observed variability is much larger than what a Poisson model would predict, particularly during peak hours (09:00 to 16:00). 
In this setting, the standard square-root safety rule would lead to \emph{under-staffing} that fails to meet the target service level.

This naturally raises three research questions: 
\begin{enumerate}[label=(\arabic*)]
    \item \emph{Is there a fundamental law governing the observed over-dispersion, especially in heavy traffic?}
    \item \emph{If so, can we develop a simple arrival model that captures it?}
    \item \emph{What staffing rule is appropriate under such a law?}
\end{enumerate}

To answer these questions, it is crucial to quantify the degree of over-dispersion in the arrivals. 
Analysis reveals that arrivals to large-scale service systems follow a distinct \emph{dispersion scaling} phenomenon: the standard deviation ($\Sigma$) and mean ($M$) of arrival counts follow the power-law relationship $\Sigma^2 \propto M^{\alpha+1}$, where $0\leq\alpha<1$. 
This phenomenon, known as \emph{Taylor's law}, was first discovered in ecological population studies \citep{Taylor61}. The Poisson process represents a special case where $\alpha=0$.
Empirical validation using NYC 311 Call Center data confirms Taylor's law through linear regression of $\log(\Sigma^2)$ on $\log M$ across varying time intervals. The model achieves remarkable fit ($R^2=97.74\%$) with an over-dispersion coefficient $\alpha\approx0.21$ (Figure~\ref{fig:taylor} in Section~\ref{motivation}).

\subsection{Main Contributions}

While Taylor's law characterizes the magnitude of over-dispersion relative to the mean arrival level, it does not provide a model that describes the arrival dynamics and facilitates staffing decisions.
Our first contribution fills this gap through a novel doubly stochastic Poisson process (DSPP) model.
In this model, the arrival intensity is modeled by a generalized Cox--Ingersoll--Ross (CIR) process \citep{CIR1985} with a dispersion scaling parameter that reflects Taylor's law. 
The model also captures temporal correlation decay, which is common in practice. This dependence can contribute to over-dispersion and has been a main focus of arrival modeling in the literature \citep{ibrahim2016modeling}. 
Finally, our model is parsimonious and includes several prior models as special cases.

Second, we establish a functional central limit theorem (FCLT) for this arrival model in heavy traffic. 
The FCLT reveals an unconventional scaling of $\lambda^{-(\alpha+1)/2}$, rather than the usual square-root scaling. 
The limiting process has an additive form, combining an integrated Ornstein--Uhlenbeck (OU) process  (which inherits parameters from the generalized CIR intensity) with an independent Brownian motion.
This decomposition separates two sources of randomness: intensity uncertainty (captured by the OU process) and conditional Poisson distribution (captured by the Brownian motion). 
The parameter $\alpha$ determines their relative contributions, yielding three over-dispersion scenarios: intensity-dominated, Poisson-dominated, and balanced. 
Furthermore, 
since FCLTs for arrival processes are central to heavy-traffic queueing analysis \citep{Whitt02}, our result provides theoretical basis for extending classical queueing results to systems with over-dispersed arrivals.

Our third contribution is the development of staffing rules for large service systems with over-dispersed arrivals under our DSPP model. 
Using an infinite-server approximation, we obtain a safety level of order $\lambda^{(\alpha+1)/2}$ when arrivals follow Taylor's law. 
Temporal correlation affects the staffing level through the leading constant multiplying this term. 
Thus, in heavy traffic, Taylor's law becomes dominant and largely determines staffing rules.
This formulation generalizes existing staffing principles: when $\alpha=0$, it reduces to a square-root safety rule, and as $\alpha$ varies between 0 and 1, it bridges the gap between square-root and linear safety rules. This unified framework adapts to varying degrees of arrival over-dispersion across different operational settings. We complement the theoretical development with a practical heuristic refinement for implementing the staffing rule in finite-server systems.

Our final contribution provides a practical framework for staffing decisions under non-stationary and over-dispersed arrivals. 
We extend our model to non-stationary settings and develop a computationally efficient maximum likelihood estimation (MLE) method based on heavy-traffic approximations. 
This method, which may be of independent interest, dramatically reduces the well-known computational burden compared to standard MLE of DSPP models.
Extensive numerical experiments demonstrate that our model and staffing rules consistently outperform alternatives across various scenarios, including stationary and non-stationary arrivals, misspecified arrival processes, different staffing interval lengths, and varying service time distributions.

\subsection{Literature Review}

\paragraph{Modeling Over-dispersed Arrivals.}
A growing body of queueing literature addresses arrival modeling that captures over-dispersion, including DSPP models with random business factors \citep{Whitt99}, Poisson mixture models \citep{JongbloedKoole}, extensions incorporating time-varying rates and intraday dependence \citep{Avramidis04,channouf2012normal,OreshkinReegnardLEcuyer16}, and the survey by \cite{ibrahim2016modeling}.

However, most existing models fail to  capture Taylor's law, which reveals a deeper characteristic: a power-law scaling relationship between dispersion and mean arrival rate. 
A notable exception is the DSPP model used by \cite{BassambooRandhawaZeevi10} and \cite{HuChanDong25}. 
Their model specifies intensity as $\lambda + \lambda^{(\alpha+1)/2}Y$, where $Y$ is a zero-mean random variable.
Yet, since this intensity is modeled as a static random variable independent of time, it is better suited for single-stage rather than sequential multi-stage decision-making, as evidenced by their newsvendor-type problem formulations. 
Our paper fills this gap by proposing a continuous-time model that captures Taylor's law while remaining parsimonious (requiring only one additional parameter over their model) and analytically tractable.

Our model also closely relates to that of 
\cite{ZhangHongGlynn}, 
who developed a DSPP model with a CIR arrival rate. Their analysis of data from Israeli and US bank call centers showed that over-dispersion could be explained by 
mesoscopic time scale behavior. We extend their model by incorporating a dispersion scaling parameter to capture Taylor's law.

\paragraph{Staffing under Over-dispersion.}
The square-root safety rule under Poisson arrivals has been extensively studied \citep{HalfinWhitt81,jennings1996server,GarnettMandelbaumReiman02,borst2004dimensioning,mandelbaum2009staffing}.
Staffing under over-dispersed arrivals has gained increasing attention \citep{Whitt06Staffing,BassambooData-driven,KocagaArmonyWard15}, as has staffing under other non-standard arrival processes \citep{GaoZhu18,DawPender18,HeemskerkMandjesMathijsen22,ChenHong23,DawHampshirePender25}.
However, few prior studies developed staffing rules that accommodate Taylor's law in the arrivals.

Our staffing rules share a similar analytical form with those developed by \cite{Maman09}, \cite{BassambooRandhawaZeevi10},  \cite{HuChanDong25} under the DSPP model with intensity $\lambda + \lambda^{(\alpha+1)/2} Y $. However, their staffing objectives differ from our goal of achieving a target quality of service. For example, \cite{BassambooRandhawaZeevi10} and   \cite{HuChanDong25} focus on cost minimization formulated as newsvendor-type problems. 
For fair comparison, we derive a staffing rule with our same objective under their model in the e-companion (Appendix~\ref{StaticDSPP}). Extensive numerical experiments show this staffing rule is significantly less competitive than ours.
Our staffing rule's superior performance stems from our model's more accurate representation of temporal correlation structure in arrivals. Our model, driven by a mean-reverting stochastic process, exhibits decay of temporal correlation: the correlation between arrival counts in separate time intervals decreases with their temporal separation, a pattern widely documented in practice \citep{ibrahim2016modeling}. In contrast, their model's static intensity implies constant correlation regardless of temporal separation.

\cite{SunLiu21} examined staffing under the DSPP model of \cite{ZhangHongGlynn}, showing that the square-root safety rule persists despite over-dispersion. This is because their model's intensity variability does not scale with $\lambda$, so it does not capture Taylor's law.
In contrast, our model's variability scales as a power law in $\lambda$, yielding a safety level of order $\lambda^{(\alpha+1)/2}$ that generalizes the square-root rule. While both analyses use weak convergence arguments \citep{Whitt02,Whitt2007}, the proofs for $\alpha=0$ do not extend to $\alpha \in (0, 1)$. We also extend to non-stationary arrivals and develop a new estimation method based on heavy-traffic approximations.

\smallskip
The rest of the paper is organized as follows. Section \ref{motivation} shows that arrivals in large-scale service systems follow Taylor's law. Section \ref{themodel} introduces a DSPP arrival model with a generalized CIR process to capture Taylor’s law, and establishes its heavy-traffic stochastic-process limit. Section \ref{sec:staffing} develops a staffing rule based on the DSPP and an infinite-server approximation, and then adapts it to finite-server queues. Section~\ref{sec:model-comparison} compares our model with related over-dispersion models and proposes a heavy-traffic-based MLE method for efficient estimation and information-criteria model selection. Sections~\ref{sec:synthetic} and \ref{sec:casestudy} present experiments with stationary and non-stationary arrivals that support the theory and show our staffing rule outperforms alternatives. Section~\ref{conclusion} concludes, and the e-companion contains technical proofs and additional experiments.

\section{Over-dispersion and Taylor's Law}\label{motivation}

The term ``over-dispersion'' in queueing theory literature addresses the stochastic variability of an arrival process within a queueing system. 
It refers to the phenomenon in which the variance of the arrival count within a given time period exceeds its mean. 
Formally, let $A = \{A(t):t\geq 0\}$ represent an arrival process, such that $A(0)=0$ and $A(t)$ denotes the arrival count in the time interval $(0, t]$. For any $0 \leq s < t$, let $A(s, t] = A(t) - A(s)$ represent the arrival count in the time interval $(s, t]$. 
We define $A(s, t]$ as \emph{over-dispersed} if its \emph{coefficient of dispersion} (CoD) is greater than one, i.e., \[\CD[A(s, t]] := \frac{\Var[A(s,t]]}{\E[A(s, t]]} > 1.\] 
If this property holds for any $0 \leq s< t$, we consider the arrival process $A$ to be over-dispersed. 
In contrast, if $A$ is a Poisson process, the arrival count in any time period follows a Poisson distribution, resulting in equal variance and mean (i.e., $\CD[A(s, t]] =1$ for any $0 \leq s < t$).

It is well-established that arrival processes in real-world queueing systems---particularly those with a high intensity of customer arrivals in modern service industries, such as telephone call centers---often exhibit a considerable degree of over-dispersion \citep{Avramidis04,OreshkinReegnardLEcuyer16}. The presence of over-dispersion carries significant implications for operations management within the corresponding queueing systems.
Intuitively, to maintain the same quality of service, an increased service capacity is necessary to handle the elevated stochastic variability in the system. Assuming that individual servers possess essentially identical service efficiency---indicating that the service capacity is primarily determined by the number of servers---an over-dispersed arrival process would necessitate a higher staffing level compared to a Poisson arrival process in order to preserve equivalent system performance (e.g., delay probability and mean waiting time).

In general, a practical and often optimal rule of thumb for staffing that balances operating costs and quality of service takes the form of the sum of a ``base level'' and a ``safety level''. 
The former is set to match the mean demand for service, while the latter represents the additional number of servers deployed to mitigate stochastic variability in the system.

To make this concrete, consider an infinite-server queue with a stationary arrival process and exponential service times with unit rate. 
Let $N$ denote the arrival count in one unit of time. 
In this setting, the steady-state number of customers in the system has the same distribution as $N$. 
This link is useful for staffing an $n$-server system, because the delay probability can be approximated by $\pr(N \geq n)$. 
We can therefore choose $n$ to meet a target delay probability, resulting in the base-plus-safety staffing rule
\begin{equation}\label{eq:base-plus-safety}
\underbrace{\E[N]}_{\text{base}} \ +\ \underbrace{\beta \sqrt{\Var[N]}}_{\text{safety}}.
\end{equation}
Here, $\beta>0$ is a constant that determines the degree of hedging against the stochastic variability of $A$, and it is specified according to an appropriate trade-off between staffing costs and quality of service. 
This managerial insight is further supported by \cite{BassambooRandhawaZeevi10} through a newsvendor-type model. Evidently, the staffing rule \eqref{eq:base-plus-safety} is closely related to the level of over-dispersion measured by $\CD[N]$, which we will now demonstrate using two examples.

\smallskip
\begin{example}[Square-root Safety Rule]\label{example:sr-staffing}
Assume that the arrival process $A$ is a Poisson process with a constant rate $\lambda$. In this case, $\Var[N] = \E[N] = \lambda$, which results in $\CD[N] = 1$ and reduces the staffing rule \eqref{eq:base-plus-safety} to $\lambda + \beta \sqrt{\lambda}$. This is the square-root safety rule.
Its asymptotic optimality has been established in a heavy-traffic regime as $\lambda\to\infty$ across various settings, extending beyond standard multi-server queues \citep{HalfinWhitt81,GarnettMandelbaumReiman02,borst2004dimensioning}. 
\end{example}

\smallskip
\begin{example}[Linear Safety Rule]\label{example:linear-staffing}
Assume that the arrival process $A$ is a doubly stochastic Poisson process with a random rate $\Lambda = \lambda G$, where $G$ is a random variable with unit mean and finite variance \citep{Whitt99,BassambooRandhawaZeevi10}. 
Using the tower property of conditional expectations, 
it is straightforward to demonstrate that $\E[N] = \lambda$ and $\Var[N] = \lambda^2 \Var[G] + \lambda$. Consequently, both $\CD[N]$ and the safety level of the staffing rule \eqref{eq:base-plus-safety} grow at a rate of $\lambda$ as $\lambda\to\infty$. This is referred to as the \emph{linear safety rule}.
\end{example}

\smallskip

Despite their widespread use, the arrival models employed in these two examples fail to accurately capture the degree of over-dispersion experienced in heavy-traffic environments. 
Instead, the coefficient of dispersion (CoD) of the arrival count within a given time period often exhibits a growth rate of $\lambda^\alpha$, where $\alpha$ is a constant within the range $(0, 1)$. 
For instance, in the arrival process of the NYC 311 Call Center, we observe this growth rate of the CoD with $\alpha\approx 0.21$. 

\begin{figure}[t]
    \FIGURE{
    \includegraphics[width=0.45\textwidth]{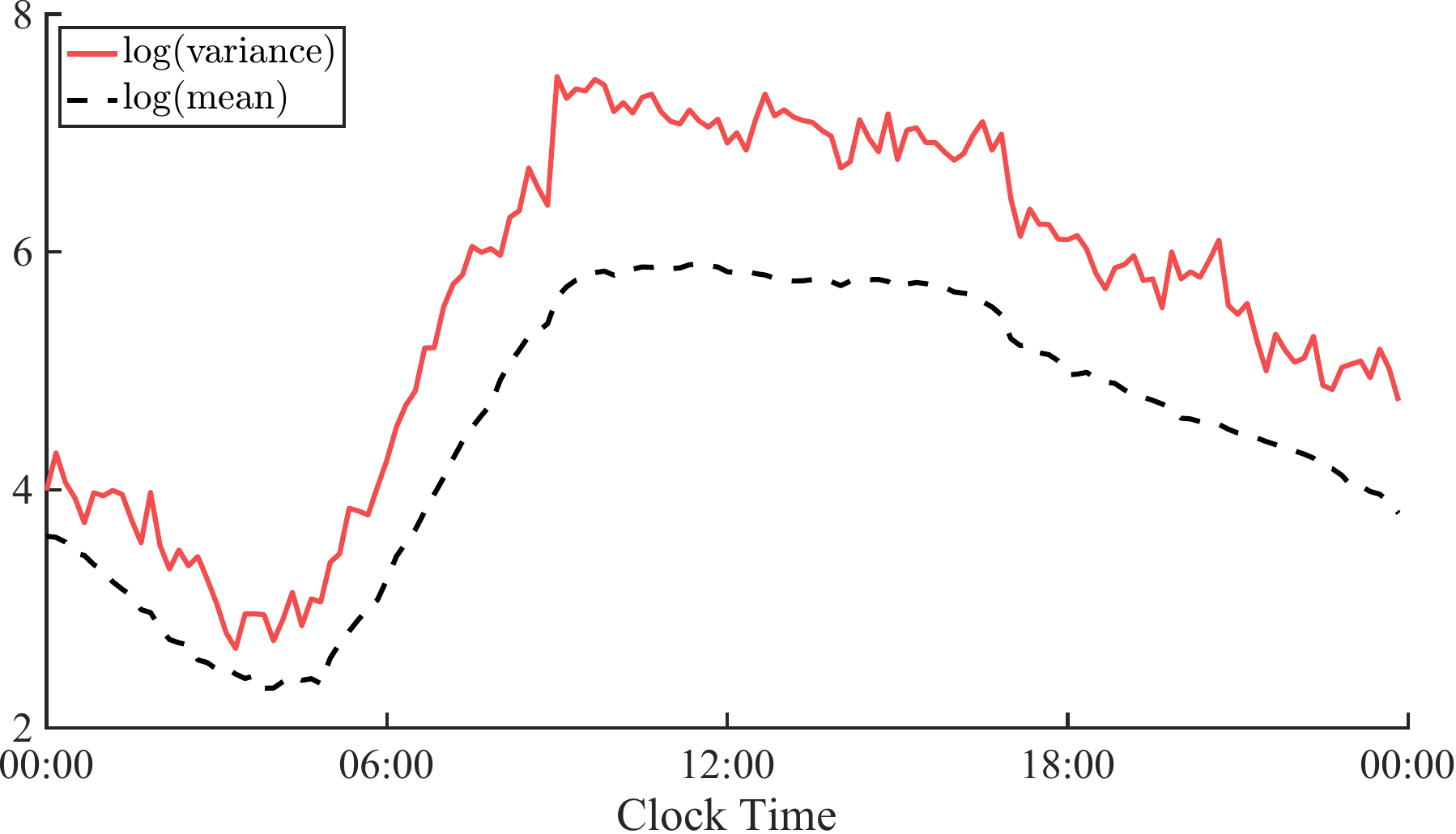} 
    \includegraphics[width=0.45\textwidth]{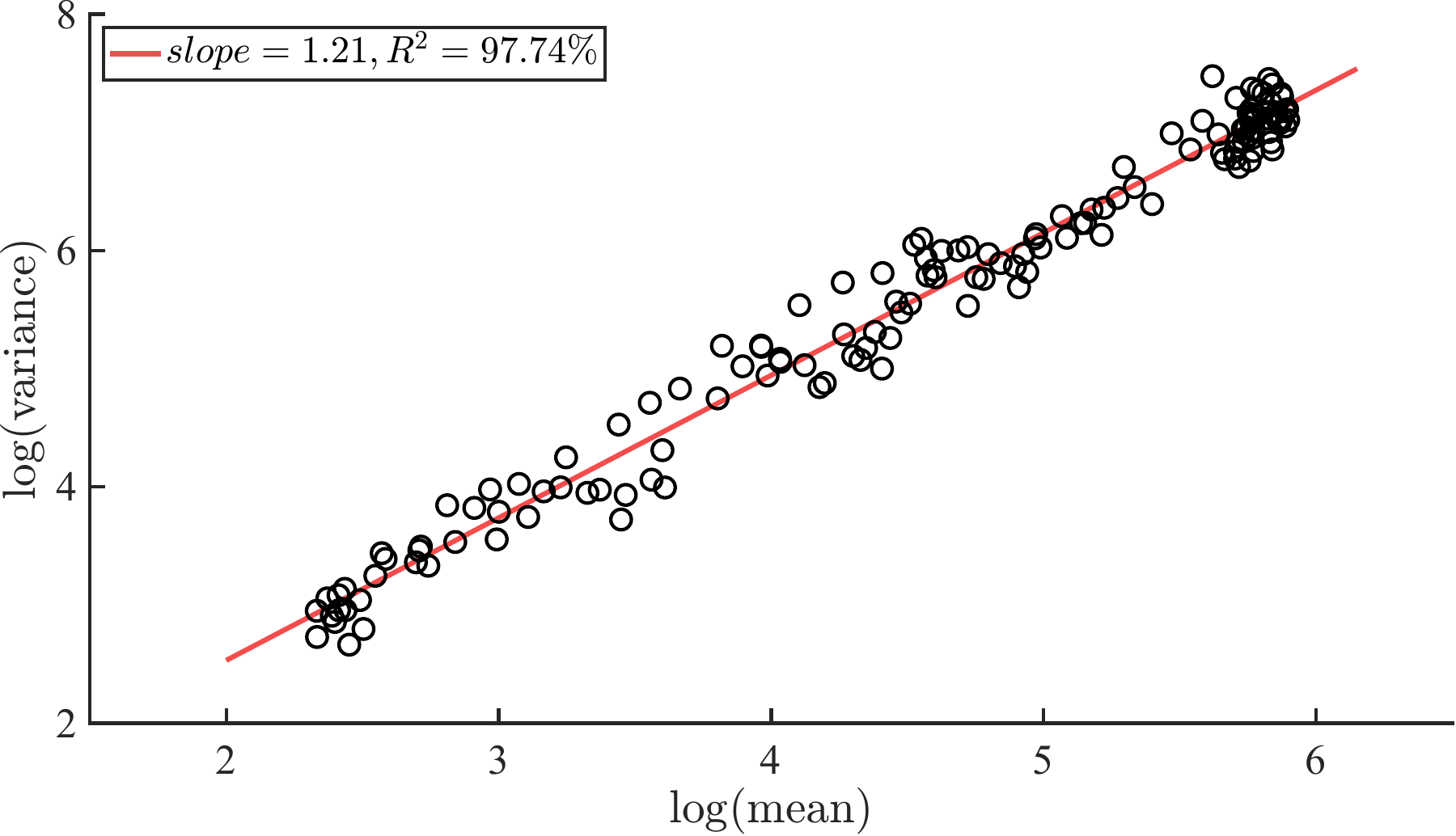}
    }
    {Taylor's Law and the Arrival Process of the NYC 311 Call Center. \label{fig:taylor}}
    {Left: The variance and mean of the arrival count in 10-minute time periods, plotted on a logarithmic scale. Right: Linear regression between the variance and mean after performing a logarithmic transformation.} 
\end{figure}

Taking into consideration the notable time-of-day effect displayed by the arrival process (refer to Figure~\ref{fig:time-of-day}), we treat the arrivals within a short 10-minute time period as stationary---a common practice in the field \citep{GreenKolesarWhitt07}. 
Each day is divided into 10-minute intervals, resulting in 144 periods in total. 
Subsequently, the mean and variance of the arrival count in each time period are calculated based on the 96 weekdays in our dataset from July to November 2017. 
These values are plotted on a logarithmic scale in the left pane of Figure~\ref{fig:taylor}. 
Upon performing a logarithmic transformation on both the means and variances, we conduct a linear regression analysis with the log-variance as the dependent variable. Remarkably, the linear regression exhibits a strong fit, as evidenced by the coefficient of determination, $R^2 = 97.74\%$ (see the right pane of Figure~\ref{fig:taylor}).

Let $N_{\Delta}$ represent the arrival count within a time period of length $\Delta$. The linear regression result implies the following relationship:
\begin{equation}\label{LR}
\log({\Var[N_{\Delta}]}) = (1 + \alpha)\log({\E[N_{\Delta}]}) + c,
\end{equation}
where $\alpha$ and $c$ are constants, with $\alpha \approx 0.21$ in Figure~\ref{fig:taylor}. 
Consequently,
\begin{equation}\label{eq:Taylor}
    \Var[N_{\Delta}] \propto (\E[N_{\Delta}])^{\alpha+1} \qq{or}  \CD[N_{\Delta}] \propto (\E[N_{\Delta}])^{\alpha}.
\end{equation}

The power-law relationship described above is known as Taylor's law, which was originally identified in the field of ecology by \cite{Taylor61}. This law characterizes the variance in the number of individuals belonging to a species distributed across space and/or time as a power-law function of the mean. Over the years, Taylor's law has been observed in a wide range of disciplines, including life sciences, social networks, internet traffic, and stock markets \citep{Taylor19}.
Empirical evidence similar to that depicted in Figure~\ref{fig:taylor} has been discovered in various healthcare facilities \citep{Maman09,HuChanDong25,HongLiuLuoXie22}.

Let $\lambda = \E[N_{\Delta}]/\Delta$ represent the mean arrival count per time unit. In accordance with the base-plus-safety staffing rule given by Equation~\eqref{eq:base-plus-safety}, it is plausible that if the arrival process adheres to Taylor's law (as expressed in Equation~\eqref{eq:Taylor}), an effective staffing rule for queueing systems in which servers provide service at a unit rate would take the form:
\begin{equation}\label{eq:taylor-staffing}
\lambda + \beta \lambda^{\frac{\alpha+1}{2}}.    
\end{equation}
Here, $\beta>0$ is a constant determined according to a specific target quality of service. The primary objective of this paper is to provide a theoretical foundation for this staffing rule.

\section{A Doubly Stochastic Poisson Process Model} \label{themodel}

In this section, we propose a new model for arrivals that exhibit the over-dispersion characterized by Taylor's law.\footnote{We focus on modeling \emph{stationary} arrivals. Although the arrival rate is stochastic, it does not exhibit ``time-varying'' characteristics in the conventional sense.  
In queueing literature, the term ``time-varying'' typically refers to arrival rates that are \emph{deterministic} functions, often capturing seasonality like ``time-of-day'' or ``day-of-week'' effects \citep{ibrahim2016modeling}.
Instead, our stochastic arrival rate has a stationary distribution, fluctuating around an equilibrium level. Nevertheless, our model can be extended to incorporate deterministic trends (Section~\ref{sec:casestudy}).} 
We begin by fixing a complete probability space $(\Omega, \mathscr{F}, \pr)$ and a right-continuous, complete information filtration $\mathbb{F}=\{\mathscr{F}_t: t\geq 0\}$. 
We model the arrival process $A=\{A(t):t\geq 0\}$ as a doubly stochastic Poisson process with intensity $X=\{X(t):t\geq 0\}$ (refer to Definition~\ref{DSPPdef} for a formal definition). 
Intuitively, $X(t)$ serves as the conditional arrival rate at time $t$, such that $\E[A(t, t+\Delta]|\mathscr{F}_t] \approx X(t) \Delta$ for any sufficiently small $\Delta>0$.

\smallskip
\begin{definition}[Doubly Stochastic Poisson Process]\label{DSPPdef}
Let $A$ be an $\mathbb{F}$-adapted counting process, and let $X$ be a non-negative, integrable, $\mathbb{F}$-predictable process. If, conditional on $\mathscr{F}_s$, $A(s, t]$ follows the Poisson distribution with mean $\int_s^t X(r) \dd{r}$ for all $0 \leq s < t$, then $A$ is referred to as a doubly stochastic Poisson process with intensity $X$. 
\end{definition}

\begin{figure}[t]
    \FIGURE{	\includegraphics[width=0.33\textwidth]{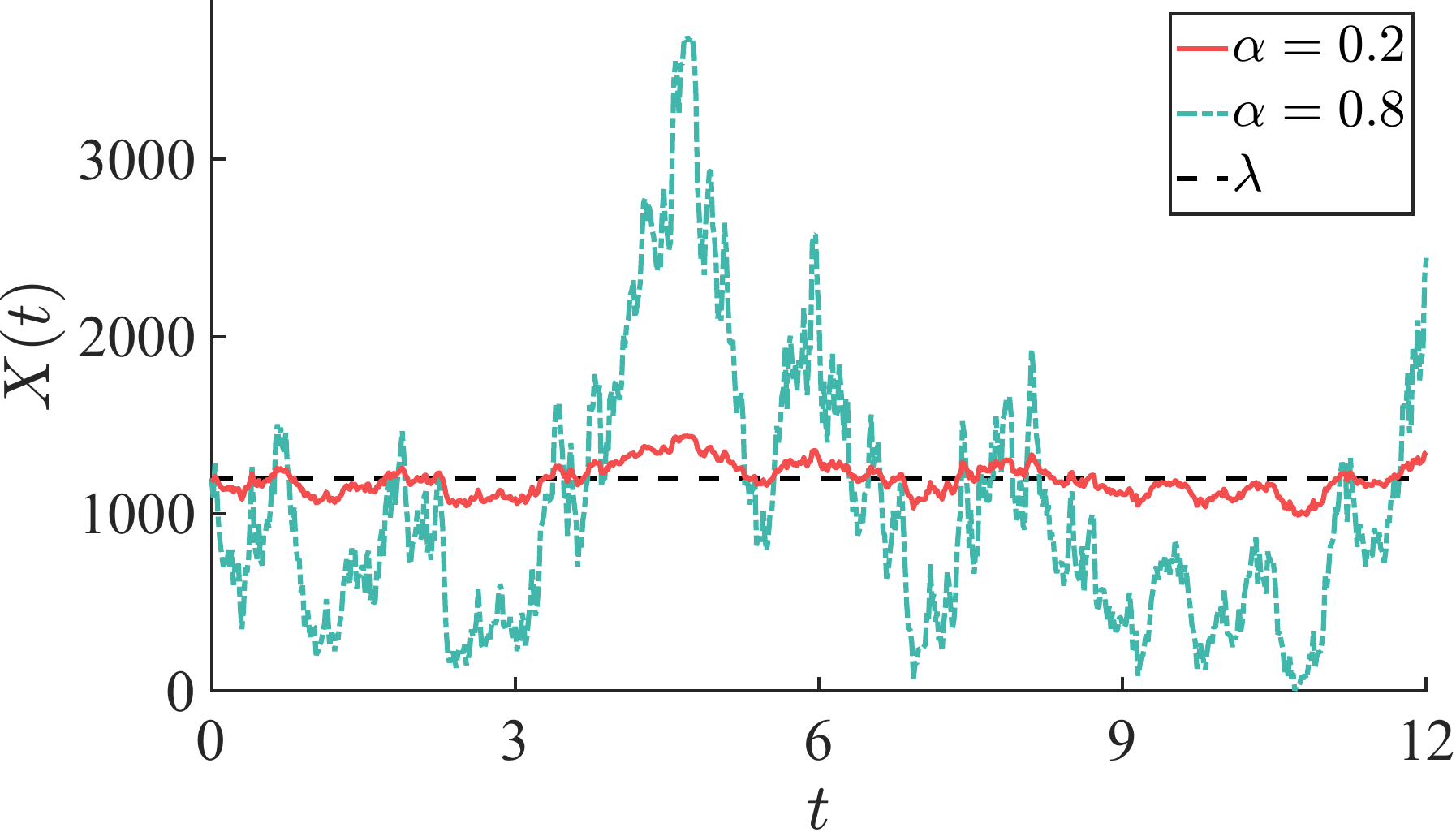}
	\includegraphics[width=0.33\textwidth]{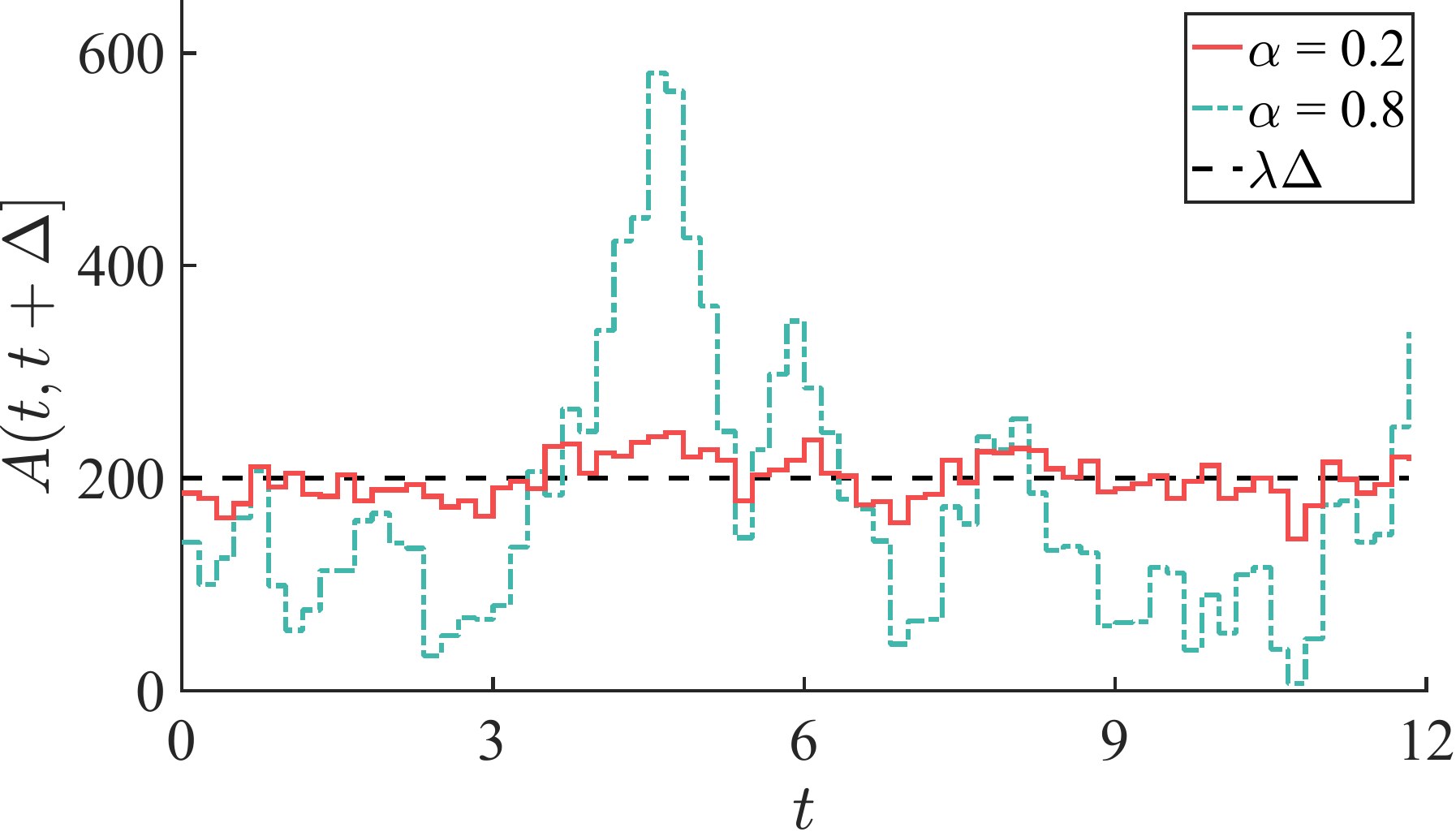}
	\includegraphics[width=0.33\textwidth]{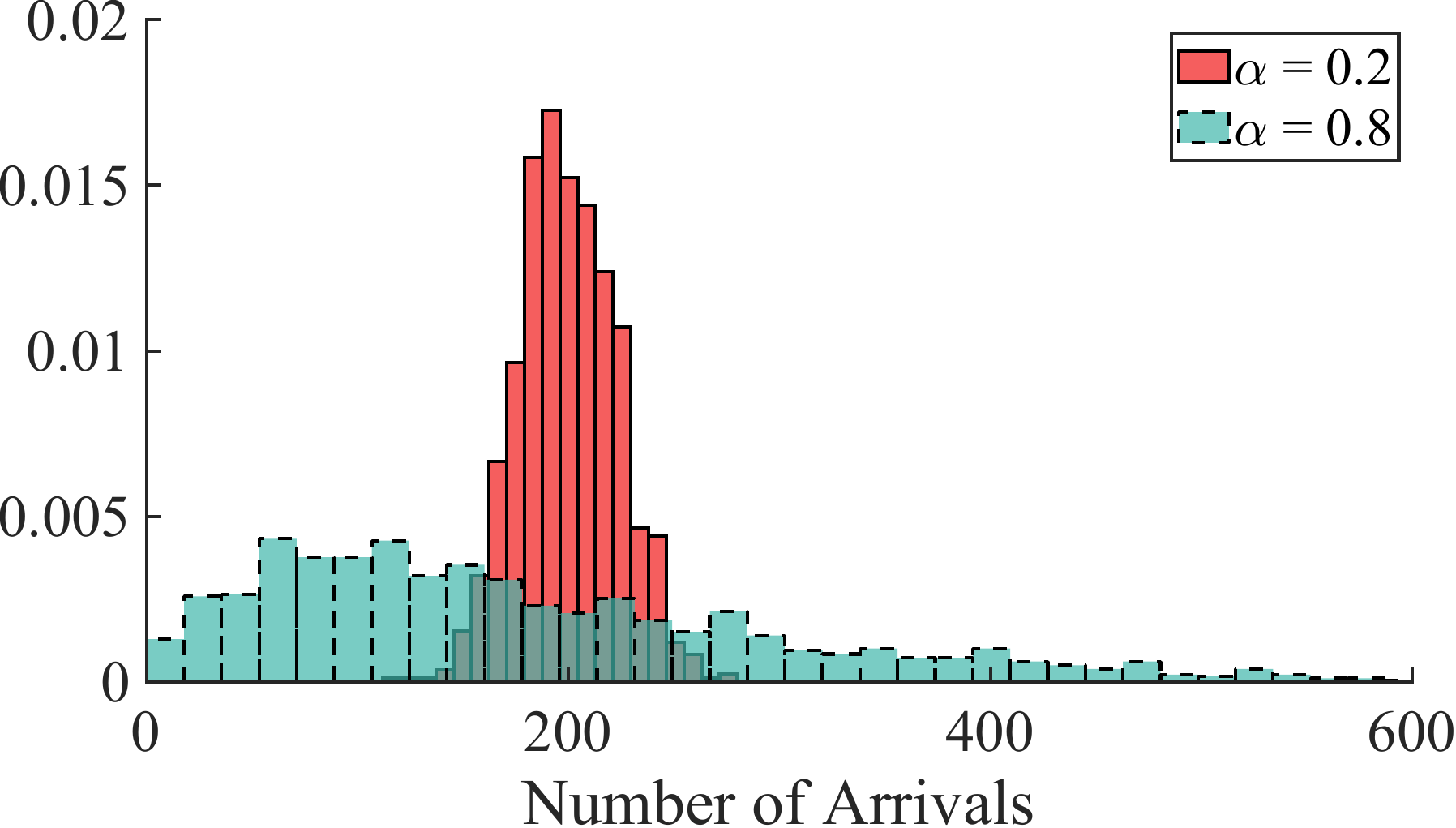}
	}
    {
    Simulated Intensity Process $X(t)$ and 
    Arrival Count $A(t, t+\Delta]$. 
    \label{fig:sample-path}    }
    {Left: Sample path of $X(t)$. Middle: Sample path of $A(t, t+\Delta]$. Right: Distribution of $A(t, t+\Delta]$ when $X(0)$ follows the stationary distribution of $X$ (see Proposition~\ref{prop:squre-diffussion}). The parameters are  $\lambda=1200$, $\kappa = 1.0$, $\sigma = 2.2$, and $\Delta=1/6$ (10 minutes).
    }
\end{figure}

\smallskip
We assume that $X$ is a diffusion process satisfying the following stochastic differential equation:
\begin{equation}\label{eq:CIR}
\mathrm{d} X(t) = \kappa\left( \lambda -X(t) \right) \dd{t} + \sigma  \sqrt{\lambda^\alpha X(t)} \dd{B(t)},
\end{equation}
with $X(0) >0$ almost surely. Here, $B=\{B(t):t\geq 0\}$ is a standard Brownian motion adapted to $\mathbb{F}$, with parameters $\lambda>0$, $\kappa>0$, $\sigma>0$, and $\alpha\geq 0$, each having a clear interpretation. (We may use $X_\lambda$ when emphasizing dependence on $\lambda$.) Figure~\ref{fig:sample-path} shows typical sample paths of $X(t)$ and the corresponding arrival counts over intervals of length $\Delta$.

We call the process \eqref{eq:CIR} a \emph{generalized CIR} process, because it reduces to the CIR process when $\alpha=0$. 
Note that the drift term $\kappa(\lambda - X(t))$ in  Equation~\eqref{eq:CIR} is \emph{mean-reverting}---it is positive when $X(t) < \lambda$ and negative when $X(t) > \lambda$. 
We can interpret $\lambda$ as an equilibrium level (i.e., long-run average) of $X(t)$ and $\kappa$ as the speed at which $X(t)$ is driven towards $\lambda$. 
The mean-reversion property ensures the process to be ergodic, thereby having a stationary distribution.

As $X(t)$ approaches zero, the diffusion term $\sigma \sqrt{\lambda^\alpha X(t)}$ in Equation~\eqref{eq:CIR} vanishes, causing random fluctuations to become negligible compared to the mean-reverting force. This prevents $X(t)$ from becoming negative, a desirable feature for modeling arrival rates. For rigorous mathematical statements, see Proposition~\ref{prop:squre-diffussion}.

Evidently, the over-dispersion of the arrival process $A$ stems from its doubly stochastic structure. 
Indeed, applying the law of total variance, we obtain the following inequality for any $0\leq s < t$:
\begin{align*}
\Var\big[A(s, t]\big] \geq 
\E\biggl[\Var\biggl[A(s,t]\bigg|\int_s^t X(r)\dd{r}\biggr]\biggr] = \E\biggl[\int_s^t X(r)\dd{r}\biggr] = \E\big[A(s, t]\big].
\end{align*}
The two identities arise from the conditional Poisson distribution and the tower property of conditional expectations, respectively. 
Our model exhibits a notable feature in comparison to other over-dispersed arrival models found in the literature \citep{Avramidis04,Zhang13,OreshkinReegnardLEcuyer16}: the inclusion of the parameter $\alpha$. 
This parameter, referred to as the \emph{dispersion scaling} parameter, explicitly captures Taylor's law \eqref{eq:Taylor}, which governs the relationship between the increase in mean arrival rate and the accompanying level of over-dispersion.

Loosely speaking, 
$X(t)$ has an order of magnitude of $\lambda$ in the long run. 
Consequently, the diffusion term $\sigma \sqrt{\lambda^\alpha X(t)}$ in Equation \eqref{eq:CIR}, which represents the standard deviation of the process, has an order of magnitude of $\sqrt{\lambda^{\alpha+1}}$. 
This observation implies that the contribution of $X(t)$ to the variance of $A(t)$ has an order of magnitude of $\lambda^{\alpha+1}$.
Drawing from this heuristic argument, we can infer that the CoD of the arrival count in one time unit increases at a rate of $\lambda^\alpha$ as $\lambda\to\infty$, which corroborates Taylor's law. Theorem~\ref{theo:CoD} provides a detailed explanation of this result.

In addition to its parsimonious nature, ease of interpretation, and ability to capture Taylor's law of dispersion scaling, the newly proposed arrival model offers significant analytical tractability.
Specifically, it enables explicit calculations of both the finite-time CoD and the asymptotic distribution of the arrival process $A$ (upon appropriate scaling). 
Furthermore, the model facilitates the derivation of simple staffing rules and provides ample opportunities for extensions.

Before delving into these calculations, we first introduce several properties of the intensity process $X$, which will be employed in the subsequent analysis of the arrival process $A$. 
Throughout the remainder of the paper, the following assumption will be imposed.

\begin{assumption}\label{asp:Feller}
$2\kappa \lambda^{1-\alpha} \geq \sigma^2$. 
\end{assumption}

\begin{remark}
The significance of Assumption~\ref{asp:Feller} can be described as follows. 
It ensures that when $X(t)$ is in close proximity to the origin, the drift tendency towards $\lambda$ is sufficiently strong to prevent the origin from being reachable. 
Without this guarantee, $X(t)$ might approach the origin within a finite time frame, necessitating additional specification of the boundary behavior at the origin for the process to be well-defined. 
To conduct this analysis, one can verify the boundary classification criteria for diffusion processes as detailed in \citet[Chapter 15]{KarlinTaylor81}.
\end{remark}

\begin{proposition}\label{prop:squre-diffussion}
If Assumption~\ref{asp:Feller} holds, then 
the following statements are valid. 
\begin{enumerate}[label=(\roman*)]
    \item \label{prop:SD-positive} If $X(0)>0$, then almost surely $X(t)>0$ for all $t>0$.
    \item \label{prop:SD-stationary} $X$ has a unique stationary distribution $\pi$, which is $\mathsf{Gamma}(\frac{2\kappa}{\sigma^2\lambda^{\alpha-1}}, \frac{2\kappa}{\sigma^2 \lambda^\alpha})$, that is, the gamma distribution having the shape parameter $\frac{2\kappa}{\sigma^2\lambda^{\alpha-1}}$ and the rate parameter $\frac{2\kappa}{\sigma^2 \lambda^\alpha}$. 
    \item \label{prop:SD-ergodic} $\pr(X(t)\in \cdot |X(0)) \to \pi$ in total variation as $t\to\infty$.  
\end{enumerate}
\end{proposition}

\subsection{Coefficient of Dispersion}

We present the CoD of the arrival process using the Laplace transform of $A(t)$,  leveraging our model's analytical tractability. 
Full derivations and the proof of Theorem~\ref{theo:CoD} appear in  Appendix~\ref{ec-sec:laplace}.

\begin{theorem}\label{theo:CoD}
Suppose Assumption~\ref{asp:Feller} holds. 
If $X(0)\sim \pi$, then for any $t>0$, we have 
    \[
    \E_\pi[A(t)] = \lambda t \qq{and} 
    \Var_\pi[A(t)] = \lambda t +   \frac{\sigma^2 \lambda^{\alpha +1} t}{\kappa^2}\left( 1- \frac{1-e^{-\kappa t}}{\kappa t}\right), 
    \]
    and thus, 
    \[
    \CD_\pi[A(t)] = 1 + \frac{\sigma^2 \lambda^{\alpha} }{\kappa^2}\left( 1- \frac{1-e^{-\kappa t}}{\kappa t}\right).
    \]
\end{theorem}

\smallskip
\begin{remark}
Calculations for the scenario where $X(0)$ is fixed at an arbitrary $x>0$ can be carried out in a similar manner. 
For instance, we have:
\[
\E[A(t) | X(0) = x] = \lambda t + \kappa^{-1}(1-e^{-\kappa t})(x-\lambda).
\]
However, the expressions for the variance and CoD in this scenario are more cumbersome and lack simple structure. 
Consequently, we have chosen to omit the details for the sake of brevity.
\end{remark}
\smallskip

It is worth noting that in Theorem~\ref{theo:CoD}, the expression for $\Var[A(t)]$ aligns with the law of total variance. 
The first term is equal to $\E\bigl[\Var\bigl[A(t)\big|\int_0^t X(s)\dd{s}\bigr]\bigr]$, reflecting the contribution from the (conditional) Poisson distribution, while the second term is equal to $\Var\bigl[\E\bigl[A(t)\big|\int_0^t X(s)\dd{s}\bigr]\bigr]$, representing the contribution from the randomness of the intensity process.

In a light-traffic environment where $\lambda$ is small, the magnitude of the first term dominates the second term, implying that the stochastic variability of $A(t)$ is essentially ``Poisson-like''. 
In such an environment, the square-root safety rule is expected to perform well.
However, in a heavy-traffic environment---which is often the case for service systems in practice---the randomness of the intensity process becomes the dominant factor in the composition of $\Var[A(t)]$. 
Consequently, the square-root safety rule may result in a significantly understaffed system. 
A conceptually similar discussion, framed within a newsvendor-type model, can be found in \cite{BassambooRandhawaZeevi10}. 
We refer to Section~\ref{sec:staffing} for further comments on this work.

We now connect our model to Taylor's law. 
Consider the case where the intensity process $X$ starts from its stationary distribution $\pi$. Theorem~\ref{theo:CoD} implies that for all sufficiently large $\lambda$,
\begin{equation}\label{eq:CoD-model}
\CD_\pi[A(t)] \approx \frac{\sigma^2 \lambda^\alpha}{\kappa^2} \left(1 - \frac{\kappa t - (\kappa t)^2 /2 }{\kappa t}\right) = \frac{\sigma^2   t^{1-\alpha}}{2\kappa} (\E_\pi[A(t)])^\alpha,
\end{equation}
recovering the relationship \eqref{eq:Taylor}. 
(Here, we have also applied Taylor's expansion to approximate $1-e^{-\kappa t}$.) 
Thus, the parameter $\alpha$ in the intensity process \eqref{eq:CIR} is essentially identical to the parameter that defines Taylor's law \eqref{eq:Taylor}, and we do not differentiate between them in the sequel.

\subsection{Value Range of $\alpha$}

Our theoretical analysis up to this point has been conducted without explicit restrictions on the value of $\alpha$, as long as the values of the other parameters $\lambda$, $\kappa$, and $\sigma$ satisfy Assumption~\ref{asp:Feller}. 
However, in the sequel, we will impose constraints on the value of $\alpha$.

\begin{assumption}\label{asp:alpha}
$0\leq \alpha < 1$. 
\end{assumption}

We impose this assumption for several reasons. First, extensive empirical evidence from both natural sciences \citep{EislerBartosKertesz08} and social sciences \citep{HongLiuLuoXie22} shows that when a random quantity follows Taylor's law, its $\alpha$ typically lies in $(0,1)$. The NYC 311 Call Center data exemplifies this pattern, yielding $\alpha = 0.21$ (see Figure~\ref{fig:taylor}).

Second, the objective of this paper is to develop a staffing rule for large-scale service systems experiencing heavy traffic. 
To achieve this, we consider an asymptotic regime that sends $\lambda \to \infty$, which implies that Assumption~\ref{asp:Feller} must hold for all sufficiently large $\lambda$. 
As a result, $\alpha$ cannot take values greater than one, because otherwise, $2\kappa \lambda^{1-\alpha} \to 0$, which is smaller than $\sigma^2$, as $\lambda \to \infty$.

Third, while Assumption~\ref{asp:Feller} allows for $\alpha=1$, the asymptotic distribution of $A(t)$ in this boundary case is distinct and warrants a separate treatment that may be of theoretical interest in its own right. 
To simplify the exposition of this paper, we choose to focus on the more common case where $\alpha<1$. However, it is worth noting that our analysis in the subsequent sections can be extended to the boundary case without any significant difficulty.

\subsection{Asymptotic Normality}\label{AN}

In order to formulate staffing rules for heavy-traffic environments, we first analyze the asymptotic distribution of the arrival process $A(t)$ as $\lambda\to\infty$. 
As the expectation of $A(t)$ grows unbounded with increasing $\lambda$, it is crucial to apply appropriate scaling to obtain a non-degenerate limiting distribution. 
A common practice involves scaling a ``centered'' arrival process by a factor of $\lambda^{-1/2}$.
For example, if $A(t)$ is a Poisson process with a constant intensity $\lambda$, then the scaled process (i.e., $\lambda^{-1/2}[A(t) - \lambda t]$) converges in a proper sense to a standard Brownian motion as $\lambda\to\infty$.
Consequently, in this example, for any fixed $t$ and  sufficiently large $\lambda$,  $A(t)$ approximately follows a normal distribution with a standard deviation of $(\lambda t)^{1/2}$.

However, when $A(t)$ is a doubly stochastic Poisson process with its intensity process specified by Equation~\eqref{eq:CIR}, the application of the conventional scaling factor $\lambda^{-1/2}$ generally proves ineffective, resulting in a degenerate limiting distribution. 
This failure can be attributed to the significantly increased variance, which exhibits an order of $\lambda^{\alpha+1}$, as shown by Theorem~\ref{theo:CoD}. 
This contrasts with the order-$\lambda$ variance observed in the standard Poisson process example.
In light of this observation, we propose a new scaling scheme, which involves scaling the centered $A(t)$ by a factor of the same order as $1/\sqrt{\Var[A(t)]}$. 
This corresponds to a scaling factor of order $\lambda^{-(\alpha+1)/2}$, which is anticipated to yield a non-degenerate limit, as detailed in Theorem~\ref{thm:FCLTofA} (with proof in Appendix~\ref{sec:proof-scaled-arrival} and  empirical evidence in Appendix~\ref{sec:self-norm}).

To enhance clarity, we introduce a subscript to explicitly indicate the dependence on $\lambda$ for stochastic processes.
Henceforth, $A(t)$ and $X(t)$ will be denoted as $A_\lambda(t)$ and $X_\lambda(t)$, respectively.
Moreover, we define their scaled versions as:
\begin{equation}\label{eq:hat_processes}
\widehat{X}_\lambda(t)  =  \lambda^{-\frac{\alpha+1}2}[X_\lambda(t)-\lambda] \quad \textrm{ and } \quad \widehat{A}_\lambda(t)  =  \lambda^{-\frac{\alpha+1}2}[A_\lambda(t)-\lambda t].
\end{equation}
In addition, for any $T>0$, we define the function space $\mathcal{D}[0,T]$ as the set of all right-continuous real-valued functions with left limits on the interval $[0,T]$, equipped with the standard Skorohod $J_1$ topology. We use  ``$\Rightarrow$'' to denote weak convergence (i.e., convergence in distribution). 
For more details on the space $\mathcal{D}$ and weak convergence, we refer to \cite{Billingsley99} and \cite{Whitt02}.

\begin{theorem}\label{thm:FCLTofA}
Let $B$ be the Brownian motion  that drives the intensity process $X_\lambda(t)$, as defined in Equation~\eqref{eq:CIR}.
Suppose that
\begin{enumerate*}[label=(\roman*)]
    \item Assumptions~\ref{asp:Feller}--\ref{asp:alpha} hold,
    \item $X_\lambda(0)$ is independent of $B$, 
    \item $\sup_{\lambda \geq 1} \E[\widehat{X}_\lambda^4(0)]<\infty$, and 
    \item $\widehat{X}_\lambda(0)\Rightarrow X_\dagger $ as $\lambda\to\infty$ for some random variable $X_\dagger$.
\end{enumerate*} 
Let $U$ be the following Ornstein–Uhlenbeck (OU) process 
\begin{equation}\label{OU}
\dd{U(t)} = -\kappa U(t) \dd{t} +\sigma \dd{B(t)}, 
	\end{equation}
with the initial condition $U(0)=X_\dagger$. 
Then, for any $T>0$, we have  
\begin{equation}\label{thm:FCLT:eq1}
\widehat{A}_\lambda \Rightarrow \widehat{A}_\infty\quad \mbox{in $\mathcal{D}[0,T]$}
\end{equation}
as $\lambda\rightarrow\infty$,
where 
\begin{equation}\label{FCLTofA}
\widehat{A}_\infty(t) = \int_0^t U(s) \dd{s} + \widetilde{B}(t) \mathbb{I}\{\alpha = 0\},
\end{equation}
and $\widetilde{B}$ is a standard Brownian motion independent of $B$ (and thus independent of $\int_0^\cdot U(s) \dd{s}$).
\end{theorem}

Although $X_\dagger$ (defined in condition (iii) of Theorem~\ref{thm:FCLTofA}) does not appear explicitly in $\widehat A_\infty$, it affects the distribution of $\widehat A_\infty$ through $U$. For any $t > 0$, the distribution of  $U(t)$ depends on both the Gaussian transition distribution $\pr(U(t)\in\cdot | U(0))$ and the OU process's initial condition. The transition distribution is characterized by $\kappa$ and $\sigma$ from Equation~\eqref{OU}, while Theorem~\ref{thm:FCLTofA} shows that $X_\dagger$ determines the initial condition. Thus, the distribution of $\widehat A_\infty$ depends on $X_\dagger$.

It can be readily demonstrated that if $U(0)=X_\dagger$ follows a normal distribution, then the integrated OU process, $\bigl\{\int_0^t U(s)\, \dd{s} : t\geq 0\bigr\}$, is a Gaussian process. 
Furthermore, due to the independence between this integrated OU process and the Brownian motion $\widetilde{B}$, as well as the expression \eqref{FCLTofA}, the limiting process $\widehat{A}_\infty$ is also a Gaussian process.

For example, 
let us examine a representative case regarding the initial condition of the intensity process $X_\lambda$. 
Suppose $X_\lambda$ is initialized with the stationary distribution $\pi$, given by $\mathsf{Gamma}(\frac{2\kappa}{\sigma^2\lambda^{\alpha-1}}, \frac{2\kappa}{\sigma^2 \lambda^\alpha})$, as specified in Proposition~\ref{prop:squre-diffussion}. 
By calculating the Laplace transform of $\widehat{X}_\lambda(0)$, it can be readily shown that $\widehat{X}_\lambda(0) \Rightarrow \mathsf{Normal}(0, \frac{\sigma^2}{2\kappa})$ as $\lambda\to\infty$. This limiting distribution is also the stationary distribution of $U$. 
Through a routine calculation as demonstrated on page 113 of \cite{Glasserman03}, we obtain the following result for all $t\geq 0$:
\begin{equation}\label{eq:iOU-dist-pi}
    \int_0^t U(s)\dd{s} \sim \mathsf{Normal}\biggl(0, \frac{\sigma^2 t}{\kappa^2}\left( 1- \frac{1-e^{-\kappa t}}{\kappa t}\right)\biggr).
\end{equation}
It then follows from Theorem~\ref{thm:FCLTofA} that for any $t\geq 0$,
\begin{equation}\label{eq:A_inf-dist-pi}
    \widehat{A}_\infty(t) \sim \mathsf{Normal}\biggl(0, \frac{\sigma^2 t}{\kappa^2}\left( 1- \frac{1-e^{-\kappa t}}{\kappa t}\right) + t \mathbb{I}\{\alpha = 0\} \biggr).
\end{equation}

\subsection{Over-dispersion Trichotomy} \label{sec:trichotomy}

The parameter $\alpha$ typically takes values in $(0,1)$ in practice, but the case $\alpha=0$ has special theoretical significance as it bridges our model with the Poisson model for arrival processes. When $\alpha=0$, the arrival process $A_\lambda$ is scaled by $\lambda^{-1/2}$, and the limiting process $\widehat{A}_\infty$ in Equation~\eqref{FCLTofA} comprises two independent components: an integrated OU process and a Brownian motion $\widetilde{B}$.

When $A_\lambda$ degenerates to a Poisson process with constant intensity $\lambda$, its scaled process $\widehat{A}_\lambda$ (with $\alpha=0$) converges to a Brownian motion. This suggests that the Brownian motion $\widetilde{B}$ in Equation~\eqref{FCLTofA} originates from the conditional Poisson distribution in our arrival model, while the integrated OU process in Equation~\eqref{FCLTofA} stems from the intensity process $X_\lambda$. 

For $\alpha=0$ and $\sigma>0$, the stochastic variability of both the conditional Poisson distribution and the intensity process are comparable in magnitude as $\lambda\to\infty$. We refer to this situation as the \emph{balanced} scenario. 
When $\sigma=0$, the intensity process \eqref{eq:CIR} degenerates to the ordinary differential equation $\dd{X(t)} = \kappa(\lambda - X(t))\dd{t}$, so $X(t) = \lambda + (X(0) - \lambda) e^{-\kappa t}$. With $X(0) = \lambda$, this reduces our model to the Poisson model. In this \emph{Poisson-dominated} scenario, the limiting process $\widehat{A}_\infty$ (with $\alpha=0$) reduces to $\widetilde{B}$.
For typical values $\alpha \in (0,1)$ and $\sigma>0$, the intensity process dominates the doubly stochastic structure as $\lambda \to \infty$, so the Brownian motion $\widetilde{B}$ drops out and $\widehat{A}_\infty$ reduces to the integrated OU process. We refer to this as the \emph{intensity-dominated} scenario.

Table~\ref{tab:arrivals} summarizes the three scenarios. In the balanced scenario ($\alpha=0$, $\sigma>0$), the square-root scaling $\lambda^{-1/2}$ applies as for Poisson arrivals, but over-dispersion inflates the asymptotic variance, requiring a larger safety level (see Section~\ref{sec:staffing}).

\begin{table}[t]
\TABLE{Trichotomy of Over-dispersion in Arrivals. \label{tab:arrivals}}
{\begin{tabular}{@{}ccccccc@{}}
\toprule
  && \multicolumn{1}{c}{\begin{tabular}[c]{@{}c@{}}Intensity-dominated\\ ($0 < \alpha < 1$; $\sigma>0$)\end{tabular}} && \multicolumn{1}{c}{\begin{tabular}[c]{@{}c@{}}Balanced\\ ($\alpha = 0$; $\sigma>0$)\end{tabular}} && \multicolumn{1}{c}{\begin{tabular}[c]{@{}c@{}}Poisson-dominated\\ ($\alpha=0$; $\sigma= 0$, $X(0) = \lambda$)\end{tabular}} \\ 
  \midrule
~Scaling Factor  && $\lambda^{-(\alpha+1)/2}$ && $\lambda^{-1/2}$ && $\lambda^{-1/2} \vphantom{\bigg[1\bigg]}$ \\  
~$\widehat{A}_\infty(t)$ && $\displaystyle \int_0^t U(s) \dd{s}$            && $\displaystyle\int_0^t U(s) \dd{s} + \widetilde{B}(t)$    &&   $\widetilde{B}(t) \vphantom{\bigg[1\bigg]}$  \\
\bottomrule
\end{tabular}
}
{}
\end{table}

\smallskip
\begin{remark}
Beyond revealing the over-dispersion trichotomy,
the FCLT for our DSPP model in Theorem~\ref{thm:FCLTofA} has two key implications. 
First, it facilitates heavy-traffic queueing analysis under our model. FCLTs for arrival processes are fundamental to establishing diffusion limits in queueing systems \citep{Whitt02}, as shown by our infinite-server queue analysis in Section~\ref{sec:IS-apprx}. Given the Poisson process's central role in queueing theory, Theorem~\ref{thm:FCLTofA} lays a foundation for the extension of classical queueing results to over-dispersed arrivals. 
Second, Theorem~\ref{thm:FCLTofA}  enables efficient model estimation. We develop a new MLE method in Section~\ref{sec:model-estimation} that leverages this FCLT to circumvent a key computational challenge: the high-dimensional integration over the intensity process's joint distribution required by standard MLE for DSPPs. 
Additionally, this estimation method facilitates model selection using information criteria (Section~\ref{sec:model-selection}).
\end{remark}

\section{Alpha Safety Rules}\label{sec:staffing}

In this section, we derive staffing rules to achieve target service quality under our proposed DSPP arrival model. Section~\ref{sec:IS-apprx} develops an infinite-server approximation for the delay probability, our service quality metric. Section~\ref{sec:staff_inf} uses this to derive the basic alpha safety rule, a simple formula for safety staffing. This rule yields key managerial insights, confirming that safety levels should scale as $\lambda^{(\alpha+1)/2}$ in heavy traffic when arrivals follow Taylor's law with dispersion parameter $\alpha$, as suggested by Expression~\eqref{eq:taylor-staffing}. Section~\ref{sec:staff_refine} presents the refined alpha safety rule, a heuristic enhancement that maintains the basic rule's structure but modifies the leading constant of the $\lambda^{(\alpha+1)/2}$ term.

\subsection{Infinite-server Approximation} \label{sec:IS-apprx}

Our objective is to find the staffing level $n$ such that the delay probability---the probability that an arriving customer finds all $n$ servers busy---approximately equals a target $\varepsilon$, for a given mean arrival rate $\lambda$.
Since the exact distribution of the number-in-system for a finite-server queue is generally intractable, we apply an infinite-server approximation, a standard approach for large-scale systems \citep[Chapter~10]{Whitt02}.

Let $Q_\lambda(t)$ denote the number of customers at time $t$ in an infinite-server system where customers arrive according to our DSPP model \eqref{eq:CIR}, with $Q_\lambda(0) = 0$.
We use $\pr(Q_\lambda(t) > n)$ to approximate the delay probability of the $n$-server counterpart, so the staffing problem reduces to finding $n$ such that
\begin{equation}\label{eq:inf-serv-approx-obj}
    \pr(Q_\lambda(t) > n) \approx \varepsilon.
\end{equation}
This requires characterizing the distribution of $Q_\lambda(t)$ for large $\lambda$, which we develop as follows.

To characterize its asymptotic distribution in a heavy-traffic regime, we consider a scaled version of $Q_\lambda(t)$, similar to our analysis of the arrival process $A_\lambda(t)$ in Theorem~\ref{thm:FCLTofA}.
Specifically, we subtract from $Q_\lambda(t)$ a quantity having the same order as the mean of $Q_\lambda(t)$ as $\lambda \to \infty$, and then  
divide the difference by another quantity having the same order as the standard deviation of $Q_\lambda(t)$.  

Let us first examine the mean of $Q_\lambda(t)$. 
Note that 
\begin{equation}\label{Q}
Q_\lambda(t) = A_\lambda(t) - D_\lambda(t),
\end{equation}
where $D_\lambda(t)$ denotes the departure process, that is, the number of customers that have completed service and departed from the system by time $t$. 
It follows from Theorem~\ref{theo:CoD} that under the stationary distribution $\pi$ of the intensity process $X_\lambda(t)$, 
we have $\E_\pi[A_\lambda(t)] = \lambda t$. 
Moreover, let $S$ denote the service time of a customer and let $F$ denote its cumulative distribution function (CDF). 
Then,
\begin{equation}\label{eq:mean_Q}
\E_\pi[D_\lambda(t)] = \E_\pi\left[\int_0^t \mathbb{I}\{S\leq t-z\} \dd{A_\lambda(z)} \right] = \int_0^t \pr(S\leq t-z) \lambda\dd{z} = \lambda\int_0^t F(t-z)\dd{z}.
\end{equation}
Therefore, with $\bar F  \coloneqq 1-F$, 
\[\E_\pi[Q_\lambda(t)]=\E_\pi[A_\lambda(t)]-\E_\pi[D_\lambda(t)]= \lambda \int_0^t \bar{F}(t-z)\dd{z}.\]

While straightforward calculations of the standard deviation of $Q_\lambda(t)$ under the distribution $\pi$ is unavailable, it is intuitively reasonable to expect that the arrival and departure processes share the same level of stochastic variability as $\lambda \to \infty$. 
This is because the service time distribution is independent of $\lambda$.
As a result, following the implications of Theorem~\ref{theo:CoD}, we scale $Q_\lambda(t)$, after centering it according to Equation~\eqref{eq:mean_Q}, by a factor of $\lambda^{-(\alpha+1)/2}$: 
\begin{equation}\label{eqn:hatQ}
\widehat{Q}_\lambda(t) \coloneqq \lambda^{-\frac{\alpha+1}2}  \left[ Q_\lambda(t) -\lambda \int_0^t \bar{F}(t-z)\dd{z} \right]. 
\end{equation}

Theorem~\ref{thm:QQ} shows that the scaled process $\widehat{Q}_\lambda(t)$ converges to a non-degenerate limit as $\lambda\to\infty$, verifying the above heuristic argument. 
The proof (Appendix~\ref{sec:proof-scaled-num-in-sys}) relies on the heavy-traffic analysis of infinite-server queueing systems with general arrival processes in \cite{Borovkov1967}.

\begin{theorem}\label{thm:QQ}
Let $B$ be the Brownian motion  that drives the intensity process $X_\lambda(t)$, as defined in Equation~\eqref{eq:CIR}.
Suppose that
\begin{enumerate*}[label=(\roman*)]
    \item Assumptions~\ref{asp:Feller}--\ref{asp:alpha} hold,
    \item $X_\lambda(0)$ is independent of $B$, 
    \item $\sup_{\lambda \geq 1} \E[\widehat{X}_\lambda^4(0)]<\infty$, 
    \item $\widehat{X}_\lambda(0)\Rightarrow X_\dagger $ as $\lambda\to\infty$ for some random variable $X_\dagger$, and 
    \item the function $G(t)\coloneqq \int_0^t F(t-z)\dd{z}$ is H\"older continuous, that is, $|G(t_1)-G(t_2)|< c_1 |t_1-t_2|^{c_2}$ for some constants $c_1>0$ and $c_2>0$.
\end{enumerate*} 
Let $\widehat{A}_\infty$ be the process defined in Equation~\eqref{FCLTofA}, and let
$\xi$ be an independent, zero-mean Gaussian process with the following covariance function: 
\begin{equation}\label{eq:xi-cov}
\E[\xi(t)\xi(t+u)] = \int_0^t F(t-z)\bar{F}(t+u-z)\dd{z}, 
\end{equation}
for all $t, u\geq 0$.
Then, for any $T>0$, we have
\begin{equation}\label{eqn:eq13}
\widehat{Q}_\lambda  \Rightarrow  \widehat{Q}_\infty  \quad \mbox{in $\mathcal{D}[0,T]$} 
\end{equation}
as $\lambda\rightarrow\infty$, where 
\begin{equation}\label{limitofQ}
\widehat{Q}_\infty(t) = \int_0^t \bar{F}(t-z)\dd{\widehat{A}_\infty(z)} + \xi(t)\mathbb{I}\{\alpha = 0\}.  
\end{equation}
\end{theorem}
\smallskip

Applying the definition of $\widehat{A}_\infty(t)$ in Equation~\eqref{FCLTofA} to Equation~\eqref{limitofQ}, we obtain
\begin{equation}\label{Qinfty}
\widehat{Q}_\infty(t) = \underbrace{\int_0^t \bar{F}(t-z) U(z) \dd{z}}_{\coloneqq \,\widehat{Q}_\infty^{\mathsf{GCIR}^+}(t)} + \biggl[\underbrace{\int_0^t \bar{F}(t-z) \dd{\widetilde{B}(z)} + \xi(t)}_{\coloneqq \,\widehat{Q}_\infty^{\mathsf{Pois}}(t)}\biggr] \mathbb{I}\{\alpha = 0\}.
\end{equation}
Equation~\eqref{Qinfty} generalizes classical heavy-traffic limits for infinite-server queues \citep{Whitt82}. Consistent with the over-dispersion trichotomy in Table~\ref{tab:arrivals}, the structure of $\widehat{Q}_\infty(t)$ varies across the three scenarios:
\begin{enumerate*}[label=(\roman*)]
\item In the Poisson-dominated scenario ($\alpha=0$, $\sigma=0$), the OU process $U$ degenerates to 0, so $\widehat{Q}_\infty(t) = \widehat{Q}_\infty^{\mathsf{Pois}}(t)$, recovering the classical results of \cite{Iglehart65} and \cite{Borovkov1967}.
\item In the balanced scenario ($\alpha=0$, $\sigma>0$), both terms contribute, and the over-dispersion term $\widehat{Q}_\infty^{\mathsf{GCIR}^+}(t)$ amplifies the variance beyond the Poisson case, requiring a higher safety level under the same square-root scaling.
\item In the intensity-dominated scenario ($\alpha\in(0,1)$), $\widehat{Q}_\infty(t) = \widehat{Q}_\infty^{\mathsf{GCIR}^+}(t)$ (The sign ``+" refers to the positive $\alpha$), confirming that the queueing variability is driven primarily by the intensity process, and staffing is governed by Taylor's law.
\end{enumerate*}

\subsection{Basic Alpha Safety Rule}\label{sec:staff_inf}  

To determine an appropriate staffing level $n$ such that the delay probability is $\varepsilon$ (i.e., $\pr(Q_\lambda(t) > n) \approx \varepsilon$),
we must characterize the distribution of $Q_\lambda(t)$. 
For large values of $\lambda$, its scaled version $\widehat{Q}_\lambda(t)$ can be approximated in distribution by $\widehat{Q}_\infty(t)$, as deduced from Theorem~\ref{thm:QQ}. 
We now show that under a mild condition, $\widehat{Q}_\infty(t)$ follows a normal distribution. 

\begin{corollary}\label{coro:Gauss}
    Under the assumptions of Theorem~\ref{thm:QQ}, if $X_\dagger$ is a constant or an independent normal random variable, then $\widehat{Q}_\infty(t)$ is a Gaussian process. 
\end{corollary}

We know from Equation~\eqref{eqn:hatQ}, Theorem~\ref{thm:QQ}, and Corollary~\ref{coro:Gauss} that 
\begin{equation}\label{eq:Q-approx-dist-0}
\begin{aligned}
Q_\lambda(t) - \lambda \int_0^t \bar{F}(t-z)\dd{z} ={} \lambda^{\frac{\alpha+1}{2}} \widehat{Q}_\lambda(t) 
\stackrel{d}{\approx} {}  \lambda^{\frac{\alpha+1}{2}} \widehat{Q}_\infty(t) 
\stackrel{d}{=}  \mathsf{Normal}\left(0, \lambda^{\alpha+1}  \Var[\widehat{Q}_\infty(t)]\right),
\end{aligned}
\end{equation}
where $\stackrel{d}{\approx}$ (resp., $\stackrel{d}{=}$) denotes approximate equality  (resp., equality) in distribution.  
To obtain a time-independent staffing rule, we further consider a steady-state approximation by sending $t\to\infty$:
\begin{equation}\label{eq:Q-approx-dist}
    Q_\lambda(t) \stackrel{d}{\approx} \lim_{t\to\infty}  \biggl\{\lambda  \int_0^t \bar{F}(t-z)\dd{z} + \mathsf{Normal}\left(0, \lambda^{\alpha+1} \Var[\widehat{Q}_\infty(t)]\right) \biggr\}.
\end{equation}

Because of the mutual independence among $U(t)$, $\widetilde{B}(z)$, and $\xi(t)$, 
we have 
\begin{equation}\label{eq:Var-decomp}
\Var[\widehat{Q}_\infty(t)] =  \underbrace{\Var\biggl[\int_0^t \bar{F}(t-z)U(z)\dd{z}\biggr]}_{\coloneqq \, V_1(t) }  +   
\biggl\{\underbrace{\Var\biggl[ \int_0^t \bar{F}(t-z) \dd{\widetilde{B}(z)} \biggr] }_{\coloneqq \, V_2(t)}
+ \underbrace{\vphantom{\bigg[1\bigg]} \Var[\xi(t)]}_{\coloneqq \, V_3(t)} \biggr\} \mathbb{I}\{\alpha = 0\}. 
\end{equation}
In general, $V_1(t)$ does not admit a simplified expression and must be evaluated numerically. 
However, we will discuss a special case in Example~\ref{ex:exp} later. 
Regarding $V_2(t)$ and $V_3(t)$, we have 
\begin{equation}\label{eq:V2+V3}
V_2(t) + V_3(t)    =  \int_0^t \bar{F}^2(t-z) \dd{z} + \int_0^t F(t-z)\bar{F}(t-z)  \dd{z}   = \int_0^t \bar{F}(t-z) \dd{z},
\end{equation}
where the first step is derived from  the It\^{o} isometry and Equation~\eqref{eq:xi-cov}, while the second step holds because $\bar F = 1 - F$ by definition. 
In addition, since $F$ is the CDF of the service time $S$, it is straightforward to show that 
\begin{equation}\label{eq:mean-service}
    \lim_{t\rightarrow\infty} \int_0^t \bar{F}(t-z)\dd{z} = \E[S] = \frac{1}{\mu},
\end{equation}
where $\mu$ denotes the service rate. 
It then follows from %
Equations~\eqref{eq:Q-approx-dist}--\eqref{eq:mean-service} that 
\begin{equation}\label{eqn:eq140}
Q_\lambda(t) \stackrel{d}{\approx}  
\mathsf{Normal} \left( \frac{\lambda}{\mu}, \, \lambda^{\alpha+1} \biggl(  V_1(\infty) + \frac{1}{\mu}   \mathbb{I}\{\alpha = 0\} \biggr) \right),
\end{equation}
for large values of $\lambda$ and $t$, where $V_1(\infty) = \lim_{t\to\infty} V_1(t)$. 

This approximation offers an immediate staffing rule for large-scale service systems. 
To satisfy $\pr(Q_\lambda(t) > n) \approx \varepsilon$, we can configure the number of servers to be the $(1-\varepsilon)$-quantile of the normal distribution outlined in the approximation~\eqref{eqn:eq140}.
More specifically, let $\Phi$ represent the CDF of the standard normal distribution and define $\beta = \Phi^{-1}(1-\varepsilon)$. 
We propose a staffing rule as:
\begin{equation}\label{eq:basic-alpha}
n^\ast_{\mathsf{Basic}}  =   \frac{\lambda}{\mu} + \beta   \lambda^{\frac{\alpha +1}{2}} \sqrt{  V_1(\infty) +  \frac{1}{\mu}   \mathbb{I}\{\alpha = 0\} }.
\end{equation}
We  refer to this staffing rule as the basic alpha safety rule.

The safety level in the basic alpha safety rule \eqref{eq:basic-alpha} varies across the three scenarios in Table~\ref{tab:arrivals}. In the intensity-dominated scenario ($\alpha\in(0,1)$), the safety level is $\beta \lambda^{(\alpha +1)/2} \sqrt{ V_1(\infty) }$, which is of higher order than $\lambda^{1/2}$ and thus significantly exceeds the square-root safety rule. For instance, with $\alpha=0.5$ and $\lambda/\mu= 100$, the basic alpha safety level is $37$ versus $17$ under \cite{Whitt92}'s square-root rule (see Table~\ref{tab:staffing-levels}).
In the balanced scenario ($\alpha=0$, $\sigma>0$), the safety level is $\beta \lambda^{1/2} \sqrt{ V_1(\infty) + 1/\mu }$---a square-root rule, but with a larger constant due to over-dispersion.
In the Poisson-dominated scenario ($\alpha=0$, $\sigma= 0$), $V_1(\infty)=0$, recovering the classical square-root safety rule $\beta \sqrt{\lambda /\mu}$ of \cite{Whitt92}.
Table~\ref{tab:staffing} summarizes these results, illustrating that as $\alpha$ varies from 0 to 1, the basic alpha safety rule bridges the square-root rule (Example \ref{example:sr-staffing}) and the linear rule (Example \ref{example:linear-staffing}).

\begin{table}[t]
\TABLE{Trichotomy of Over-dispersion's Impact on Staffing. \label{tab:staffing}}
{\begin{tabular}{@{}ccccccc@{}}
\toprule
  && \multicolumn{1}{c}{\begin{tabular}[c]{@{}c@{}}Intensity-dominated\\ ($0 < \alpha < 1$; $\sigma>0$)\end{tabular}} && \multicolumn{1}{c}{\begin{tabular}[c]{@{}c@{}}Balanced\\ ($\alpha = 0$; $\sigma>0$)\end{tabular}} && \multicolumn{1}{c}{\begin{tabular}[c]{@{}c@{}}Poisson-dominated\\ ($\alpha=0$; $\sigma = 0$, $X(0) = \lambda$)\end{tabular}} \\ 
  \midrule
~$\widehat{Q}_\infty(t)$   &&  $\displaystyle\int_0^t \bar{F}(t-z) U(z) \dd{z}$    &&                  $\displaystyle \int_0^t \bar{F}(t-z) \dd{\widehat{A}_\infty(z)} + \xi(t)$    &&        $\displaystyle\int_0^t \bar{F}(t-z) \dd{\widetilde{B}(z)} + \xi(t) \vphantom{\bigg[1\bigg]}$ \\
~Safety Level && $\mathcal{O}(\lambda^{(\alpha+1)/2})$ && $\mathcal{O}(\lambda^{1/2})$  && $\mathcal{O}(\lambda^{1/2}) \vphantom{\bigg[1\bigg]}$ \\
\bottomrule
\end{tabular}
}
{}
\end{table}

The staffing-level trichotomy presented by \cite{BassambooRandhawaZeevi10} shares conceptual similarities with Table~\ref{tab:staffing}, but differs in a fundamental way.
Their analysis attributes over-dispersion to parameter uncertainty: arrivals follow a Poisson process with a constant but \emph{unknown} rate $\lambda$ that must be inferred from data. Their decision-maker must therefore account for uncertainty about $\lambda$ when determining staffing levels. Their three scenarios are distinguished by the relative contributions of parameter uncertainty versus conditional Poisson randomness to the overall arrival variation.

This interpretation suggests that over-dispersion diminishes as estimation of $\lambda$ improves with more data, potentially being eliminated entirely. In contrast, our DSPP model treats the arrival rate as a stochastic process, which is inherently random.  This randomness persists regardless of data volume, unlike parameter uncertainty. Consequently, our model implies that even with abundant arrival data---sufficient to eliminate parameter uncertainty---significantly higher staffing levels than those suggested by square-root safety rules may still be necessary to ensure quality service.

We conclude this subsection by examining the computation of $V_1(\infty)$.
It is easy to show that
\begin{equation}\label{eqn:V_1-expand}
V_1(t) = \Var\biggl[ \int_0^t \bar{F}(t-z)U(z)\dd{z} \biggr] 
 =  2 \int_0^t \int_0^s \bar{F}(t-s) \bar{F}(t-v) \Cov[U(s), U(v)] \dd{v} \dd{s}. 
\end{equation}
We can approximate $V_1(\infty)$ using $V_1(t)$ for sufficiently large $t$ and evaluate it through numerical integration, provided we can compute $\Cov[U(s), U(v)]$. 
This covariance's expression depends on the OU process's initial condition and has a closed form when $U(0)$ is constant or follows the stationary distribution \citep[Section~3.3]{Glasserman03}. 
Further, with exponential service times, $V_1(\infty)$ has a closed form, simplifying the basic alpha safety rule.

\smallskip
\begin{example}[Basic Alpha Safety Rule with Exponential Services]\label{ex:exp} 
Assume $U(0)$ is either constant or follows the OU process's stationary distribution, and service times are exponentially distributed with rate $\mu$ (i.e., $\bar F(t) = e^{-\mu t}$). 
Elementary calculations yield
$V_1(\infty) = \frac{\sigma^2}{2\kappa\mu(\mu+\kappa)}$, so the staffing formula \eqref{eq:basic-alpha} simplifies to:
\begin{align*}
    n^\ast_{\mathsf{Basic}}  =   \frac{\lambda}{\mu} + \beta   \lambda^{\frac{\alpha +1}{2}} \sqrt{  \frac{\sigma^2}{2\kappa\mu(\mu+\kappa)} +  \frac{1}{\mu}   \mathbb{I}\{\alpha = 0\} }.
\end{align*}
\end{example}

\subsection{Refined Alpha Safety Rule}\label{sec:staff_refine}

While the infinite-server approximation provides useful managerial insights, it tends to under-staff because an infinite-server system has fewer customers in system than its finite-server counterpart, leading to a lower staffing recommendation (see Section~\ref{sec:synthetic} for numerical illustration).

The seminal work of \cite{HalfinWhitt81} addressed this for Poisson arrivals by analyzing finite-server systems in a regime where the number of servers scales as $\lambda/\mu + c\sqrt{\lambda}$, yielding a refined square-root safety rule. A key insight is that the finite-server analysis preserves the structure of the infinite-server rule, only refining the leading constant. Inspired by this, we propose a refined alpha safety rule that retains the structure of the basic rule \eqref{eq:basic-alpha} but adjusts the coefficient:
\begin{equation}\label{eq:general-alpha}
n_\lambda = \frac{\lambda}{\mu} + \delta \lambda^{\frac{\alpha +1}{2}},
\end{equation}
where $\delta>0$ is a design parameter calibrated to achieve a target delay probability $\varepsilon$.

Let $Q^{(n)}_\lambda(t)$ denote the number of customers at time $t$ in an $n$-server system with DSPP arrivals. Analogous to the Halfin--Whitt analysis, consider the scaled process
\begin{equation}\label{eq:scaled-num-in-sys-finite-server}
\widehat{Q}_\lambda^{(n_\lambda)}(t) \coloneqq
\lambda^{-\frac{\alpha+1}{2}}\bigl[Q_\lambda^{(n_\lambda)}(t) - n_\lambda \bigr].
\end{equation}
We anticipate that $\widehat{Q}_\lambda^{(n_\lambda)}$ converges in distribution to a non-degenerate limit with a stationary distribution as $\lambda\to\infty$. This can be established by extending the approach of \cite{HalfinWhitt81} and \cite{SunLiu21} (who treated the special case $\alpha=0$); a complete proof is deferred to future research.
We assume this convergence holds.
Based on this, we propose a simulation-based method (Algorithm~\ref{alg:refined-alpha}) to determine $\delta$ in staffing formula \eqref{eq:general-alpha} for achieving a delay probability of $\varepsilon$.

\begin{algorithm}[t]
\caption{Solving for $\delta^*$ via Simulation and Stochastic Approximation}\label{alg:refined-alpha}
\begin{algorithmic}[1]
\State \textbf{Input:} Model parameters  ($\alpha$, $\kappa$, $\sigma$, and $\mu$) and algorithm parameters ($\lambda$, $t$, $m$, and $a_i$)
\State Initialize $\delta_0$ as specified by Equation~\eqref{eq:delta-init}
\While{$i=0,1,\ldots$}
\State Set $n_{\lambda,i} = \lambda / \mu + \delta_i \lambda^{(\alpha+1)/2}$
\State Simulate $m$ replications of the arrival process: $A_{\lambda,i,j}(t)$, $j=1,\ldots,m$ 
\State Simulate the number-in-system $\widehat{Q}_{\lambda,i,j}^{(n_{\lambda,i})}(t)$ in an $n_{\lambda,i}$-server system for each replication $j$
\State Calculate $M(\delta_i)$ as defined in \eqref{eq:PASTA-estimate}: 
$M(\delta_i) \coloneqq \frac{1}{m}\sum_{j=1}^m \mathbb{I}\bigl\{\widehat{Q}_{\lambda,i,j}^{(n_{\lambda,i})}(t) > 0 \bigr\}$
\State Update: $\delta_{i+1} =  \delta_i - a_i [M(\delta_i) - \varepsilon] $
\EndWhile
\end{algorithmic}
\end{algorithm}

Let $\widehat{Q}^*(t)$ denote the limiting process to which  $\widehat{Q}_\lambda^{(n_\lambda)}(t)$ converges weakly as $\lambda\to\infty$. 
The steady-state delay probability in the heavy-traffic limit is $\pr(\widehat{Q}^*(\infty) > 0)$, which depends on $\delta$ and approximates $\pr(\widehat{Q}_\lambda^{(n_\lambda)}(t) > 0 )$ (i.e., the delay probability of a finite-server system) for large values of $\lambda$ and $t$. 
As a result, we seek a value of $\delta$ such that 
$\pr\bigl(\widehat{Q}^*(\infty) > 0\bigr)  = \varepsilon$.    
Let $\delta^*$ solve this equation. We call the following staffing rule the refined alpha safety rule: 
\begin{equation}\label{eq:RefinedStaffing}
n^\ast_{\mathsf{Refined}} = \frac{\lambda}{\mu} + \delta^* \lambda^{\frac{\alpha +1}{2}}.    
\end{equation}

Despite not having an explicit characterization of the distribution of $\widehat{Q}^*(\infty)$, 
we can use simulation and apply the stochastic approximation method \citep[Chapter~4]{Spall03} to approximately solve for $\delta^*$. 
Specifically, 
we begin with an initial estimate, denoted as $\delta_0$.
A good choice for $\delta_0$ is the coefficient of the $\lambda^{(\alpha+1)/2}$ term in the basic alpha safety rule \eqref{eq:basic-alpha}, that is, 
\begin{equation}\label{eq:delta-init}
\delta_0 = \beta  \sqrt{  V_1(\infty) +  \frac{1}{\mu}   \mathbb{I}\{\alpha = 0\} }.    
\end{equation}

Then, we update the estimate via the recursion 
\begin{equation}\label{eq:SA-recursion}
    \delta_{i+1} = \delta_i - a_i [M(\delta_i) - \varepsilon],
\end{equation}
where $a_i>0$ is the step size that satisfies $\sum_{i=0}^\infty a_i = \infty$ and $\sum_{i=0}^\infty a_i^2 < \infty$, and 
$M(\delta)$ is a random variable that satisfies $\E[M(\delta)] \approx  \pr(\widehat{Q}^*(\infty) > 0)$. 
A common choice for the step size is $a_i = b/(i+c)^d$ for some constants $b>0$ and $c>0$, and $1/2< d \leq 1$. 

The most straightforward choice for $M(\delta)$ is the indicator function $\mathbb{I}\{\widehat{Q}_\lambda^{(n_\lambda)}(t) > 0 \}$ for some sufficiently large values of $\lambda$ and $t$. (Note that $n_\lambda$ depends on $\delta$.)
However, due to the strong intertemporal dependencies introduced by the doubly stochastic Poisson structure, this indicator function tends to exhibit a significant variance. 
This would make the recursion in Equation~\eqref{eq:SA-recursion} slow to converge.
Therefore, we propose an alternative choice for $M(\delta)$ to reduce variance, based on multiple replications of simulating the number-in-system process:
\begin{equation}\label{eq:PASTA-estimate}
  M(\delta) \coloneqq \frac{1}{m}\sum_{j=1}^m \mathbb{I}\bigl\{\widehat{Q}_{\lambda,j}^{(n_\lambda)}(t) > 0 \bigr\},  
\end{equation}
where the subscript $j$ denotes the $j$-th replication; see
\cite{GieseckeKimZhu11} for methods for simulating our DSPP model. 

In our experiments (Sections~\ref{sec:synthetic}--\ref{sec:casestudy}), we stop Algorithm~1 when $|M(\delta_i)-\varepsilon|\leq 0.01$, which typically happens within dozens of iterations, with a run time of approximately 10 minutes on a standard personal computer. 
When available, real arrival data should replace simulation data in the step~5.

\section{Model Comparison}\label{sec:model-comparison}

In this section, we compare our DSPP model with several existing arrival models that both account for over-dispersion and yield analytical staffing rules. Our model serves as a parsimonious generalization of these models. We analyze which arrival features these models capture and their corresponding safety staffing rules for delivering target service quality. Furthermore, we develop a MLE method based on heavy-traffic approximation. This approach dramatically reduces the computational burden of standard MLE for DSPP models while enabling data-driven, in-sample model selection through information criteria.

We analyze the following five models,\footnote{\cite{Avramidis04} and \cite{OreshkinReegnardLEcuyer16} developed models for over-dispersed arrivals that focus on the joint distribution of arrival counts across different time intervals. While these models are valuable for queueing system simulation, they do not yield analytical staffing rules, and thus are not included in our model comparison. However, we use one of these models in Appendix~\ref{sec:additional-exp} to evaluate the robustness of our DSPP model and its staffing rule under model misspecification.} denoted as $\mathcal{M}_1$ through $\mathcal{M}_5$:

\begin{itemize}[noitemsep]
    \item $\mathcal{M}_1$: Poisson process model.
    
    \item $\mathcal{M}_2$:     
    DSPP model with intensity $ \lambda G$, where $G$ is a unit-mean random variable with variance~$\sigma_G^2$. 

\item $\mathcal{M}_3$: DSPP model with intensity $\lambda + \lambda^{(\alpha+1)/2} Y$, where $\alpha\in[0, 1]$ and $Y$ is a zero-mean random variable with variance $\sigma_Y^2$.

\item $\mathcal{M}_4$: DSPP model with intensity following a CIR process (Equation~\eqref{eq:CIR} with $\alpha=0$).  

\item $\mathcal{M}_5$:  Our proposed DSPP model. 

\end{itemize}

\subsection{Model Features}

Unlike the standard Poisson process model ($\mathcal{M}_1$), the other four models capture arrival over-dispersion through a random intensity component in their doubly stochastic Poisson structure. Among these, $\mathcal{M}_2$ represents a simple approach where a single random variable $G$ models the arrival rate uncertainty \citep{Whitt99,Avramidis04,BassambooRandhawaZeevi10}. 
As shown in Example~\ref{example:linear-staffing}, this model yields a linear relationship between the variance and mean of arrival counts within any time interval, corresponding to Taylor's law \eqref{eq:Taylor} with $\alpha=1$, and consequently produces a linear safety rule.

$\mathcal{M}_3$ extends $\mathcal{M}_2$ by introducing an additional parameter $\alpha\in[0,1]$. When $\alpha=1$, $\mathcal{M}_3$ reduces to $\mathcal{M}_2$, with their random components related by $G = Y+1$. 
With the presence of this parameter, 
like our proposed DSPP model ($\mathcal{M}_5$), $\mathcal{M}_3$ captures Taylor's law for dispersion scaling. 
Moreover, $\mathcal{M}_3$ can be understood as a special case of $\mathcal{M}_5$. When $\sigma=0$, the generalized CIR intensity~\eqref{eq:CIR} becomes $X(t) = \lambda + ( X(0) - \lambda) e^{-\kappa t}$. Setting $\kappa=0$ further simplifies this to $X(t) = X(0)$ for all $t\geq 0$. Thus, with $\kappa=\sigma=0$, $\mathcal{M}_5$'s intensity~\eqref{eq:CIR} becomes the static random variable $X(0)$. Finally, setting $X(0)=\lambda+\lambda^{(\alpha+1)/2}Y$ reduces $\mathcal{M}_5$ to $\mathcal{M}_3$.\footnote{Since Theorem~\ref{thm:FCLTofA} and its resulting distribution~\eqref{eq:A_inf-dist-pi} require strictly positive parameters ($\kappa,\sigma>0$), they cannot be directly applied to $\mathcal{M}_3$. We address this by deriving a separate asymptotic distribution for the arrival process under $\mathcal{M}_3$ (see Section~\ref{StaticDSPP} of the e-companion).}

\cite{Maman09}, \cite{BassambooRandhawaZeevi10}, and \cite{HuChanDong25} used $\mathcal{M}_3$ in their staffing analyses under arrival rate uncertainty, albeit with a different objective (e.g., cost-minimization in a newsvendor-type formulation for the latter two studies). Our approach, in contrast, focuses on achieving target quality of service through heavy-traffic queueing analysis. To ensure fair comparison between $\mathcal{M}_3$ and $\mathcal{M}_5$, we derive the safety rule for $\mathcal{M}_3$ using our heavy-traffic framework (see Appendix~\ref{StaticDSPP}). 
While the resulting safety rules share similar analytical forms, their actual staffing levels may differ when these models are fitted to arrival data, due to different parameter estimates.

$\mathcal{M}_4$, proposed by \cite{ZhangHongGlynn}, is a special case of $\mathcal{M}_5$ with $\alpha=0$. As shown in Table~\ref{tab:arrivals}, this model represents a ``balanced'' scenario where, in the heavy-traffic asymptotic regime, the impact of the intensity's randomness matches that of the conditional Poisson distribution on the arrival process's overall variability, bridging the gap between the Poisson process model and our DSPP model.  Table~\ref{tab:staffing} shows that while $\mathcal{M}_4$ yields a square-root safety rule similar to $\mathcal{M}_1$, it prescribes higher staffing levels to account for the additional randomness in the intensity; see also \cite{SunLiu21} for detailed heavy-traffic analysis.

Among the four models for arrival over-dispersion, $\mathcal{M}_2$ and $\mathcal{M}_3$ have a \emph{static} intensity, where randomness is independent of time $t$, while $\mathcal{M}_4$ and $\mathcal{M}_5$ have intensities that are stochastic processes, making them \emph{dynamic}.
While static intensity models benefit from simplicity, they have a critical shortcoming: they imply constant correlation between arrival counts across different time intervals, regardless of their temporal separation. 
This characteristic contradicts empirical evidence of temporal correlation decay \citep{Avramidis04,ibrahim2016modeling,OreshkinReegnardLEcuyer16}. In contrast, the two DSPP models with dynamic intensities capture this decay through the mean-reverting property of the (generalized) CIR process.\footnote{The covariance/correlation between arrival counts in two intervals $A(t_1, t_2]$ and $A(t_3, t_4]$ under our DSPP model can be derived analytically using the Laplace transform approach, following calculations similar to those in Section~\ref{ec-sec:laplace} of the e-companion, where we derive the arrival count variance in Theorem~\ref{theo:CoD}.} 

The significance of our DSPP model's ability to capture temporal correlation decay is reflected in two aspects. First, compared to $\mathcal{M}_3$, which also captures Taylor's law and yields an alpha safety rule, the fundamental difference in temporal correlation structure (constant versus decaying) leads to distinct parameter estimates when both models are fitted to arrival data. 
This results in different staffing levels, despite their safety rules sharing similar analytical forms.

Second, by modeling temporal correlation structure more accurately, our DSPP model is better suited for arrival forecasting \citep{ibrahim2016modeling}. Thus, our model could potentially support dynamic adjustments of staffing levels and workforce schedules \citep{ShenHuang08_Pois_SVD,GansShenZhouKorolevMcCordRistock15}. However, this application lies beyond the scope of this paper and remains for future research.

Table~\ref{tab:5models} summarizes the main features of these models and their corresponding safety rules. 
Figure~\ref{fig:model-relation} illustrates how these models can be transformed into one another.

\begin{table}[t]
	\TABLE{Features and Safety Rules of $\mathcal{M}_1$ through $\mathcal{M}_5$. \label{tab:5models}}
	{	
	\begin{tabular}
		[c]{lccccccccccc}%
		\toprule 
		~Model~  & ~\makecell{Model \\ Parameters}~ & ~\makecell{Doubly \\ Stochastic}~ & ~\makecell{Taylor's \\ Law}~ & 
		~\makecell{Temporal \\ Correlation Decay}~ & ~\makecell{Safety \\ Rule}~ \\
		\midrule 
~$\mathcal{M}_1$ & $\lambda$ & No & No  & No & Square-root\\
~$\mathcal{M}_2$  & $\lambda,\sigma_G$ & Yes & No & No & Linear  \\
~$\mathcal{M}_3$ & $\lambda,\alpha,\sigma_Y$ & Yes & Yes & No & Alpha\\
~$\mathcal{M}_4$  & $\lambda,\kappa,\sigma$ & Yes & No & Yes & Square-root \\
~$\mathcal{M}_5$ (this paper) & $\lambda,\alpha,\kappa,\sigma$ & Yes & Yes & Yes & Alpha \\
\bottomrule  
	\end{tabular}}
{A ``Yes'' for Taylor's law indicates the model captures this relationship for any $\alpha \in [0,1]$. Models may yield the same type of safety rule (e.g., $\mathcal{M}_1$ and $\mathcal{M}_3$ both produce square-root rules) but with different safety levels due to distinct parameter estimates from data fitting (see Table~\ref{tab:staffing-levels}).}
\end{table}

\begin{figure}[t]
	\begin{center}
		\caption{Nested Relationships Among $\mathcal{M}_1$ through $\mathcal{M}_5$. \label{fig:model-relation}}
		\begin{tikzpicture}[
scale=0.85,
transform shape,
			node distance=2cm,
			box/.style={draw, rounded corners, minimum width=2cm, minimum height=1cm},
			myarrow/.style={-{Triangle[length=2mm,width=2mm]}, thick}
			]
			\node at (0,3) {};
			
			\node[box] (M5) at (-4,0) {$\mathcal{M}_5$};
			
			\node[box] (M3) at (0,2) {$\mathcal{M}_3$};
			\node[box] (M4) at (2,0) {$\mathcal{M}_4$};
			
			\node[box] (M2) at (4,2) {$\mathcal{M}_2$};
			
			\node[box] (M1) at (8,0) {$\mathcal{M}_1$};
			
			\draw[myarrow] (M5) -- (M3);
			\draw[myarrow] (M5) -- (M4);
			\draw[myarrow] (M3) -- (M2);
			\draw[myarrow] (M2) -- (M1);
            \draw[myarrow] (M4) -- (M1);
			
			\node[above] at (-1,0) {$\alpha=0$};
			\node[above] at (-3.5,0.9) {$\begin{array}{c}\kappa=\sigma=0 \\ X(0)=\lambda + \lambda^{\frac{\alpha+1}{2}} Y \end{array}$};
			\node[above] at (2,2) {$\alpha=1$};
			\node[below] at (2,2) {$G=Y+1$};
			\node[above] at (6.4,1.1) {$\sigma_G = 0$};
            \node[above] at (5,0) {$\sigma=0$};
            \node[below] at (5,0) {$X(0)=\lambda$};
			
		\end{tikzpicture}
	\end{center}
\end{figure}

\subsection{Model Estimation}\label{sec:model-estimation}

The five models form a hierarchy of increasing complexity, from $\mathcal{M}_1$ to $\mathcal{M}_5$, both in terms of parameter count and model structure (progressing from static to dynamic intensities). For models with doubly stochastic Poisson structure, significant estimation challenges arise due to unknown parameters in the \emph{unobservable} intensity. The challenge is most acute for $\mathcal{M}_4$ and $\mathcal{M}_5$, where the intensity follows a stochastic process.

To illustrate these challenges, consider observing arrival counts in sequential time periods of length~$\Delta$, denoted by $N_{\Delta,1}, \ldots, N_{\Delta,k}$. The likelihood of these observations under our DSPP model is
\begin{equation}\label{eq:likelihood-exact}
\begin{aligned}
& \pr\Big(A_\lambda(\Delta) = N_{\Delta,1}, A_\lambda(2\Delta) - A_\lambda(\Delta) = N_{\Delta,2},\ldots, A_\lambda(k\Delta) - A_\lambda((k-1)\Delta) = N_{\Delta,k}\Big) \\ 
={}&  \E\left[\prod_{i = 1}^k \pr\left(A_\lambda(i\Delta) - A_\lambda((i-1)\Delta) = N_{i, \Delta}\,\bigg|\, \int_{(i-1)\Delta}^{i\Delta} X_\lambda(s)\dd{s} \right)\right],
\end{aligned}
\end{equation}
where the expectation is taken with respect to the joint distribution of the integrated intensities $\int_{(i-1)\Delta}^{i\Delta}X(s)\dd{s}$. While the conditional likelihood has an analytical form (being a Poisson distribution), the full likelihood requires $k$-dimensional integration lacking an analytical form. 
Since MLE requires numerical optimization over the unknown parameters, this integration must be performed in each iteration, rendering the computation prohibitive as $k$ increases.
An expectation-maximization algorithm could be used, but its expectation steps require sampling-based computation through Markov chain Monte Carlo (MCMC). 
This approach remains computationally intensive, and its numerical performance is highly sensitive to MCMC tuning parameters \citep[Chapter~6]{Duffie11}.

To address these challenges, we introduce an efficient, MLE method for our DSPP model, leveraging the heavy-traffic approximation developed in Section~\ref{themodel}. Since the other four models are special cases of our model, this method extends naturally to them. 

From Theorem~\ref{thm:FCLTofA}, when $\lambda$ is large, the arrival process $A_\lambda(t)$ under our DSPP model is approximately distributed as a Gaussian process: when $\lambda$ is large, 
\begin{align*}
A_\lambda(t) \stackrel{d}{\approx}  {}& 
\lambda^{\frac{\alpha+1}{2}} \widehat{A}_\infty(t) + \lambda t
= \lambda^{\frac{\alpha+1}{2}}\left(\int_0^t U(s) \dd{s} + \widetilde{B}(t) \mathbb{I}\{\alpha = 0\} \right)+ \lambda t,
\end{align*}
where $U$ is the OU process \eqref{OU} and $\widetilde{B}$ is a standard Brownian motion independent of $U$. This heavy-traffic approximation implies that the arrival counts $N_{\Delta,1}, \ldots, N_{\Delta,k}$ are approximately multivariate normal, denoted by $\mathsf{MVNormal}(\mathbf{m}, \mathbf{\Sigma})$, 
where $\mathbf{m} = (\lambda \Delta) \mathbf{1}$ is the mean vector with $\mathbf{1}$ a vector of ones, and $\mathbf{\Sigma}$ is the covariance matrix with elements $\Sigma_{i,\ell}$ given by
\begin{align*}
\Sigma_{i,\ell}= \left\{
\begin{array}{ll}
     \displaystyle\lambda \Delta +  \frac{\sigma^2 \lambda^{\alpha+1} \Delta}{\kappa^2} \left( 1- \frac{1-e^{-\kappa \Delta}}{\kappa \Delta} \right),&  \mbox{if } i = \ell, \\[2ex] 
    \displaystyle \lambda^{\alpha+1} \frac{\sigma^2}{2\kappa^3} \left[ e^{-\kappa (\ell-1) \Delta} - e^{-\kappa \ell \Delta} \right] \left[ e^{\kappa i \Delta} - e^{\kappa (i-1) \Delta} + e^{-\kappa i \Delta} - e^{-\kappa (i-1) \Delta} \right], & \mbox{if } i < \ell, \\[2ex] 
    \displaystyle \lambda^{\alpha+1} \frac{\sigma^2}{2\kappa^3} \left[ e^{-\kappa (i-1) \Delta} - e^{-\kappa i \Delta} \right] \left[ e^{\kappa \ell \Delta} - e^{\kappa (\ell-1) \Delta} + e^{-\kappa \ell \Delta} - e^{-\kappa (\ell-1) \Delta} \right], & \mbox{if } i > \ell. 
\end{array}
\right.
\end{align*}

These analytical expressions for the mean vector and covariance matrix are derived using properties of OU processes and Brownian motions, following calculations similar to those for the marginal distribution of $\widehat{A}_\infty(t)$ in Equation~\eqref{eq:A_inf-dist-pi}; see \citet[Chapter~3.3]{Glasserman03}.\footnote{In Section~\ref{sec:casestudy}, we extend our DSPP model to accommodate non-stationary arrivals and adapt the multivariate normal approximation accordingly in Appendix~\ref{sec:MLE}. The application of this approach to estimating models $\mathcal{M}_2$ and $\mathcal{M}_3$ (the DSPP models with static intensities) is detailed in Appendix~\ref{StaticDSPP}.}

In practice, we may observe $m$ independent sample paths of the arrival process, yielding count data $\{\mathbf{N}_{\Delta,j}\coloneqq (N_{\Delta,1,j}, \ldots, N_{\Delta, k, j})^\intercal: j=1,\ldots,m\}$. Under the multivariate normal approximation, the estimation reduces to fitting parameters for $m$ independent samples of a $k$-dimensional normal vector. The log-likelihood function of the unknown parameters $\Theta = (\lambda, \alpha, \kappa, \sigma)$,  summing individual log-likelihoods of $\mathbf{N}_{\Delta,j}$, is approximated using the multivariate normal density $\mathsf{MVNormal}(\mathbf{m},\mathbf{\Sigma})$:
\begin{align}
    \log\mathcal{L}(\Theta) \coloneqq {}& \sum_{j=1}^m \log \pr\Big(A_\lambda(i\Delta) - A_\lambda((i-1)\Delta) = N_{\Delta,i,j}, \, i = 1,\ldots,k\Big)   \nonumber \\ 
    \approx{}&  \sum_{j=1}^m \log\left( (2\pi)^{-k/2} \det(\mathbf{\Sigma})^{-1/2} \exp\left(-\frac{1}{2} (\mathbf{N}_{\Delta,j} - \mathbf{m})^\intercal
     \mathbf{\Sigma}^{-1} (\mathbf{N}_{\Delta,j} - \mathbf{m})  \right) \right) \nonumber \\ 
     ={}  & - \frac{k m}{2} \log(2\pi) - \frac{m}{2} \det(\mathbf{\Sigma}) - \frac{1}{2}\sum_{j=1}^m (\mathbf{N}_{\Delta,j} - \mathbf{m})^\intercal
     \mathbf{\Sigma}^{-1} (\mathbf{N}_{\Delta,j} - \mathbf{m}), \label{eq:log-likelihood-approx}
\end{align}
where $\det(\mathbf{\Sigma})$ is the determinant of $\mathbf{\Sigma}$.  
This approximation \eqref{eq:log-likelihood-approx} enables efficient MLE by avoiding the high-dimensional numerical integration in Equation~\eqref{eq:likelihood-exact}.

\subsection{Model Selection}\label{sec:model-selection}

The key features in Table~\ref{tab:5models} provide qualitative guidelines for model selection in practice. When analyzing arrival data, one can first examine whether the data shows significant over-dispersion, dispersion scaling consistent with Taylor's law, or temporal correlation decay. These empirical features offer managers a quick screening of suitable models.

Model selection can be further refined through comprehensive data analysis using both training and test datasets: the former for model estimation and the latter for evaluation. The selection can be based on either in-sample or out-of-sample performance.

\subsubsection{In-sample Performance}
The likelihood approximation in Section~\ref{sec:model-estimation} facilitates model selection through the Akaike information criterion (AIC) and Bayesian information criterion (BIC):
\begin{gather*}
  \textsf{AIC} = 2 q - 2 \log{}(\widehat{\mathcal{L}}) \quad\mbox{ and } \quad 
  \textsf{BIC} = q \log(m) - 2 \log{} (\widehat{\mathcal{L}}),
\end{gather*}
where $q$ is the number of unknown parameters, $m$ is the number of independent sample paths, and $\widehat{\mathcal{L}} \coloneqq \max_\Theta \mathcal{L}(\Theta)$ is the maximized likelihood value. 

Both criteria balance model fit (captured by $\widehat{\mathcal{L}}$) against model complexity (measured by $q$), with BIC imposing a stronger penalty that increases with sample size $m$. 
Lower values of AIC and BIC indicate better models, with differences exceeding 10 providing strong evidence favoring the model with the smaller value \citep{Raftery95,BurnhamAnderson02}.

\subsubsection{Out-of-sample Performance}

In-sample performance measures how well an arrival model fits the data. However, for managers, a model's value lies in its ability to predict system performance using its derived staffing rules. Since we focus on delay probability as the quality-of-service metric, we evaluate out-of-sample performance by computing the realized delay probability of the staffing level derived from the fitted model, as demonstrated in our numerical experiments.

\section{Numerical Experiments with Stationary Arrivals}\label{sec:synthetic}

In this section, we validate the proposed safety rules through numerical experiments using our DSPP model as ground truth. 
We conduct experiments in two settings: infinite-server systems for the basic alpha safety rule and finite-server systems for the refined alpha safety rule. 
Since the true model is known, we focus on out-of-sample performance.
We assess convergence of quality of service to theoretical values under increasing arrival rates, reflecting the heavy-traffic regime where these rules are derived.
We also compare both alpha safety rules against four alternative staffing rules from different arrival models (see Table~\ref{tab:5models}).

\subsection{General Setup}\label{sec:exp-setup}

We examine quality of service targets $\varepsilon = 0.05$ and $0.15$, with corresponding standard normal quantiles $\beta=\Phi^{-1}(1-\varepsilon) = 1.64$ and $1.04$. 
For the generalized CIR intensity~\eqref{eq:CIR} of our DSPP model, we set $\alpha=0.5$, $\kappa=0.1$, and $\sigma=0.5$. We set service time to be log-normal distribution with mean $1/6$ and standard deviation $1/6$. 
This choice of log-normal distribution is supported by \cite{Brown} and \cite{ShenBrown06}, who show that call center service times typically follow this distribution with standard deviation similar to mean. Additional numerical studies using gamma-distributed service times yield consistent results, presented in Appendix~\ref{app:sec:gamma}.

To evaluate the staffing rules' performance in smaller systems and demonstrate how dispersion scaling affects staffing and service quality as systems grow, we adopt a data-driven approach using a small-to-large system progression. This evaluation is particularly relevant since all staffing rules under $\mathcal{M}_1$ through $\mathcal{M}_5$ are derived through heavy-traffic asymptotic analysis.

First, we generate training data from $\mathcal{M}_5$ with $\lambda=100$ customers per hour, comprising $m=1000$ independent sample paths with parameters $(\alpha,\kappa,\sigma)$ as specified earlier. 
We fit models $\mathcal{M}_1$ through $\mathcal{M}_4$ to this small-system data using the estimation method from Section~\ref{sec:model-estimation}. 
Using the estimated parameters (e.g., $\sigma_G$ for $\mathcal{M}_2$ and $(\kappa,\sigma)$ for $\mathcal{M}_4$), we compute staffing levels under the four models for larger systems with $\lambda=150, 600,$ and $2400$.\footnote{For both basic and refined alpha safety rules, staffing levels are computed using the true model parameters of $\mathcal{M}_5$. All safety rules in Table~\ref{tab:5models} use the same $\beta$ value, except that the refined alpha safety rule is computed via Algorithm~\ref{alg:refined-alpha} with tuning parameters $\lambda=100$, $t=24$, $m=100$, and $a_i=20/(i+20)$.} 
Finally, we generate additional sample paths from our DSPP model with each $\lambda$ as test data to evaluate the performance of different staffing rules.

\begin{table}[t]
	\TABLE{Staffing Levels Under Different Fitted Models. \label{tab:staffing-levels}}
	{	
	\begin{tabular}
		[c]{l@{\extracolsep{30pt}}
        c@{\extracolsep{25pt}}
        c@{\extracolsep{25pt}}
        c@{\extracolsep{25pt}}
        c@{\extracolsep{25pt}}
        c@{\extracolsep{25pt}}
        c@{\extracolsep{25pt}}}%
		\toprule 
		& \multicolumn{3}{c}{$\varepsilon=0.05$} &  \multicolumn{3}{c}{$\varepsilon=0.15$} \\
		\cmidrule(l){2-4} \cmidrule(l){5-7} 
		Model  &
$\lambda=150$ & $\lambda=600$ & $\lambda=2400$ &  $\lambda=150$ & $\lambda=600$ & $\lambda=2400$ \\
\midrule 
$\mathcal{M}_1$ (Square-root) & $34$ & $117$   &  $433$  & $31$ &  $111$  &  $421$ \\

$\mathcal{M}_2$ (Linear)   &  $36$ & $141$ &    $564$ &  $32$ & $126$ &  $504$  \\

$\mathcal{M}_3$ (Alpha)  & $35$ & $125$ &  $468$ &  $31$ & $116$ &  $443$ \\

$\mathcal{M}_4$ (Square-root) &  $34$ & $118$   &  $435$  &   $31$ & $112$  &   $423$ \\

$\mathcal{M}_5$ (Basic alpha) & $38$ & $137$ & $504$ &  $34$ & $124$ &  $466$ \\

$\mathcal{M}_5$ (Refined alpha) & $42$ & $147$ &  $532$  & $37$ & $134$ &  $496$ \\
		\bottomrule  
	\end{tabular}}
{Parameters: $\alpha=0.5$, $\kappa=0.1$, $\sigma=0.5$. 
 Service time distribution: Log-normal with mean $1/6$ and standard deviation $1/6$. 
}
\end{table}

Table~\ref{tab:staffing-levels} shows the staffing levels determined by each of these six rules. For sufficiently large $\lambda$, the levels increase in the order of $\mathcal{M}_1$ (square-root), $\mathcal{M}_4$ (square-root), $\mathcal{M}_3$ (alpha), $\mathcal{M}_5$ (basic alpha), $\mathcal{M}_5$ (refined alpha), and $\mathcal{M}_2$ (linear). 
For example, when $\lambda=2400$ and $\varepsilon = 0.05$, the staffing levels of these six rules (in the same order) are $433$, $435$, $468$, $504$, $532$, and $564$, respectively.
This ordering aligns with their theoretical safety levels: $\mathcal{O}(\lambda^{1/2})$ for $\mathcal{M}_1$ and $\mathcal{M}_4$, $\mathcal{O}(\lambda^{(\alpha+1)/2})$ with $\alpha=0.5$ for $\mathcal{M}_3$ and $\mathcal{M}_5$, and $\mathcal{O}(\lambda)$ for $\mathcal{M}_2$. 

The square-root safety rule under $\mathcal{M}_4$ prescribes higher staffing than that under $\mathcal{M}_1$ because it incorporates intensity randomness through a CIR process, while $\mathcal{M}_1$ assumes constant intensity.
Among the three alpha safety rules, $\mathcal{M}_3$ yields lower levels than both rules under $\mathcal{M}_5$ because its static intensity misspecifies the temporal correlation structure, underestimating overall intensity uncertainty. The refined alpha safety rule suggests higher staffing than the basic alpha safety rule because the latter's infinite-server approximation underestimates delay probability in finite-server systems, as explained in Section~\ref{sec:staff_refine}.

\subsection{Infinite-server Systems}\label{sec:inf-ser-sym}

In this subsection, we numerically validate the basic alpha safety rule. Using the infinite-server queueing simulation detailed in Section~\ref{sec:IS-apprx}, we estimate the steady-state probability of $\{Q_\lambda(t) > n\}$, the event that the number-in-system exceeds the staffing level prescribed by the basic alpha safety rule and the rules under $\mathcal{M}_1$ through $\mathcal{M}_4$. We expect relation \eqref{eq:inf-serv-approx-obj} to hold: this probability should approach $\varepsilon$ as $\lambda$ increases when $n$ follows the basic alpha safety rule.

Given a mean arrival rate $\lambda$, we independently simulate $R=200$ sample paths of our DSPP model arrivals with parameters as detailed in Section~\ref{sec:exp-setup}. We feed each of these to the infinite-server system and record the number-in-system process, denoted by $Q_{\lambda,r}(t)$, for $r=1,\ldots,R$.
Since the system starts from an empty state---a condition significantly divergent from the steady state---we use the first 24 hours of simulation clock time as a warm-up period for each sample path. 
We estimate the steady-state probability $\pr( Q_\lambda(\infty) > n)$ based on an additional 24-hour period of the sample path, as follows:
\begin{equation*}%
\pr\big( Q_\lambda(t) > n \big)  \approx  \frac{1}{R} \sum_{r=1}^{R} \mathbb{I}\big\{ Q_{\lambda,r}(t) > n \big\}, \quad t\in[24, 48].
\end{equation*}

\begin{figure}[t]

\FIGURE{    \includegraphics[width=0.5\textwidth]{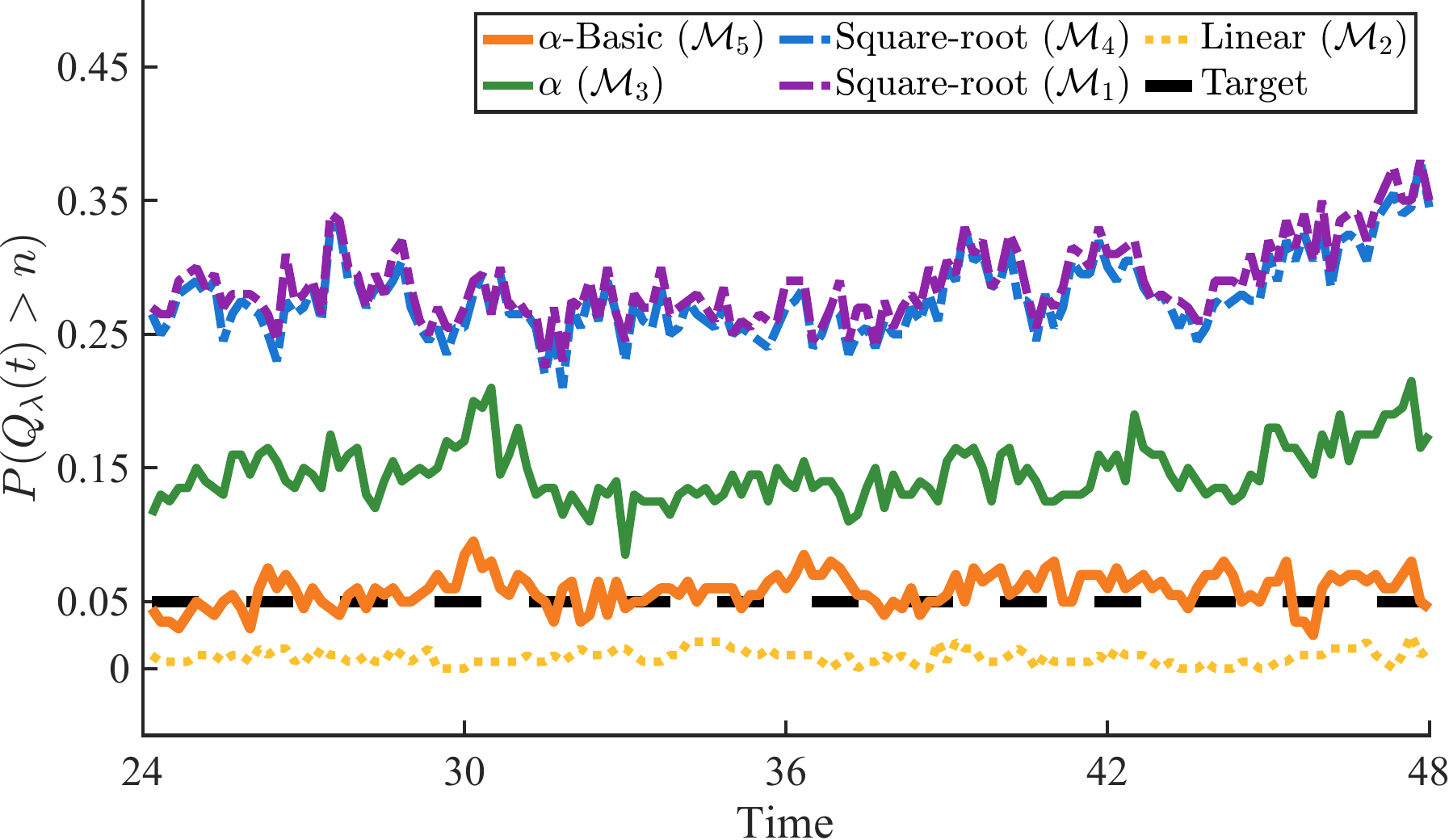}
	\includegraphics[width=0.5\textwidth]{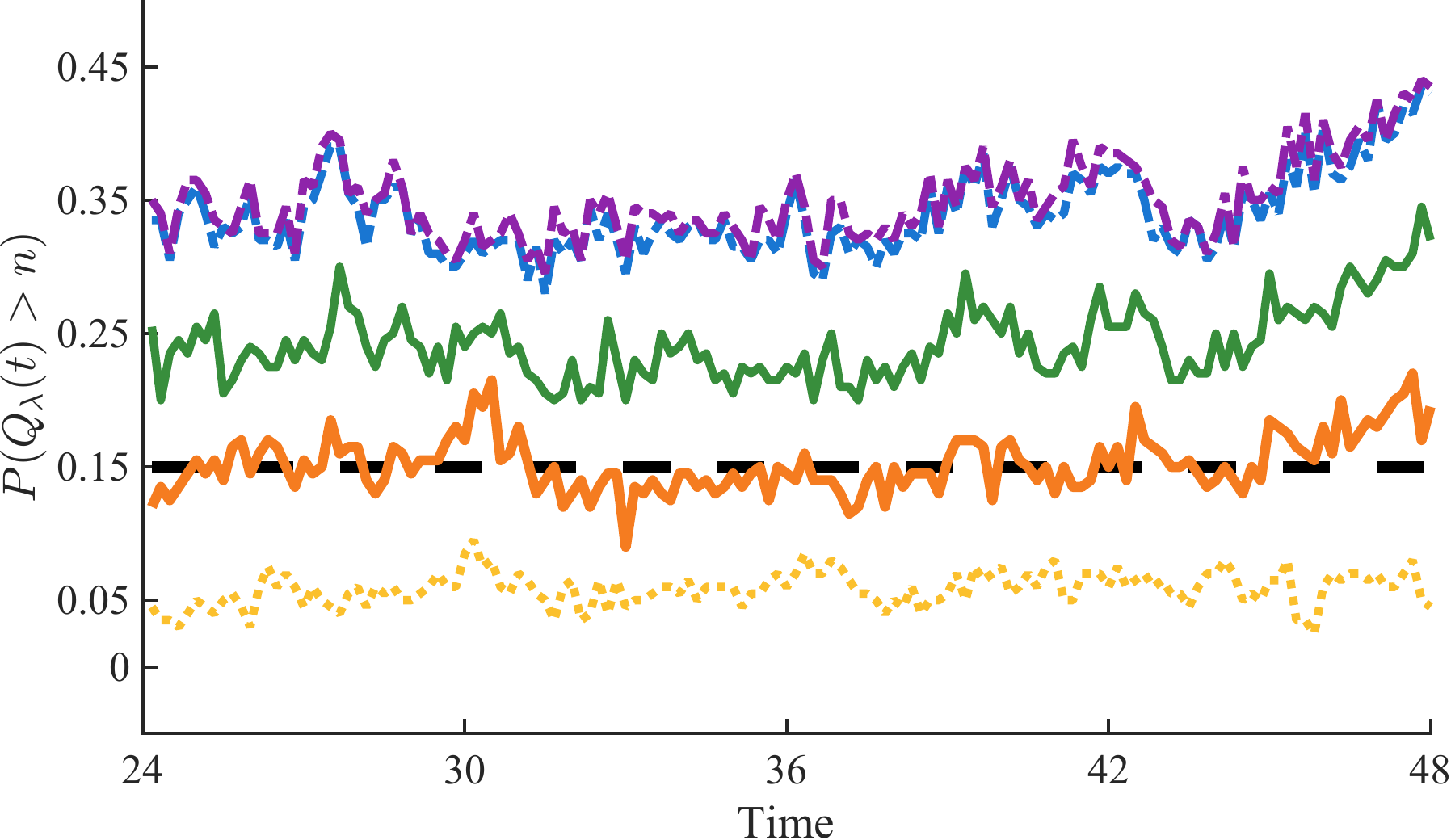}}
{Performance of Staffing Rules in Infinite-Server Systems with Arrival Model $\mathcal{M}_5$  ($\lambda=2400$). \label{fig:sim-inf-large-lambda}}
{The probability that the number of customers in an infinite-server system exceeds $n$, where $n$ varies depending on the staffing rule. Results exclude the first 24 hours (warm-up period).
Left: $\varepsilon = 0.05$. Right: $\varepsilon = 0.15$.}
\end{figure}

We run experiments with $\lambda=150$, $600$, and $2400$.
Figure~\ref{fig:sim-inf-large-lambda} presents results for $\lambda=2400$, 
and the remaining results ($\lambda=150$ and $600$) are provided in Section~\ref{sec:add-res-Sec6}. 
Under the basic alpha safety rule (i.e.,  $n = n^\ast_{\mathsf{Basic}}$), 
the estimates of $\mathbb{P}(Q_\lambda(t) > n)$ after warm-up closely match the target value, even for relatively small mean arrival rates ($\lambda=150$, Figure~\ref{fig:sim-inf-small-lambda}). This alignment improves as $\lambda$ increases, with estimates becoming nearly identical to $\varepsilon$ when $\lambda=2400$ (Figure~\ref{fig:sim-inf-large-lambda}).

The alpha safety rule under $\mathcal{M}_3$ yields estimates of $\pr(Q_\lambda(t) > n)$ significantly higher than the target $\varepsilon$ across all settings, indicating under-staffing. 
While $\mathcal{M}_3$ captures Taylor's law, it fails to account for temporal correlation decay in arrival data generated from $\mathcal{M}_5$, causing under-estimation of intensity uncertainty through mis-estimated dispersion parameter $\alpha$.

The square-root safety rules under $\mathcal{M}_1$ and $\mathcal{M}_4$ result in substantial under-staffing, with $\mathbb{P}(Q_\lambda(t) > n)$ significantly exceeding the target $\varepsilon$. 
This gap widens as $\lambda$ increases, indicating these rules systematically underestimate stochastic variability in over-dispersed arrivals, an issue that intensifies with system scale.

With the linear safety rule under $\mathcal{M}_2$, 
the estimates of $\mathbb{P}(Q_\lambda(t) > n)$ exhibit distinct behaviors: exceeding target $\varepsilon$ for small $\lambda$, but falling below it---eventually to zero---as $\lambda$ increases. 
In small systems, $\mathcal{M}_2$ understaffs for two reasons. 
First, it underestimates intensity uncertainty by failing to capture temporal correlation decay. 
Second, its safety level coefficient, which is $\mathcal{O}(\lambda)$, is calibrated for heavy-traffic limits ($\lambda\to\infty$) and may be smaller than coefficients in other models like $\mathcal{M}_5$, resulting in inadequate safety levels for small systems. 
However, in large systems, this rule's linear scaling with $\lambda$ leads to overstaffing, contradicting the dispersion scaling described by Taylor's law.

\subsection{Finite-server Systems}\label{sec:fin-ser-sym}

\begin{figure}[t]

\FIGURE{ 
    \includegraphics[width=0.5\textwidth]{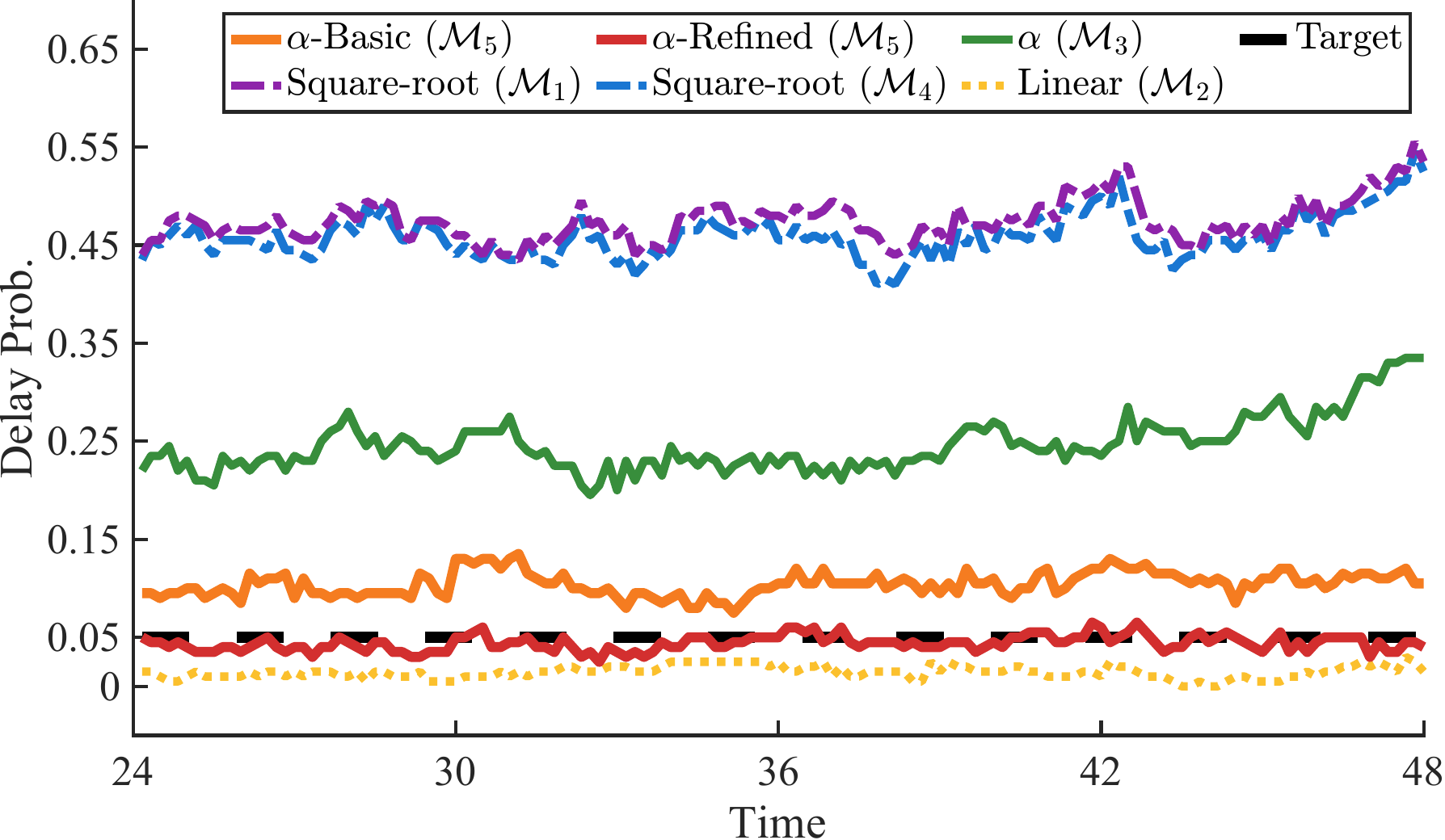}
	\includegraphics[width=0.5\textwidth]{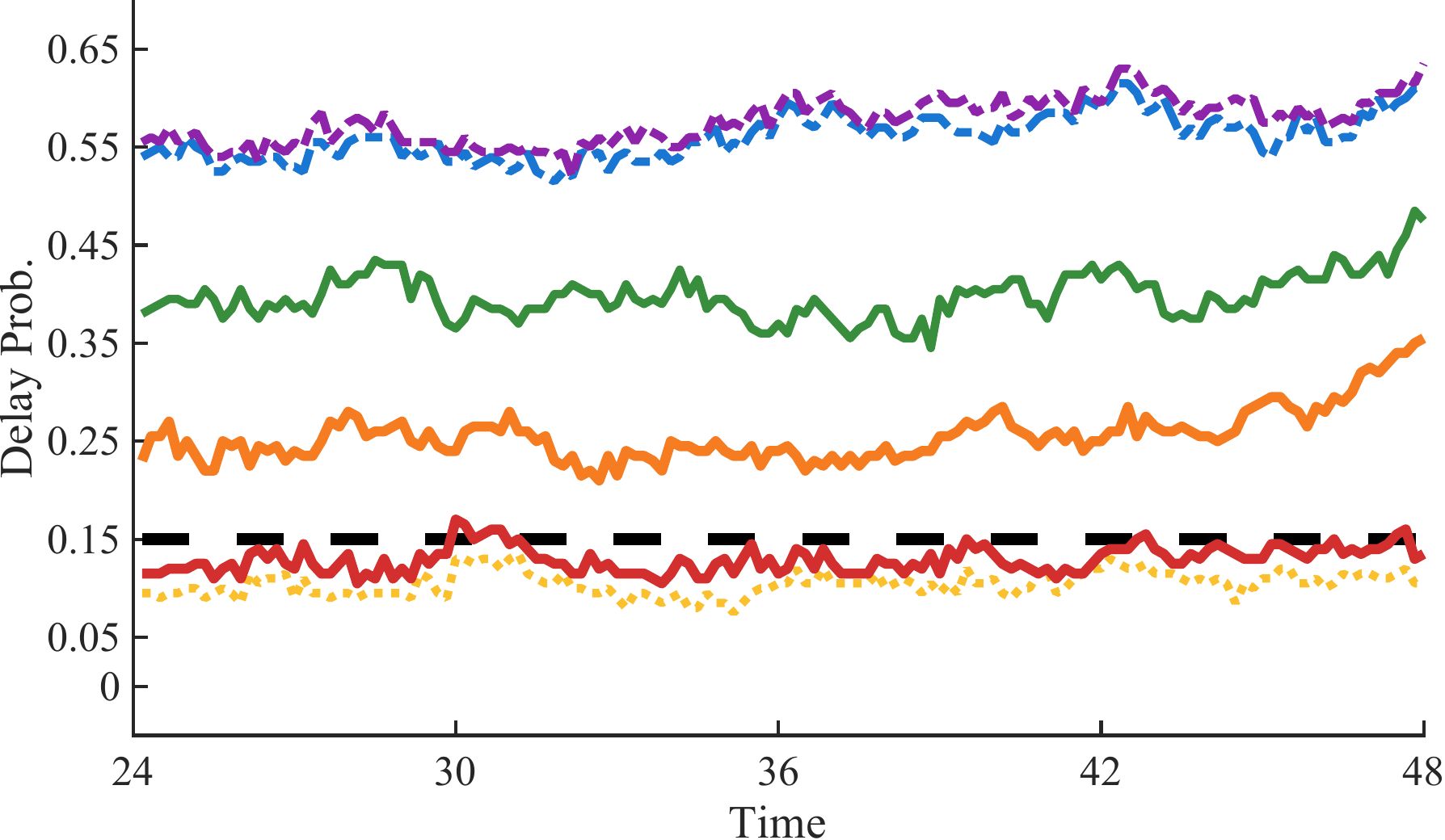}
	}
{Performance of Staffing Rules in Finite-server Systems with Arrival Model $\mathcal{M}_5$  ($\lambda=2400$). \label{figsim121}}
{The delay probability as defined by 
$\pr( Q_\lambda^{(n)}(t) > n)$, where $n$ varies depending on the staffing rule. 
Results exclude the first 24 hours (warm-up period).
Left: $\varepsilon = 0.05$. Right: $\varepsilon = 0.15$. }
\end{figure}

In this subsection, we validate the refined alpha safety rule through numerical simulations of finite-server queueing systems. The experimental setup mirrors that of Section~\ref{sec:inf-ser-sym}. Here, we focus on estimating the delay probability $\pr(Q_\lambda^{(n)}(t) > n)$, rather than the tail probability of the number-in-system of an infinite-server queueing system. 
Figures~\ref{figsim121} present results for $\lambda=2400$. Additional results (for $\lambda=150$ and $600$) are provided in Section~\ref{sec:add-res-Sec6}. We have the following findings.

The refined and basic alpha safety rules under $\mathcal{M}_5$ show contrasting performance. The refined alpha safety rule achieves its design objective, producing delay probabilities that match the target across all settings, with negligible deviations even for small systems ($\lambda=150$, Figure~\ref{figsim120}). 
In contrast, the basic alpha safety rule produces delay probabilities exceeding the target  because it relies on infinite-server approximation.
As detailed in Section~\ref{sec:staff_refine}, this approximation tends to underestimate delay probabilities in finite-server systems, resulting in under-staffing.

The staffing rules under $\mathcal{M}_1$ through $\mathcal{M}_4$ perform consistently with infinite-server results in Section~\ref{sec:inf-ser-sym}.
The square-root safety rules ($\mathcal{M}_1$ and $\mathcal{M}_4$) exhibit significant under-staffing that intensifies with system scale. 
For example, when $\varepsilon=0.05$, delay probabilities increase from approximately $0.3$ to $0.4$ to $0.5$ as $\lambda$ increases from $150$ to $600$ to $2400$ (Figures~\ref{figsim120}, \ref{fig:sim-finite-server-med}, and \ref{figsim121}, left pane). This  pattern reflects these rules' failure to account for dispersion scaling by Taylor's law.

The linear safety rule under $\mathcal{M}_2$ transitions from under-staffing to over-staffing as $\lambda$ increases. For $\varepsilon=0.05$, delay probabilities decrease from approximately $0.2$ to $0.1$ to $0.0$ as $\lambda$ increases from $150$ to $600$ to $2400$ (Figures~\ref{figsim120}, \ref{fig:sim-finite-server-med}, and \ref{figsim121}, left pane).

The alpha safety rule under $\mathcal{M}_3$ results in significant under-staffing, albeit less severe than the two square-root safety rules. 
For instance, the delay probability associated with this rule is about $0.25$ and $0.4$ for $\varepsilon=0.05$ and $0.15$, respectively. 
As discussed in Section~\ref{sec:inf-ser-sym}, this stems from $\mathcal{M}_3$'s misspecification of the arrivals' temporal correlation structure. 
While both $\mathcal{M}_3$ and $\mathcal{M}_5$ capture Taylor's law, this misspecification leads to different estimates of the dispersion scaling parameter.

\subsection{Taylor's Law versus Temporal Correlation Decay}
\label{sec:two characteristics}

Both Taylor's law and temporal correlation structure of the arrival process affect staffing. 
As shown in Equation~\eqref{eq:basic-alpha}, 
Taylor's law determines the order of magnitude of the safety level through $\alpha$ and is therefore the main driver in heavy traffic as $\lambda\to\infty$. 
Temporal correlation affects the safety level through the leading constant multiplying $\lambda^{\frac{\alpha+1}{2}}$.
When $\lambda$ is not large or $\alpha$ is small, this effect can be substantial and may even outweigh the effect of Taylor's law.

To better understand which effect matters more for a given dataset, 
we compare the performance of $\mathcal{M}_3$ and $\mathcal{M}_4$, since each captures only one of the two effects: Taylor’s law or temporal correlation decay.
We use $\mathcal{M}_5$ as the data-generating model and vary $\alpha$ and $\kappa$. 
Larger $\kappa$ means faster temporal correlation decay. 
Specifically, we consider three cases, all with $\sigma=0.5$:
\begin{enumerate}[label=(\Roman*)]
    \item \emph{High-$\alpha$, slow decay}  $(\alpha=0.8,\kappa=0.1)$: this case favors $\mathcal{M}_3$;
    \item \emph{High-$\alpha$, fast decay} $(\alpha=0.8,\kappa=0.25)$: this case does not clearly favor either model;
    \item \emph{Low-$\alpha$, fast decay} $(\alpha=0.1,\kappa=0.25)$: this case favors $\mathcal{M}_4$.
\end{enumerate}

The experiments follow the same procedures as those in Section~\ref{sec:fin-ser-sym}. 
The results,  reported in  
Table~\ref{table:taylorvscorrelation}, show that while neither model achieves the target service quality, the relative benefit of modeling Taylor's law versus temporal correlation decay depends on the characteristics of the arrival data. 
In Case~(I), $\mathcal{M}_3$ consistently outperforms $\mathcal{M}_4$, with delay probabilities closer to the target, because the arrival data has slow temporal correlation decay and a strong Taylor's law effect. 
In contrast,
in Case~(III), $\mathcal{M}_4$ consistently outperforms $\mathcal{M}_3$ because these characteristics are reversed.  

In Case~(II), the ranking depends on the system scale. 
$\mathcal{M}_4$ performs better when the scale is relatively small ($\lambda=500$), suggesting that modeling temporal correlation decay is more important than Taylor's law  in the regime with high $\alpha$, high $\kappa$, and small $\lambda$. 
However, $\mathcal{M}_3$ becomes more competitive and eventually performs better as $\lambda$ increases ($\lambda=2000,5000$). 
This is because the safety levels under $\mathcal{M}_3$ and $\mathcal{M}_4$ are of order $\mathcal{O}(\lambda^{\frac{\alpha+1}{2}})$ and $\mathcal{O}(\lambda^{\frac{1}{2}})$, respectively. 
Thus, in heavy traffic, Taylor's law becomes dominant and $\mathcal{M}_3$ outperforms $\mathcal{M}_4$.

\begin{table}[t]
\centering
    \TABLE{Taylor's Law versus Temporal Correlation Decay. \label{table:taylorvscorrelation}}
    {
    \begin{tabular}[c]{c@{\extracolsep{10pt}}
    c@{\extracolsep{15pt}}
    c@{\extracolsep{15pt}}
    c@{\extracolsep{15pt}}
    c@{\extracolsep{10pt}}
    c@{\extracolsep{10pt}}
    c@{\extracolsep{15pt}}
    c@{\extracolsep{15pt}}
    c@{\extracolsep{10pt}}
    c@{\extracolsep{10pt}}
    c@{\extracolsep{15pt}}
    c@{\extracolsep{15pt}}
    c@{\extracolsep{10pt}}
    c@{\extracolsep{10pt}}
    }
    \toprule
 \multirow{3}{*}{Target}   & \multirow{3}{*}{$\lambda$}  &  \multicolumn{4}{c}{\shortstack{(I) High-$\alpha$, Slow decay %
 \\$(\alpha=0.8,\kappa=0.1)$}} & 
 \multicolumn{4}{c}{\shortstack{(II) High-$\alpha$, Fast decay %
 \\$(\alpha=0.8,\kappa=0.25)$}} & 
 \multicolumn{4}{c}{\shortstack{(III) Low-$\alpha$, Fast decay %
 \\$(\alpha=0.1,\kappa=0.25)$}} \\
 \cmidrule(lr){3-6}\cmidrule(lr){7-10}\cmidrule(lr){11-14}
  &  & 
\multicolumn{2}{c}{Staffing} &  \multicolumn{2}{c}{Delay prob.}  &
\multicolumn{2}{c}{Staffing} &  \multicolumn{2}{c}{Delay prob.}  &
\multicolumn{2}{c}{Staffing} &  \multicolumn{2}{c}{Delay prob.} \\
\cmidrule(lr){3-4}\cmidrule(lr){5-6}\cmidrule(lr){7-8}
\cmidrule(lr){9-10}\cmidrule(lr){11-12}\cmidrule(lr){13-14}
& & $\mathcal{M}_3$ & $\mathcal{M}_4$ & $\mathcal{M}_3$ & $\mathcal{M}_4$ & $\mathcal{M}_3$ & $\mathcal{M}_4$ & $\mathcal{M}_3$ & $\mathcal{M}_4$ & $\mathcal{M}_3$ & $\mathcal{M}_4$ & $\mathcal{M}_3$ & $\mathcal{M}_4$ \\
\midrule
\multirow{3}{*}{0.05}
 & 500 &
 128 & 113 & \underline{0.32} & 0.43 & 104 & 107  & 0.45 & \underline{0.41} & 87 & 89 & 0.61 & \underline{0.46} \\
 & 2000 & 
 457 & 392  &  \underline{0.35} & 0.51  &    389 & 380  &  \underline{0.50} & 0.54  &  341 & 344  &  0.58 & \underline{0.46}  \\
  & 5000 & 
 1077 & 926  &  \underline{0.44} & 0.62  &    939 & 907  &  \underline{0.59} & 0.66  &  846 & 850  &  0.61 & \underline{0.51}  \\
    \midrule
 \multirow{3}{*}{0.15}
 & 500 & 
 112 & 102  &  \underline{0.44} & 0.53  &   97 & 98  &  0.57 & \underline{0.56}  &  86 & 87  &  0.70 & \underline{0.61}   \\
 & 2000 & 
 412 & 371   & \underline{0.46} & 0.57   & 369 & 363   & \underline{0.60} & 0.63  &   338 & 340  &   0.71 & \underline{0.62}   \\
 & 5000 & 
 988 & 892  &  \underline{0.54} & 0.66  &    901 & 880  &  \underline{0.67} & 0.72  &  842 & 844  &  0.72 & \underline{0.66}  \\
    \bottomrule
    \end{tabular}
    }
    {%
    The underlined numbers indicate whether $\mathcal{M}_3$ or $\mathcal{M}_4$ has a delay probability closer to the target. }
\end{table}

Furthermore, comparing Cases~(I) and~(II) reveals that the performance of $\mathcal{M}_3$ deteriorates as $\kappa$ increases. 
For example, when $\epsilon=0.05$ and $\lambda=5000$, the delay probability under $\mathcal{M}_3$ increases from $0.44$ to $0.59$ as $\kappa$ becomes larger.
With faster temporal correlation decay, overall temporal dependence weakens. 
When fitting $\mathcal{M}_3$, which assumes constant temporal correlation, the model compensates for this mismatch by estimating a smaller $\alpha$ to better align with the weaker dependence. 
A smaller estimated $\alpha$ reduces the implied safety level and exacerbates under-staffing, leading to worse service quality in Case (II) than Case (I).

Comparing Cases~(III) and~(II) shows that as $\alpha$ increases, the performance of $\mathcal{M}_4$ changes little when $\lambda$ is small, but deteriorates substantially in heavy traffic. 
For example, when $\epsilon=0.05$ and $\lambda=5000$, the delay probability under $\mathcal{M}_4$ increases from $0.51$ to $0.66$ as $\alpha$ increases from $0.1$ to $0.8$.
This is because, in heavy traffic, 
Taylor's law is the main driver of the staffing level.  
Thus, the drawback of $\mathcal{M}_4$, which does not capture Taylor's law, becomes more severe as $\alpha$ increases.

These results show that distinguishing the relative importance of Taylor's law and temporal correlation decay can be complex in practice, especially when the system is not in heavy traffic. 
A key advantage of our model is that it captures both effects. 
Practitioners can therefore avoid this issue: they can fit $\mathcal{M}_5$ to arrival data, letting the parameter estimates implicitly balance the two effects, and then compute staffing levels without having to decide in advance which effect is more important or which model is more appropriate for a given dataset.

\section{Case Study: NYC 311 Call Center with Non-stationary Arrival Data}\label{sec:casestudy}

Arrival data frequently exhibit strong time-of-day patterns (Figure~\ref{fig:time-of-day}).
In this section, we generalize our DSPP model to handle non-stationary arrivals and evaluate various safety rules using real data from the NYC 311 Call Center (see Appendix~\ref{sec:data-description} for data description). 

\subsection{Non-stationary DSPP Model with Generalized CIR Intensity}\label{sec:non-stat-DSPP}

To accommodate the time-of-day effect, we adopt the stationary independent period-by-period (SIPP) approach \citep{green2001improving,GreenKolesarWhitt07}, a standard method for staffing systems with non-stationary arrivals. 
The intensity is modeled as:
\begin{equation}\label{eq:CIR_time-varying}
\dd{X(t)} = \kappa\left( \lambda(t) -X(t) \right) \dd{t} + \sigma  \sqrt{\lambda(t)^\alpha X(t)} \dd{B(t)},
\end{equation}
where $\lambda(t)$ is a periodic, piecewise-constant function. 

Let $T$ be the length of each period. We divide each period into $k$ segments of equal length $\Delta$, with $\lambda(t)$ constant within each segment:
\begin{align}\label{eq:piecewise-constant-rate}
    \lambda(t) = \lambda_i, \quad \mbox{if  $t = T(j-1) + \displaystyle(i-1)\Delta + s$, for some $j\geq 1$, $1\leq i\leq k$, and $0\leq s < \Delta$}. 
\end{align}
For example, the NYC 311 Call Center  operates 24 hours daily. 
We set $T=24$, $k=48$, and $\Delta=1/2$ to adjust staffing every 30 minutes.

For each segment $i=1,\ldots,k$, this non-stationary arrival process follows a DSPP model with intensity process \eqref{eq:CIR} parameterized by $(\lambda_i, \kappa, \sigma,\alpha)$. 
While each segment has its distinct mean arrival rate $\lambda_i$, the parameters $\alpha$ (dispersion scaling), $\kappa$ (mean-reversion speed), and $\sigma$ (volatility) remain constant across segments. 
These parameters are estimated using the new MLE method from Section~\ref{sec:model-estimation}, leveraging heavy-traffic approximation (see Appendix~\ref{sec:MLE} for details).

Similarly, we extend $\mathcal{M}_1$ through $\mathcal{M}_4$ to non-stationarity by incorporating the piecewise-constant arrival rate $\lambda(t)$ while keeping other parameters constant across segments.  
For example, $\mathcal{M}_3$ becomes a DSPP model with intensity $\lambda(t) + \lambda(t)^{(\alpha+1)/2}Y$, where $Y$ is a zero-mean random variable with variance $\sigma_Y^2$. 
Safety rules are then computed for each segment using its corresponding parameters.

In the following subsections, we analyze phone call arrival timestamps from NYC 311 call center in 2017, excluding weekends, public holidays, and the months of June and December due to their distinct arrival patterns. Our analysis uses two datasets from 2017: 
\begin{itemize}[noitemsep]
    \item Training dataset: January through May (92 business days).
    \item Test dataset: July through November (96 business days).
\end{itemize}
Figure~\ref{fig:time-of-day} shows the daily arrival pattern from July through November, which is similar to the pattern observed from January through May. We use the January--May dataset for model estimation and the July--November dataset to evaluate staffing rules. 

However, these datasets do not contain durations of each phone call. Hence, we examine two service time distributions: log-normal and gamma, which represent heavy-tailed and light-tailed distributions, respectively. Empirical studies suggest service times in various service systems, including call centers \citep{Brown} and hospitals \citep{StrumMayVargas00}, can be well modeled by log-normal distributions. We mainly focus on log-normal service times and present results for gamma service times in Appendix~\ref{sec:additional-exp}.

\subsection{In-sample Model Selection via Information Criteria}\label{sec:MCSC}

We fit the non-stationary versions of models $\mathcal{M}_1$ through $\mathcal{M}_5$ to the training dataset, 
consisting of 30-minute arrival counts collected over 92 days during 24 working hours each day. For the non-stationary Poisson process ($\mathcal{M}_1$), MLE is straightforward since the model avoids the high-dimensional numerical integration inherent in doubly stochastic Poisson structures. 
For models $\mathcal{M}_2$ through $\mathcal{M}_5$, we use a heavy-traffic approximation to compute the MLE efficiently. Details are provided in Appendix~\ref{sec:MLE} for $\mathcal{M}_4$ and $\mathcal{M}_5$, and in Appendix~\ref{StaticDSPP} for $\mathcal{M}_2$ and $\mathcal{M}_3$.
Using the maximized log-likelihood values, we compare models through AIC and BIC criteria.

\begin{table}[t]
	\TABLE{AIC and BIC for NYC 311 Call Center Arrival Data 
    (January--May 2017). \label{tab:Likelihood-realdata}}
	{	
	\begin{tabular}
		[c]{l@{\extracolsep{40pt}}
        c@{\extracolsep{40pt}}
        c@{\extracolsep{40pt}}
        c@{\extracolsep{40pt}}
        c@{\extracolsep{40pt}}}%
		\toprule 
		Model & $\textsf{AIC}$ & $\Delta\textsf{AIC}$ & $\textsf{BIC}$ & $\Delta\textsf{BIC}$ \\
		\midrule 
$\mathcal{M}_1$ & $72417.56$  & $29255.28$ & $72417.56$ & $29247.72$   \\
$\mathcal{M}_2$ & $65773.90$  &  $22611.62$ & $65776.42$ & $22606.58$   \\
$\mathcal{M}_3$ & $62238.08$  & $19075.80$ & $62243.12$ & $19073.28$  \\		 
$\mathcal{M}_4$ & $43181.16$ & $18.88$ & $43186.20$ & $16.36$  \\
$\mathcal{M}_5$ (this paper) &  $43162.28$  &  & $43169.84$ &   \\
		\bottomrule  
	\end{tabular}}
{Models $\mathcal{M}_1$ through $\mathcal{M}_5$ are the non-stationary version with piecewise-constant mean arrival rates. %
$\Delta\textsf{AIC}$ and $\Delta\textsf{BIC}$ denote the difference in $\textsf{AIC}$ and $\textsf{BIC}$, respectively, between $\mathcal{M}_i$ and $\mathcal{M}_5$, for $i=1,2,3,4$.}
\end{table}

Table~\ref{tab:Likelihood-realdata} presents these results. In order of increasing AIC (and BIC), the models rank as $\mathcal{M}_5$, $\mathcal{M}_4$, $\mathcal{M}_3$, $\mathcal{M}_2$, and $\mathcal{M}_1$, with both criteria favoring our DSPP model ($\mathcal{M}_5$). The performance gap between $\mathcal{M}_5$ and the second-best model $\mathcal{M}_4$ is substantial, with an AIC difference of $18.88$ and a BIC difference of $16.36$. In practice, a difference exceeding $10$ is considered strong evidence favoring the model with the lower AIC (or BIC) value \citep{Raftery95,BurnhamAnderson02}.

Moreover, the AIC (and BIC) difference between $\mathcal{M}_3$ and $\mathcal{M}_5$ is enormous (over $19,000$). 
As shown in Table~\ref{tab:5models}, while both $\mathcal{M}_3$ and $\mathcal{M}_5$ capture Taylor's law, $\mathcal{M}_3$ fails to model the temporal correlation structure in the arrival data. This highlights the importance of modeling the intensity as a stochastic process (as in $\mathcal{M}_4$ and $\mathcal{M}_5$) rather than a static random variable.

\subsection{Out-of-sample Model Comparison via Staffing Performance}

We evaluate our DSPP model's staffing rule performance using the test dataset. 
Following common practice \citep{green2001improving,ShenHuang08_Pois_SVD,Avramidis04,OreshkinReegnardLEcuyer16}, we set the staffing interval to 30 minutes: the piecewise-constant mean arrival rate \eqref{eq:piecewise-constant-rate} has a period of $T=24$ hours, divided into $k=28$ segments of length $\Delta=1/2$ hours. For experiments with alternative staffing intervals (10 minutes and 1 hour), see Appendix~\ref{sec:additional-exp}.

To demonstrate how our DSPP model's ability to capture Taylor's law of dispersion scaling benefits staffing decisions during system scale-up, we estimate model parameters $(\alpha, \kappa,\sigma)$ using only nighttime arrivals (9 PM--8 AM)---a period of relatively low arrival rates---from January through May 2017. We then apply these estimated parameters to compute staffing levels across all hours in the July--November test dataset, using the refined alpha safety rule for each 30-minute segment.

The model's performance during daytime hours (8 AM--9 PM), when arrival rates substantially exceed those of the training  period, provides insights into how well the refined alpha safety rule---and by extension our DSPP model---generalizes to unseen, busier settings.  
These insights can guide managers in system expansion decisions. 
We apply this same approach when computing staffing levels under other models.
This approach parallels our methodology in Section~\ref{sec:synthetic}, where we train models on small systems and evaluate them on larger ones for stationary arrivals.

To estimate delay probabilities for various staffing rules, we simulate a queueing system that runs continuously for $M=96$ days using a single stream of phone calls from the July--November test dataset, initializing the system as empty at the start. Consistent with our previous experiments, we assume log-normally distributed service times with both mean and standard deviation equal to $1/6$ hours, corresponding to a service rate of $\mu=6$ customers per hour.

\begin{figure}[t]
\FIGURE{
    \includegraphics[width=0.5\textwidth]{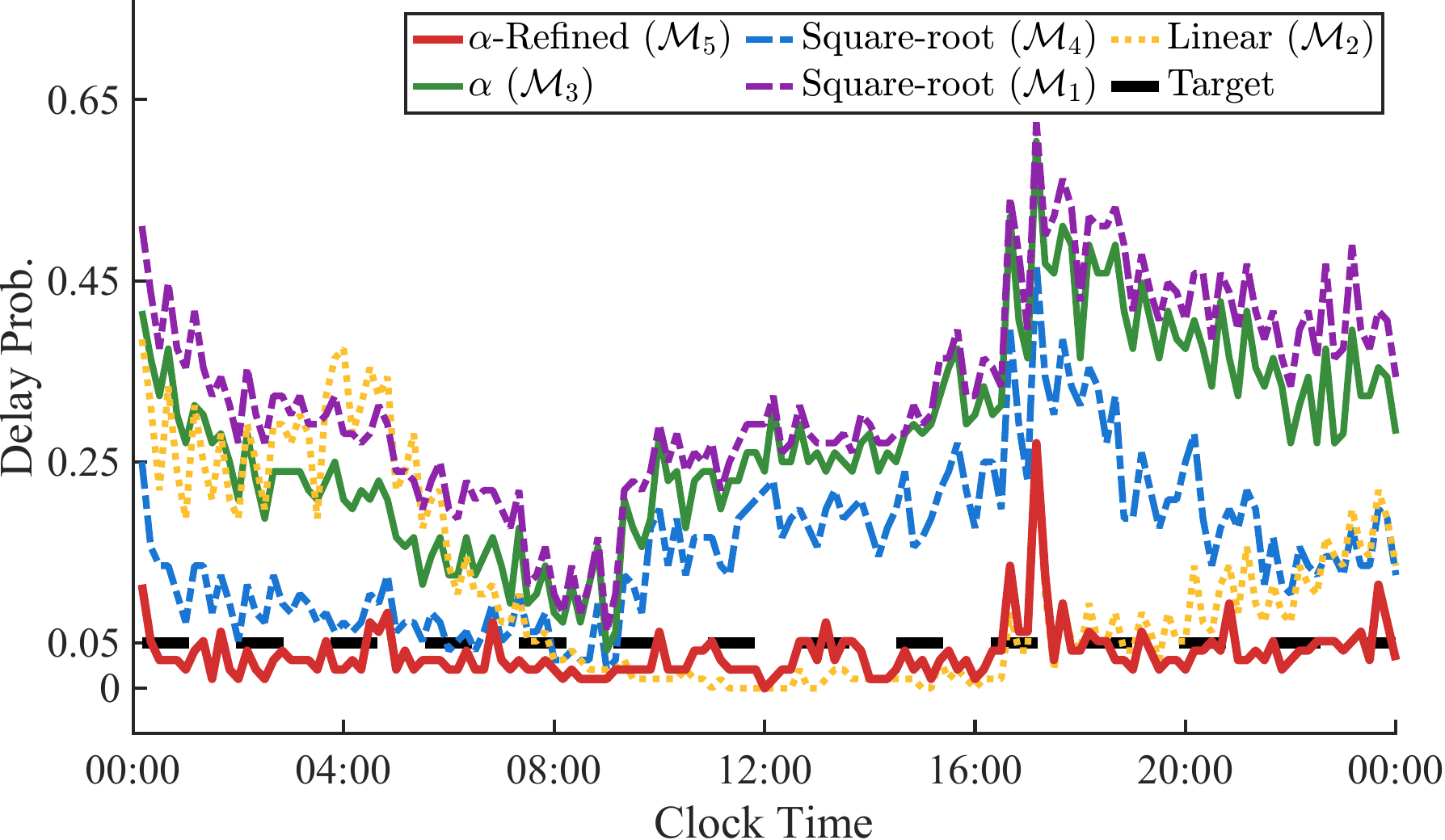}
	\includegraphics[width=0.5\textwidth]{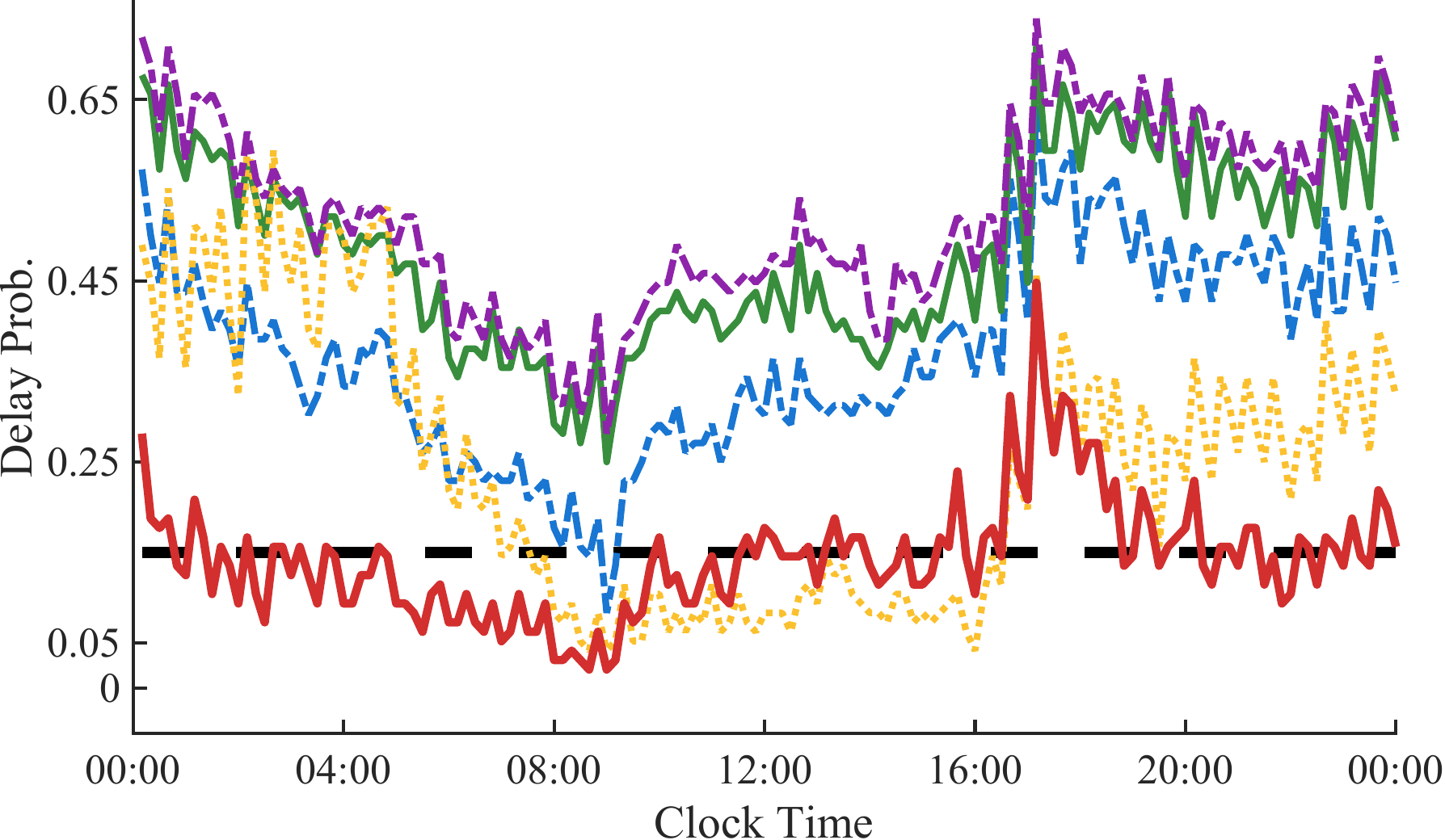}
	}
{Performance of Staffing Rules in Finite-server Systems with Real Arrival Data (SIPP Avg). \label{figsim2120}}
{Models $\mathcal{M}_1$ through $\mathcal{M}_5$ represent non-stationary processes with piecewise-constant mean arrival rates, with staffing levels computed for each 30-minute segment. The delay probability follows Equation~\eqref{delayprob00}, shown for $\varepsilon = 0.05$ (left) and $\varepsilon = 0.15$ (right). Model parameters are estimated using nighttime data (9 PM--8 AM) from January--May 2017, while delay probabilities are evaluated using full-day data from July--November 2017.  Each $\lambda_i$ is estimated using the \textsf{SIPP~Avg} approach (Appendix~\ref{sec:MLE}).}
\end{figure}

Let  $\widetilde{Q}(t )$ denote the number of customers in this queueing system at time $t$. 
We measure performance using the delay probability at different times of the day:
\begin{equation}\label{delayprob00}
\pr\big( \widetilde{Q}(t) > n(t) \big)   \approx  \frac{1}{M} \sum_{j=1}^{M} \mathbb{I}\big\{\widetilde{Q}( s ) > n(s) \big\}, \quad s=24(j-1) + t \mbox{ and } t\in[0,24],
\end{equation}
where $n(t)$ is the number of servers at time $t$, determined by a piecewise-constant staffing rule. For each time of day $t\in[0,24]$, we average observations across all $96$ days to estimate the corresponding delay probability. 
The results appear in Figure~\ref{figsim2120}.

We have three findings. Firstly, the refined alpha safety rule under our DSPP model $\mathcal{M}_5$ consistently outperforms alternative safety rules. While the square-root safety rule under $\mathcal{M}_1$ combined with SIPP serves as a standard benchmark for non-stationary arrival systems, our experiments show it generates substantial under-staffing, with delay probabilities exceeding target due to unaccounted arrival over-dispersion. This limitation also affects the square-root safety rule under $\mathcal{M}_4$.
The alpha safety rule under $\mathcal{M}_3$ similarly produces considerable under-staffing, mainly because its estimate of the parameter $\alpha$ fails to fully capture dispersion scaling in the data, a consequence of the model's static intensity not accounting for temporal correlation decay in arrivals. 
The linear safety rule under $\mathcal{M}_2$ shows time-dependent performance: it under-staffs during low-arrival nighttime hours, producing elevated delay probabilities, but over-staffs during high-arrival daytime hours. This pattern aligns with our findings for stationary arrivals in Section~\ref{sec:synthetic}.

Secondly, with a target delay probability of $\varepsilon=0.05$, the refined alpha safety rule achieves near-target delay probabilities throughout most of the day. Its performance excels during training period hours (9 PM--8 AM) but declines during daytime hours (8 AM--9 PM), an expected outcome since parameter estimation for $(\alpha, \kappa, \sigma)$ excludes daytime data. 
When $\varepsilon=0.15$, while the rule's performance shows notable deterioration, delay probabilities remain close to the target value, with modest oscillations for the majority of the time. 

Thirdly, examining the right panel of Figure~\ref{figsim2120} ($\varepsilon=0.15$) reveals two periods of unsatisfactory performance for the refined alpha safety rule. 
From 6 AM--9 AM, delay probabilities consistently fall below the target value, while from 4 PM--8 PM, they significantly exceed it. 
These performance declines stem from inaccurate mean arrival rate estimates during periods of dramatic rate changes, as shown in Figure~\ref{fig:time-of-day}: rates rise rapidly in early morning and descend quickly in late afternoon. 
While the SIPP approach assumes constant rates within each 30-minute segment, this piecewise-constant approximation poorly captures such rapid changes. 
To address this limitation, \cite{green2001improving} proposed several modifications to the SIPP approach for estimating $\lambda_i$ in each segment $i$.

To address the staffing imbalances observed in Figure~\ref{figsim2120}, we propose the \textsf{SIPP Mix} approach, which applies different SIPP variants based on time of day. We use \textsf{SIPP Min} from 6 AM--9 AM to reduce over-staffing by lowering estimated mean arrival rates, \textsf{SIPP Max} from 4 PM--8 PM to minimize under-staffing by increasing rates, and maintain \textsf{SIPP Avg} for all other periods.

\begin{figure}[t]
\FIGURE{
    \includegraphics[width=0.5\textwidth]{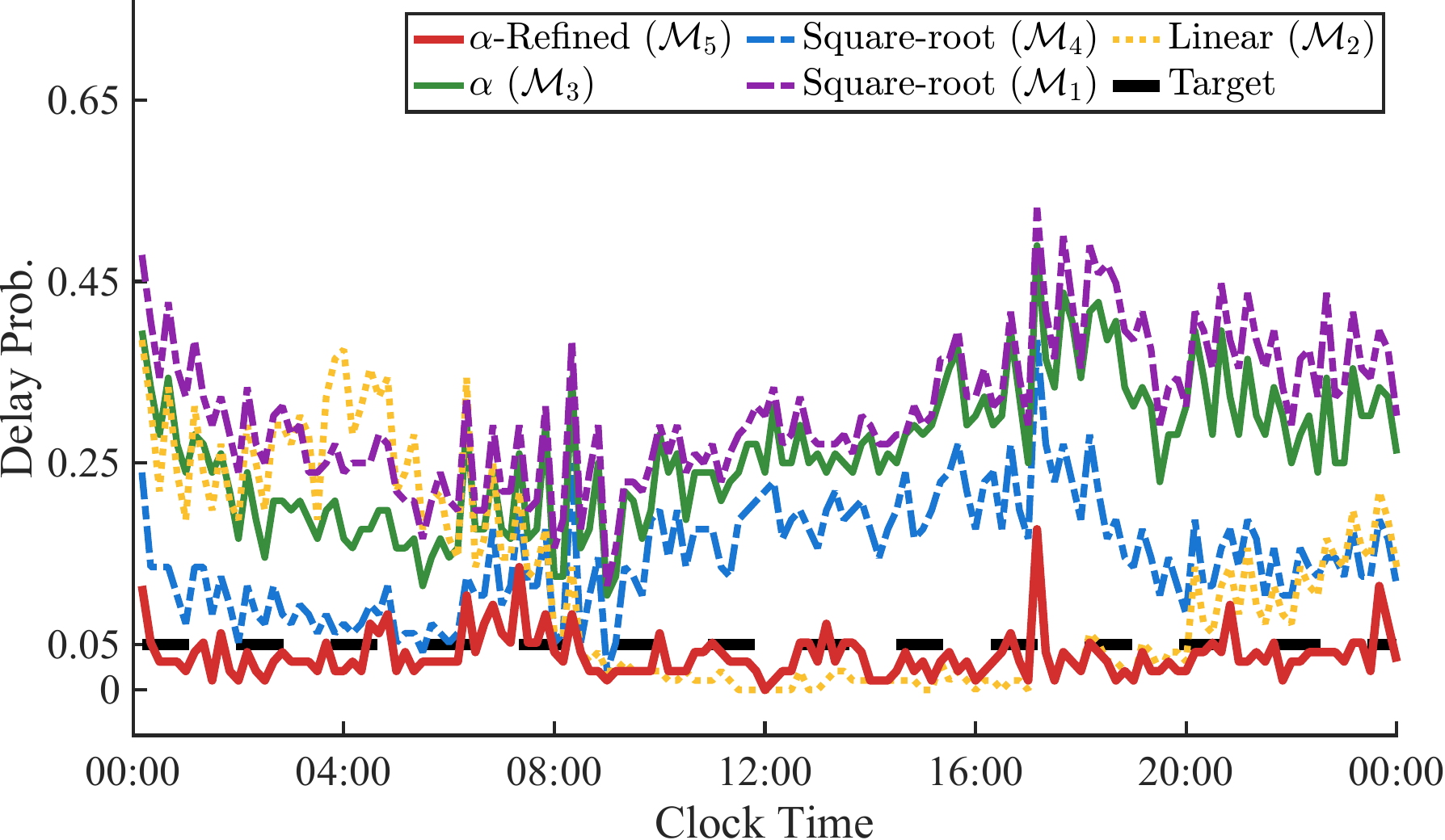}
	\includegraphics[width=0.5\textwidth]{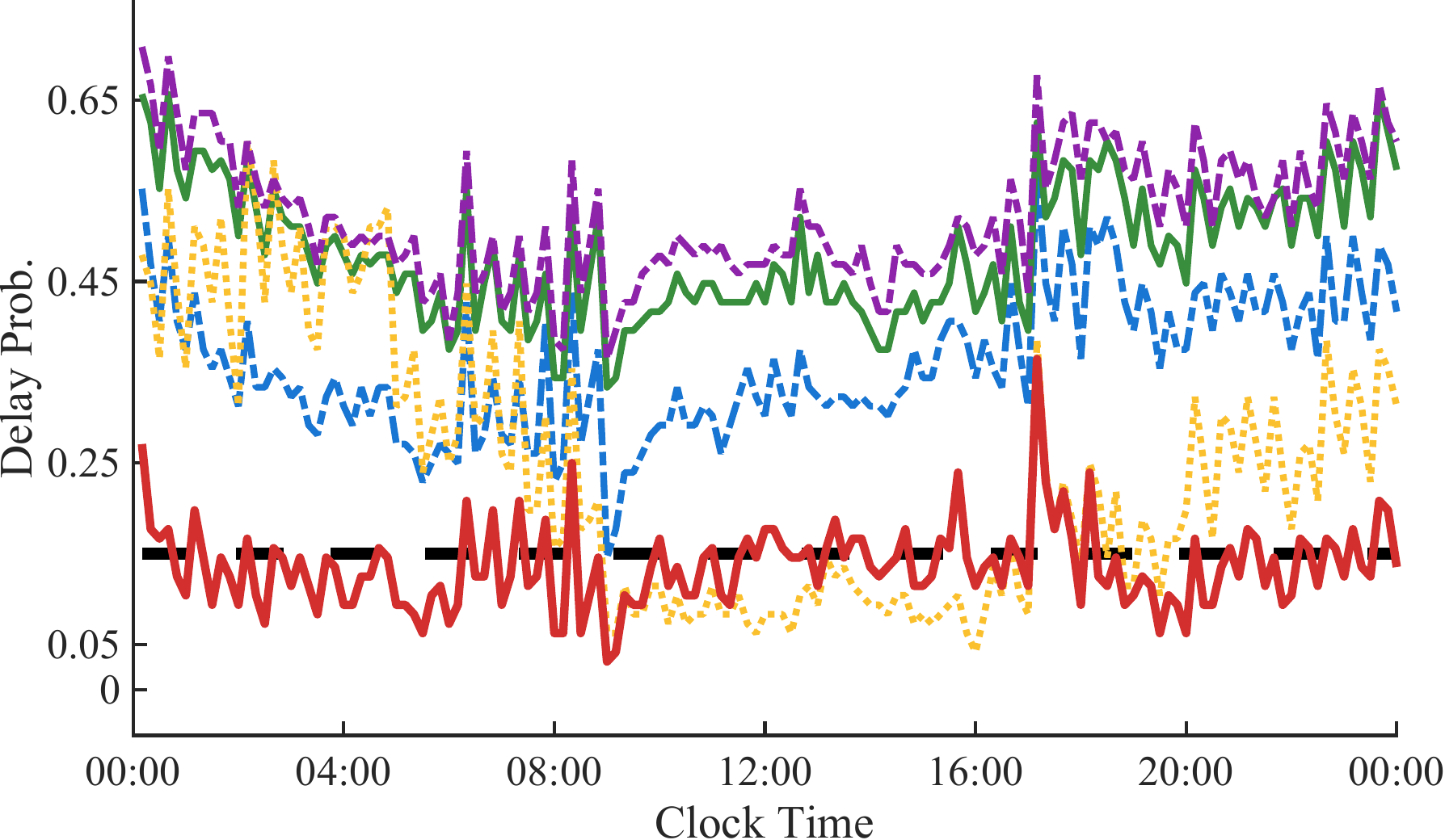}
	}
{Performance of Staffing Rules in Finite-server Systems with Real Arrival Data (SIPP Mix). \label{figsim2121}}
{Models $\mathcal{M}_1$ through $\mathcal{M}_5$ represent non-stationary processes with piecewise-constant mean arrival rates, with staffing levels computed for each 30-minute segment. The delay probability follows Equation~\eqref{delayprob00}, shown for $\varepsilon = 0.05$ (left) and $\varepsilon = 0.15$ (right). Model parameters are estimated using nighttime data (9 PM--8 AM) from January--May 2017, while delay probabilities are evaluated using full-day data from July--November 2017.  Each $\lambda_i$ is estimated using the \textsf{SIPP~Mix} approach (Appendix~\ref{sec:MLE}).}
\end{figure}

Applying \textsf{SIPP Mix} to all five staffing rules (Figure~\ref{figsim2121}) yields improved performance compared to Figure~\ref{figsim2120}, particularly for the refined alpha safety rule during the two periods (6 AM--9 AM and 4 PM--8 PM) using \textsf{SIPP Max} or \textsf{SIPP Min}. 
The refined alpha safety rule's delay probabilities now oscillate around their target values ($\varepsilon = 0.05$ or $0.15$) with a much reduced fluctuation, demonstrating both its practical utility and our DSPP model's effectiveness in handling arrival over-dispersion.

\section{Conclusion}\label{conclusion}

Motivated by empirical evidence that arrivals are over-dispersed and follow Taylor's law, we develop a novel DSPP model.
The arrival intensity is a generalized CIR process with a dispersion scaling parameter $\alpha$ that captures Taylor's law. 
The model extends the Poisson arrival model while remaining parsimonious and analytically tractable, as shown by our theoretical results.

Using this DSPP model, we identify a trichotomy in arrival over-dispersion and its impact on staffing. 
The three scenarios differ in how total stochastic variability is split between two sources: randomness in the intensity process and conditional Poisson distribution given the intensity. In the intensity-dominated scenario (the typical case where $\alpha\in(0,1)$), most variability comes from the intensity process, making it the key driver of staffing decisions needed to meet a target service level.

Building on the model, we derive alpha safety rules for staffing large-scale service systems to meet target delay probability requirements. The main feature is the safety staffing term, which grows at a rate of $\lambda^{(\alpha+1)/2}$ as the mean arrival rate $\lambda$ increases. As $\alpha$ ranges from 0 to 1, the alpha safety rules bridge the gap between the square-root safety rule (dominant in queueing literature and practice) and the linear safety rule designed for delay-free service.

This paper points to several directions for future research. Over-dispersed arrivals create long stretches of high and low demand. To maintain service quality, staffing must increase to cover the high-demand periods.  
This can leave servers idle during low-demand periods and raise operating costs.
This trade-off underscores the importance of alternative service operations, such as outsourcing \citep{KocagaArmonyWard15}, call blending \citep{PangPerry15}, and callback options \citep{AtaPeng20}. These strategies can use available capacity more effectively and improve efficiency.  
Analyzing them under our DSPP model is a natural next step, especially for systems with high stochastic variability in arrivals.

\bibliographystyle{informs2014} %
\bibliography{staffing.bib} %

\ECSwitch

\ECHead{Supplemental Material}

\EquationsNumberedBySection %

\section{Proof of Proposition~\ref{prop:squre-diffussion}} 
Part~\ref{prop:SD-positive} and part~\ref{prop:SD-stationary} can be proved using calculations similar to those found on pages 391--392 of \cite{CoxIngersollRoss85_ec}. To prove part~\ref{prop:SD-ergodic}, one can verify the Foster--Lyapunov criterion of ergodicity for continuous-time Markov processes, as outlined in Section 5 of \cite{MeynTweedie93_ec}. Theorem~1 of \cite{ZhangGlynn18_ec} also provides a relevant reference for this proof.

\section{Proof of Theorem~\ref{theo:CoD} \label{ec-sec:laplace}} 

The proof employs explicit calculations using the Laplace transform of $A(t)$, where $X(0)$ follows the stationary distribution $\pi$ established in Proposition~\ref{prop:SD-positive}.
We derive the Laplace transform of $X(t)$ before determining the Laplace transform of $A(t)$.

\begin{proposition}    \label{prop:Laplace-X}
If Assumption~\ref{asp:Feller} holds, then 
for any $\theta\geq 0$, $t\geq 0$, and $x>0$, we have
	\begin{equation}\label{eq2}
	\E\left[ \exp(-\theta \int_0^t X(s) \dd{s} ) \bigg| X(0) = x \right] = [\phi(\theta, t) ]^\gamma \exp(- \psi(\theta, t) x), 
	\end{equation}
	where $\gamma = 2\kappa \lambda^{1-\alpha}/\sigma^2$, 
	\begin{align*}
	    \phi(\theta, t) = \frac{2 \eta e^{t ( \eta+\kappa)/2} }{ \eta -\kappa + (\eta +\kappa) e^{\eta t} }, \qq{and}  
	    \psi(\theta, t) = \frac{2 \theta \left( e^{\eta t} -1 \right)}{  \eta -\kappa + (\eta +\kappa) e^{\eta t} },
	\end{align*}
	with $\eta = \sqrt{\kappa^2 + 2\sigma^2\lambda^{\alpha} \theta }$.
\end{proposition}
\proof{Proof.}
The result can be obtained through a straightforward calculation, based on Proposition~6.2.5 from \cite{LambertonLapeyre07}.
\Halmos 
\endproof

\begin{proposition}    \label{prop:Laplace-A}
If Assumption~\ref{asp:Feller} holds, then
for any $\theta\geq 0$, $t\geq 0$, and $x>0$, we have
    \begin{align*}
	\mathcal{L}_{\pi, t}(\theta) = \E_\pi\left[e^{-\theta A(t)}\right] =  \left[\phi(1-e^{-\theta}, t) \right]^\gamma \left( 1+ \frac{\sigma^2 \lambda^{\alpha}}{2\kappa} \psi(1-e^{-\theta},t) \right)^{-\gamma}, 
	\end{align*}
where $\gamma$, $\phi$, and $\psi$ are defined in Proposition~\ref{prop:Laplace-X}. 
\end{proposition}

\proof{Proof.} By virtue of the doubly stochastic Poisson structure, 
\begin{align*}
\E\left[e^{-\theta A(t)} \big| X(0) = x\right] ={}& \E\biggl[\E\biggl[e^{-\theta  A(t)} \bigg| \int_0^t X(s)\dd{s}\biggr] \bigg| X(0) = x\biggr] \\ 
={}& \E\biggl[ \exp(-(1-e^{-\theta}) \int_0^t X(s) \dd{s}) \bigg| X(0) = x\biggr] \\ 
={}& \left[\phi(1-e^{-\theta}, t) \right]^\gamma \exp(- \psi(1-e^{-\theta}, t) x),
\end{align*}
where the last step follows from Proposition~\ref{prop:Laplace-X}. 
Then, integrating over $x$ with respect to $\pi(\dd{x})$ yields 
\begin{align*}
\E_\pi[e^{-\theta A(t)}] ={}&  \int_0^\infty \left[\phi(1-e^{-\theta}, t) \right]^\gamma \exp(- \psi(1-e^{-\theta}, t) x) \pi(\dd{x}) \\
={}& \left[\phi(1-e^{-\theta}, t) \right]^\gamma \left( 1+ \frac{\sigma^2 \lambda^{\alpha}}{2\kappa} \psi(1-e^{-\theta},t) \right)^{-\gamma},
\end{align*}
where the second step holds because $\pi$ is a gamma distribution and we may apply the formula of its Laplace transform. 
\Halmos \endproof

\proof{Proof of Theorem~\ref{theo:CoD}.}
Let $\mathcal{L}_{\pi, t}(\theta)$ denote the Laplace transform of $A(t)$ conditional on $X(0)\sim\pi$, as derived from Proposition~\ref{prop:Laplace-A}. Then, we have:
\[
\E_\pi[A(t)] = -\frac{\dd }{\dd{\theta}} \mathcal{L}_{\pi, t}(0) \qq{and}
\E_\pi[A(t)^2] = \frac{\dd^2 }{\dd{\theta^2}} \mathcal{L}_{\pi, t}(0).
\]

The expressions for the mean and variance of $A(t)$ conditional on $X(0)\sim\pi$ can be determined through an elementary yet lengthy calculation. 
We omit the details for brevity.\Halmos
\endproof

\section{Proof of Theorem \ref{thm:FCLTofA}}\label{sec:proof-scaled-arrival}

In Section~\ref{sec:outline}, we provide an outline of the proof. 
In Section~\ref{lemmas}, we introduce a series of lemmas that form the foundation for the following proofs. In Section~\ref{EC12}, we establish the convergence of the scaled intensity process. 
Following this, in Section~\ref{EC13}, we provide a formal proof for Theorem~\ref{thm:FCLTofA}.

\subsection{Proof Outline}\label{sec:outline}

First, we prove that $\widehat{X}_\lambda\Rightarrow U$ in $\mathcal{C}[0,T]$ as $\lambda\to \infty$, where $\mathcal{C}[0,T]$ denotes the space of all continuous real-valued functions on $[0,T]$ and is equipped with the uniform metric \citep{Billingsley99}. 
To this end, we rewrite Equation~\eqref{eq:CIR} as 
\[ \widehat{X}_\lambda(t) = \widehat{X}_\lambda(0) + \sigma \int_0^t \sqrt{\bar{X}_\lambda(s)} \dd{B(s)} - \kappa \int_0^t \widehat{X}_\lambda(s)\dd{s}, \]
where $\bar{X}_\lambda := \lambda^{-1} X_\lambda$. 
By Theorem 4.1 of \cite{PangTalrejaWhitt07_ec}, we can regard the above integral equation as a continuous mapping $h$ that maps $\bigl(\widehat{X}_\lambda(0),  \int_0^\cdot \sqrt{\bar{X}_\lambda(s)} \dd{B(s)} \bigr)$ to $\widehat{X}_\lambda$.
Here, $h$ maps any $(b, y) \in \mathbb{R}\times \mathcal{C}[0,T]$ to $x\in \mathcal{C}[0,T]$ that solves the equation
\begin{equation}\label{eq:cont-map}
x(t) = b + \sigma y(t) - \kappa \int_0^t x(s)\dd{s}. 
\end{equation}
Then, we establish the weak convergence of $\widehat{X}_\lambda$ in three steps:
\begin{enumerate}[label=(\Roman*)]
    \item We use Prohorov's theorem to show that $\bar{X}_\lambda  \Rightarrow {\bf 1}$ in $\mathcal{C}[0,T]$ as $\lambda\to \infty$, where ${\bf 1}$ denotes a constant process taking the value $1$ for all $t$. Specifically, we verify that any finite-dimensional distributions of $\bar{X}_\lambda$ converge to a vector of ones, and that $\bar{X}_\lambda$ is tight with respect to $\lambda$. 
    \item We show that as $\lambda\to\infty$,
    \begin{equation} \label{eq:joint-convergence}
    \biggl(\widehat{X}_\lambda(0),  \int_0^\cdot \sqrt{\bar{X}_\lambda(s)} \dd{B(s)} \biggr)\Rightarrow (X_\dagger, B) \quad \mbox{in $\mathbb{R}\times \mathcal{C}[0,T]$}.
    \end{equation}
    While $\widehat{X}_\lambda(0) \Rightarrow X_\dagger$ is an assumption, we prove that $\int_0^\cdot \sqrt{\bar{X}_\lambda(s)} \dd{B(s)} \Rightarrow B$ in $\mathcal{C}[0,T]$ by employing the weak limit theory for stochastic integrals developed by \cite{kurtzprotter1991}, since the integrand $\sqrt{\bar{X}_\lambda } \Rightarrow {\bf 1}$ in $\mathcal{C}[0,T]$ as shown in Step (I). 
    We then extend the two marginal convergences to the joint convergence. 
    \item We apply the mapping $h$, as defined by Equation~\eqref{eq:cont-map}, to the weak convergence \eqref{eq:joint-convergence}.
    The continuous mapping theorem implies that $\widehat{X}_\lambda$ weakly converges in $\mathcal{C}[0,T]$ to the solution to the integral equation 
    \[U(t) = X_\dagger + \sigma B(t) - \kappa \int_0^t U(s) \dd{s},\]
    which is precisely the OU process \eqref{OU}.
\end{enumerate}
Building upon the weak convergence of $\widehat{X}_\lambda$, we prove Equation~\eqref{thm:FCLT:eq1} by  decomposing the process $\widehat{A}_\lambda$ into two processes, i.e., $\widehat{A}_\lambda= \widehat{A}_{\lambda,1}+\widehat{A}_{\lambda,2}$, where 
$$
\widehat{A}_{\lambda,1}(t) = \lambda^{-\frac{\alpha+1}2} \left[ A_\lambda(t) -  \int_0^t X_\lambda(s) \dd{s}\right] \quad \mbox{and} \quad \widehat{A}_{\lambda,2}(t) = \int_0^t \widehat{X}_\lambda(s) \dd{s}.
$$
Then, we establish the weak convergence of $\widehat{A}_\lambda$ in three steps: 
\begin{enumerate}[resume*]
    \item Consider $\hat{x}_\lambda\in\mathcal{C}[0,T]$, a convergent sequence with a limit $\hat{x}\in\mathcal{C}[0,T]$. Given $\widehat{X}_\lambda= \hat{x}_\lambda$, the process $\widehat{A}_{\lambda,1}$ becomes a square-integrable martingale when $\alpha=0$, and the process $\widehat{A}_{\lambda,2}$ becomes a deterministic process. 
    By the martingale functional central limit theorem, we can show that
\[
\left(\widehat{A}_{\lambda,1}, \widehat{A}_{\lambda,2}\right) \Big| \left\{\widehat{X}_\lambda= \hat{x}_\lambda \right\} \ \Rightarrow \ \left(\widetilde{B}\mathbb{I}\{\alpha = 0\}, \int_0^\cdot x(s) \dd s \right) \quad \text{in $\mathcal{D}[0,T]\times \mathcal{C}[0,T]$}.
\]

    \item For any bounded and continuous functional $f: \mathcal{D}[0,T]\times\mathcal{C}[0,T] \to \mathbb{R}$, we introduce two bounded functionals, namely $g_\lambda: \mathcal{C}[0,T]\to \mathbb{R}$ and $g: \mathcal{C}[0,T]\to \mathbb{R}$ as follows:
\[
g_\lambda(\hat{x})=\E\left[ f\big( \widehat{A}_{\lambda,1}, \widehat{A}_{\lambda,2} \big) \Big| \widehat{X}_\lambda= \hat{x} \right] \quad \text{and} \quad g(\hat{x})= \E\left[ f\left( \widetilde{B}, \int_0^\cdot \hat{x}(s)\dd{s} \right) \right]. 
\] 
    By the definition of weak convergence along with Step (IV), we can show that $g_\lambda(\hat{x}_\lambda) \to g(\hat{x})$ as $\lambda\to\infty$.
    
    \item It follows from the generalized continuous mapping theorem \citep[Theorem~3.4.4]{Whitt02} and the previously proven result $\widehat{X}_\lambda\Rightarrow U$ in Step~(III) that $ g_\lambda(\widehat{X}_\lambda)\Rightarrow g(U)$. 
    Then, applying the definition of weak convergence, we have 
    \[
    \left( \widehat{A}_{\lambda,1}, \widehat{A}_{\lambda,2} \right) \Rightarrow \left( \widetilde{B}\mathbb{I}\{\alpha = 0\}, \int_0^\cdot U(s) \dd s \right) \quad \mbox{ in $\mathcal{D}[0,T]\times \mathcal{C}[0,T]$}. 
    \] 
\end{enumerate}

\subsection{Technical Lemmas on Weak Convergence of Stochastic Processes}\label{lemmas}

For any $T>0$, let $\mathcal{C}[0,T]$ denote the set of all continuous real-valued functions on the interval $[0,T]$, equipped with the \emph{uniform metric} $\|\cdot\|_T$, which is defined as $\|x-y\|_T=\sup_{t\in[0,T]} |x(t)-y(t)|$ for all $x,y\in\mathcal{C}[0,T]$. 
Moreover, let $\mathcal{D}[0,T]$ denote the set of all right-continuous real-valued functions with left limits on the interval $[0,T]$, equipped with the standard Skorohod $J_1$ topology. 
We use ``$\Rightarrow$'' to denote the weak convergence. 
For more details on the space $\mathcal{C}$, the space $\mathcal{D}$, and weak convergence, we refer to \cite{Billingsley99_ec} and \cite{Whitt02_ec}.  

The following convergence results extend the fundamental properties of weak convergence from random variables to stochastic processes. As the proofs for these lemmas are straightforward, we have opted to omit them.

\begin{lemma}\label{lemma:weak-convergence}
Fix $T>0$. 
For any $\lambda>0$, 
let $S_\lambda$, $\tilde{S}_\lambda$, and $S_\dagger$ be random elements in the space $\mathcal{C}[0,T]$,
let $L_\lambda$ and  $L_\dagger$ be random elements in the space $\mathcal{D}[0,T]$, 
and let $\zeta_\lambda$ and $\zeta$ be real-valued random variables. 
Furthermore, let $\BFc_\lambda$ and $\BFc$ be constant elements in the space $\mathcal{C}[0,T]$; that is, $\BFc_\lambda(t) \equiv c_\lambda$ and $\BFc(t) \equiv c$ for some constants $c_\lambda$ and $c$, for all $t\in[0,T]$.
\begin{enumerate}[label=(\roman*)]
    \item \label{part:from-weak-to-prob} If $S_\lambda \Rightarrow \BFc$ in $\mathcal{C}[0,T]$ as $\lambda\to\infty$, then $S_\lambda \rightarrow \BFc$ in probability as $\lambda\to\infty$. 
    \item \label{part:joint-weak-R-and-C} Suppose $S_\lambda$ is independent of $\zeta_\lambda$ for all $\lambda>0$. 
    If  $\zeta_\lambda\Rightarrow \zeta_\dagger$ in $\mathbb{R}$ and $S_\lambda \Rightarrow S_\dagger$ in $\mathcal{C}[0,T]$  as $\lambda\to\infty$, then $(\zeta_\lambda, S_\lambda) \Rightarrow (\zeta_\dagger, S_\dagger)$ in $\mathbb{R}  \times \mathcal{C}[0,T]$ as $\lambda\to\infty$. 
    \item \label{part:joint-weak-D-and-C} If $L_\lambda \Rightarrow L_\dagger$ in $\mathcal{D}[0,T]$ and $\BFc_\lambda \to \BFc$ in $\mathcal{C}[0,T]$ as $\lambda \to\infty$, then $(L_\lambda, \BFc_\lambda) \Rightarrow (L_\dagger,\BFc)$ in $\mathcal{D}[0,T]\times \mathcal{C}[0,T]$ as $\lambda\to\infty$. 
    \item \label{part:joint-weak-additive} If $(\zeta_\lambda, S_\lambda) \Rightarrow (\zeta_\dagger, S_\dagger)$ in $\mathbb{R}\times \mathcal{C}[0,T]$ and $\| S_\lambda - \tilde{S}_\lambda \|_T \to 0$ in probability as $\lambda\to\infty$, 
    then $(\zeta_\lambda, \tilde{S}_\lambda) \Rightarrow (\zeta_\dagger, S_\dagger)$ in $\mathbb{R}\times \mathcal{C}[0,T]$ as $\lambda\to\infty$. 
\end{enumerate}
\end{lemma}

The following lemma originates from the classical weak limit theory of stochastic integrals. While the general theorem is constructed on the semimartingale framework, we adapt it to a specific case that involves the It\^{o} integral with respect to  Brownian motion.

\begin{lemma}[Theorem 2.2 of \cite{KurtzProtter91_ec}]\label{App1:lem3}
Fix $T>0$. 
Let $S_\lambda$ and $S_\dagger$ be random elements in the space $\mathcal{C}[0,T]$ for any $\lambda>0$. 
Suppose $S_\lambda$ is adapted to the filtration generated by a standard Brownian motion $B$. 
Define $\mathcal{I}_\lambda(t)=\int_0^t S_\lambda(s) \dd B(s)$  
and $\mathcal{I}_\dagger(t)=\int_0^t S_\dagger(s) \dd B(s)$ for any $t\in[0,T]$.  
\begin{enumerate}[label=(\roman*)]
\item If $(S_\lambda,B)\Rightarrow (S_\dagger,B)$ in $\mathcal{C}^2[0,T]$ as $\lambda\to\infty$, 
then $(S_\lambda,B,\mathcal{I}_\lambda )\Rightarrow (S_\dagger,B,\mathcal{I}_\dagger)$ in $\mathcal{C}^3[0,T]$ as $\lambda\to\infty$. 
\item If $(S_\lambda,B)\to (S_\dagger,B)$ in probability as $\lambda\to\infty$, then $(S_\lambda,B,\mathcal{I}_\lambda )\to (S_\dagger,B,\mathcal{I}_\dagger)$ in probability as $\lambda\to\infty$. 
\end{enumerate} 
\end{lemma}

The following lemma is a generalized version of the continuous-mapping theorem. While the standard continuous-mapping theorem applies to a single continuous mapping, the generalized version extends to a sequence of continuous mappings.

\begin{lemma}[Theorem 3.4.4 of \cite{Whitt02_ec}]%
\label{GCM}
Fix $T>0$. For any $\lambda>0$, 
let $g_\lambda$ and $g$ be measurable functions mapping $\mathcal{C}[0,T]$ to $\mathbb{R}$. 
Let $S_\lambda$ and $S_\dagger$ be random elements in the space $\mathcal{C}[0,T]$. 
Moreover, let $E$ be the set of $x\in \mathcal{C}[0,T]$ such that $g_\lambda(x_\lambda)\to g(x)$ fails for some sequence $\{x_\lambda \in\mathcal{C}[0,T]:\lambda >0\}$ with $x_\lambda\to x$ in $\mathcal{C}[0,T]$ as $\lambda\to\infty$. 
If $S_\lambda\Rightarrow S_\dagger$ in $\mathcal{C}[0,T]$ as $\lambda\to\infty$ and $\pr(S_\dagger \in E)=0$,
then $g_\lambda(S_\lambda)\Rightarrow g(S_\dagger)$ as $\lambda\to\infty$. 
\end{lemma}

\subsection{Weak Convergence of the Scaled Intensity Process}\label{EC12}

For the ease of reference, we restate the definition of the generalized CIR process in Equation~\eqref{eq:CIR} as follows: 
\begin{equation}\label{eq:CIR_ec}
\dd{X_\lambda(t)} = \kappa\left( \lambda -X_\lambda(t) \right) \dd{t} + \sigma  \sqrt{\lambda^\alpha X_\lambda(t)} \dd{B(t)},
\end{equation}
with $X(0) >0$ almost surely, where $B$ is a standard Brownian motion.  
In addition, we define 
\[\bar{X}_\lambda(t)  =\lambda^{-1} X_\lambda(t)\quad\mbox{ and }\quad \widehat{X}_\lambda(t)  =  \lambda^{-\frac{\alpha+1}2}[X_\lambda(t)-\lambda]. \]
In the following, we first prove that the fourth moment of $\bar{X}_\lambda(t)$ is uniformly bounded with respect to $\lambda$ and $t$, provided that the same holds true for the initial value  $\bar{X}_\lambda(0)$. 
Then, we show prove $\widehat{X}_\lambda$ weakly converges to an OU process. 

\begin{lemma}\label{App1:lem1}
Fix $T>0$. 
Suppose that Assumptions~\ref{asp:Feller}--\ref{asp:alpha} and $\sup_{\lambda\geq 1} \E[\bar{X}_\lambda^4(0)] <\infty$ hold. Then, for $k=1,2,3,4$, 
\begin{equation*}%
	\sup_{\lambda\geq 1} \sup_{0\leq t\leq T} \E[\bar{X}_\lambda^k(t)]  < \infty. 
\end{equation*}
\end{lemma}

\proof{Proof.}
Dividing both sides of Equation \eqref{eq:CIR_ec} by $\lambda$, we have 
\begin{equation}
	\bar{X}_\lambda(t) = \bar{X}_\lambda(0) + \kappa \int_0^t [1-\bar{X}_\lambda(s)]\dd{s} + \sigma_\lambda \int_0^t \sqrt{\bar{X}_\lambda(s)} \dd{B(s)}, 
	\label{eqn:intensitybyn}
\end{equation}
where $\sigma_\lambda \coloneqq  \sigma  \lambda^{(\alpha-1)/2}$. 
Notice that $\sigma_\lambda\downarrow 0$ as $\lambda\to\infty$, because $0\leq \alpha<1$.

It follows from It\^o's lemma that, for $k=1,2,3,4$, 
\[
\bar{X}_\lambda^k(t)= \bar{X}_\lambda^k(0) +  \int_0^t \left( \kappa k \bar{X}_\lambda^{k-1}(s) -\kappa k \bar{X}_\lambda^{k}(s) +\frac{\sigma_\lambda^2}2 (k-1)k \bar{X}_\lambda^{k-1}(s)  \right) \dd{s} + k \sigma_\lambda \int_0^t \bar{X}_\lambda^{k-\frac12}(s) \dd{B(s)}. 
\]
By Young's inequality, we have $k \bar{X}_\lambda^{k-1}(s) \leq 1+ (k-1) \bar{X}_\lambda^{k}(s)$. 
Therefore, 
\[
\E[\bar{X}_\lambda^k(t)] \leq \E[\bar{X}_\lambda^k(0)] + \left( \frac{\sigma_\lambda^2}2 (k-1) + \kappa \right) t +  \left( \frac{\sigma_\lambda^2}2 (k-1)^2 -\kappa \right) \int_0^t \E[\bar{X}_\lambda^{k}(s) ]  \dd{s}. 
\]
By Gr\"onwall's inequality, we have 
$$
\E[\bar{X}_\lambda^k(t)] \leq \left[ \E[\bar{X}_\lambda^k(0)] + \left( \frac{\sigma_\lambda^2}2 (k-1) + \kappa \right) t  \right] \cdot  \exp\left[\left( \frac{\sigma_\lambda^2}2 (k-1)^2 - \kappa \right) t \right]. 
$$
Because $\sup_{\lambda\geq 1} \E[\bar{X}_\lambda^4(0)]<\infty$ and $\sigma_\lambda\downarrow 0$ as $\lambda\to\infty$, we conclude the proof.
\Halmos
\endproof

\begin{theorem}\label{FCLTofintensity}
Suppose that 
\begin{enumerate*}[label=(\roman*)]
	\item \label{part:common-assump} Assumptions~\ref{asp:Feller}--\ref{asp:alpha} hold,
	\item \label{part:indepen-B} $X_\lambda(0)$ is independent of $B$, 
	\item \label{part:uniform-bound} $\sup_{\lambda \geq 1} \E[\widehat{X}_\lambda^4(0)]<\infty$, and 
	\item \label{part:converg-init} $\widehat{X}_\lambda(0)\Rightarrow X_\dagger $ as $\lambda\to\infty$ for some random variable $X_\dagger$.
\end{enumerate*} 
Then, for any $T>0$, 
\begin{equation}\label{eq:X_hat_convergence}
	\widehat{X}_\lambda \Rightarrow U \quad \mbox{in $\mathcal{C}[0,T]$}, 
\end{equation}
as $\lambda\rightarrow\infty$, where $U$ is the OU process \eqref{OU}, that is, 
$\dd{U(t)} = -\kappa U(t) \dd{t} +\sigma \dd{B(t)}$ with the initial condition $U(0)=X_\dagger$. 
\end{theorem}

\proof{Proof.} 
We first prove $X_\lambda$ satisfies a functional weak law of large numbers as follows:
\begin{equation}\label{FWLLNintensity}
    \bar{X}_\lambda  \Rightarrow {\bf 1} \quad \mbox{in $\mathcal{C}[0,T]$}, 
\end{equation}
as $\lambda\rightarrow \infty$, 
where ${\bf 1}$ denotes a constant process taking the value $1$ for all $t\in[0,T]$.

Calculations by \cite{CoxIngersollRoss85_ec} show that for any $t\geq 0$,
\begin{align}
	\E[ X_\lambda(t) | X_\lambda(0)]  ={}&  X_\lambda(0) e^{-\kappa t} + \lambda \big(1-e^{-\kappa t}\big),  \nonumber \\ %
	\Var[X_\lambda(t) | X_\lambda(0)] ={}& \ X_\lambda(0)\frac{\sigma^2\lambda^\alpha}{\kappa} \big(e^{-\kappa t}- e^{-2\kappa t}\big) + \frac{\sigma^2\lambda^{\alpha+1}}{2\kappa}\big(1-e^{-\kappa t}\big)^2. \nonumber %
\end{align}
It follows from the tower property of conditional expectations and the law of total variance that 
\begin{align}
\E[X_\lambda(t)]  ={}&  \E[X_\lambda(0)]  e^{-\kappa t} + \lambda \big(1-e^{-\kappa t}\big), \label{eqn45} \\
\Var[X_\lambda(t)]  ={}&  \E[X_\lambda(0)] \frac{\sigma^2\lambda^\alpha}{\kappa} \big(e^{-\kappa t}- e^{-2\kappa t}\big) + \frac{\sigma^2\lambda^{\alpha+1}}{2\kappa}\big(1-e^{-\kappa t}\big)^2 + \Var[X_\lambda(0) ] e^{-2\kappa t}. 
\label{eqn46}
\end{align}

Because of Condition~\ref{part:uniform-bound}, $\widehat{X}_\lambda(0)$ is uniformly integrable with respect to $\lambda$, so  
by Condition~\ref{part:converg-init},
\[\lim_{\lambda\to\infty} \E[\widehat{X}_\lambda(0)]= \E[X_\dagger].\]
Thus, $\E[X_\lambda(0)] = \lambda+ \mathcal{O}(\lambda^{(\alpha+1)/2})$. 
It then follows from Equation~\eqref{eqn45} that for any $t\geq 0$,
\begin{equation}\label{eq:X_lambda_E_order}
\E[X_\lambda(t)] =  \lambda+ \mathcal{O} \big(\lambda^{\frac{\alpha+1}2}\big).    
\end{equation}
Moreover, by Conditions~\ref{part:common-assump} and \ref{part:uniform-bound},  $\E[X_\lambda^2(0)]= \lambda^2 + \mathcal{O}(\lambda^{(\alpha+3)/2})$. 
Hence, for any $t\geq 0$,
\begin{align}
    \Var[X_\lambda(0)]={}& \E[X_\lambda^2(0)] - (\E[X_\lambda(0)])^2 \nonumber \\ 
    ={}& \Big[ \lambda^2 + \mathcal{O}\big(\lambda^{\frac{\alpha+3}2}\big) \Big] - \Big[ \lambda+ \mathcal{O}\big(\lambda^{\frac{\alpha+1}2}\big)\Big]^2  \nonumber \\ 
    ={}& \mathcal{O}\big(\lambda^{\frac{\alpha+3}2}\big),\label{eq:X_lambda_var_order}
\end{align}
where the last step holds because $\alpha < 1$. 
Combining Equations~\eqref{eqn46}--\eqref{eq:X_lambda_var_order} yields
\begin{equation*}
\Var[X_\lambda(t)]= \mathcal{O}\big(\lambda^{\frac{\alpha+3}2}\big). 
\end{equation*}
Then, we apply Chebyshev's inequality to deduce that for any $t\geq 0$  and  $\varepsilon>0$, 
\begin{equation*}
    \pr \left( \left| \bar{X}_\lambda(t) - \E[X_\lambda(t)]  \right| \geq \varepsilon \right)  \leq  \frac1{\lambda^2\varepsilon^2} \Var[X_\lambda(t)] = \mathcal{O}\big(\lambda^{\frac{\alpha-1}2}\big)\to 0 
\end{equation*}
as $\lambda\to\infty$, 
where the convergence to zero holds because $\alpha<1$. 
Therefore, $\bar{X}_\lambda(t) - \E[\bar{X}_\lambda(t)] \to 0$ in probability as $\lambda\to\infty$, for any $t\in [0,T]$.  
Note that for any $t\in[0,T]$, $\E[\bar{X}_\lambda(t)]\to 1$ as $\lambda\to\infty$, so $\bar{X}_\lambda(t) \Rightarrow 1$ as $\lambda\to\infty$. 
Furthermore, 
because the limit is a constant, 
it is straightforward to show that  for any $n\geq 1$ and $t_1,\ldots,t_n \in[0,T]$, the random vector $(\bar{X}_\lambda(t_1), \ldots, \bar{X}_\lambda(t_n))$ weakly converges to the corresponding finite-dimensional distribution of the constant process ${\bf 1}$. 

Given the convergence of the finite-dimensional distributions, 
Prohorov's theorem implies that to prove \eqref{FWLLNintensity}, 
it suffices to prove that the sequence of stochastic processes $\{\bar{X}_\lambda: \lambda >0\}$ is \emph{tight} in $\mathcal{C}[0,T]$. 
(For Prohorov's theorem and the definition of tightness, we refer to  \cite{Billingsley99_ec}.) 
According to 
Theorem~11.6.5 of \cite{Whitt02_ec},
this tightness property can be shown by verifying two conditions:
\begin{enumerate}[label=(\arabic*)]
    \item $\{\bar{X}_\lambda(0): \lambda > 0\}$  is tight in $\mathbb{R}$, and 
    \item there exists positive constants $c_1$, $c_2$, and $C$ such that for all $\lambda\geq 1$ and $0\leq s\leq t\leq 1$,
    \begin{equation}\label{eq:moment-criterion-tightness}
    \mathbb{E}\big[|\bar{X}_\lambda(t) - \bar{X}_\lambda(s)|^{c_1}\big] \leq C |t-s|^{1 + c_2}.
    \end{equation}
\end{enumerate}
The tightness of $\bar{X}_\lambda(0)$ with respect to $\lambda$ can be verified directly from the definition of tightness, Condition~\ref{part:uniform-bound}, and the application of Markov's inequality. 

We verify the condition \eqref{eq:moment-criterion-tightness} in the following.  
Applying Jensen's inequality to Equation~\eqref{eqn:intensitybyn} yields 
\begin{equation} \label{eqn:MCC1} 
     \E\left[ |\bar{X}_\lambda(t) - \bar{X}_\lambda(s)|^4 \right] 
\leq  8 \E\biggl[ \biggl| \kappa \int_s^{t} [1-\bar{X}_\lambda(u)] \dd{u} \biggr|^4 \biggr] + 8\sigma_\lambda^4 \E\biggl[ \biggl|\int_s^{t} \sqrt{\bar{X}_\lambda(u)} \dd{B(u)}\biggr|^4 \biggr].
\end{equation}
By H\"older's inequality, we have 
\begin{align}  
   \E\biggl[ \biggl| \kappa \int_s^{t} [1-\bar{X}_\lambda(u)] \dd{u} \biggr|^4 \biggr] \leq{}&  \kappa^4 (t-s)^3 \E\left[ \int_s^{t} [1-\bar{X}_\lambda(u)]^4 \dd{u} \right] \nonumber \\    
   \leq{}& \kappa^4 (t-s)^3 \E\left[ \int_s^{t} 8 [1 + \bar{X}_\lambda(u)^4] \dd{u} \right] \nonumber \\ 
   \leq{}& 8\kappa^4 (1 + \tilde{C}^4) (t-s)^4, 
   \label{eqn:MCC1-1} 
\end{align}
where the second step follows from Jensen's inequality, and 
\[
\tilde{C}\coloneqq \sup_{\lambda\geq 1} \sup_{0\leq t\leq T} \E[\bar{X}_\lambda^4(t)].  
\]
Furthermore, it follows from Condition~\ref{part:uniform-bound} and Lemma~\ref{App1:lem1} that $\tilde{C} < \infty$.

In addition, the fact that $\tilde{C} < \infty$ enables us to apply the inequality (3.25) on page~163 \cite{KaratzasShreve91_ec}.
This inequality, a consequence of the Burkholder--Davis--Gundy inequality, allows us to deduce that
\begin{align} 
    \E\biggl[ \biggl|\int_s^{t} \sqrt{\bar{X}_\lambda(u)} \dd{B(u)}\biggr|^4 \biggl] \leq{}&  36(t-s) \E\left[ \int_s^t \bar{X}_\lambda^2(u)\dd u \right] \nonumber \\ 
    \leq{} & 36 \tilde{C}^2 (t-s)^2 
    \label{eqn:MCC1-2} . 
\end{align}
Combining \eqref{eqn:MCC1}--\eqref{eqn:MCC1-2} yields that for all $0\leq s\leq t\leq 1$, 
\begin{align*}
    \E\left[ |\bar{X}_\lambda(t) - \bar{X}_\lambda(s)|^4 \right] 
\leq{}& 64 \kappa^4 (1 + \tilde{C}^4) (t-s)^4 + 288 \sigma_\lambda^4 \tilde{C}^2 (t-s)^2  \\ 
\leq{}& 64 \kappa^4 (1 + \tilde{C}^4) (t-s)^2 + 288 \sigma^4 \tilde{C}^2 (t-s)^2,  
\end{align*}
where the second step holds because $0\leq t-s\leq 1 $ and $\sigma_\lambda = \sigma \lambda^{(\alpha-1)/2} \leq \sigma$ for all $\lambda \geq 1$. 
Consequently, the inequality \eqref{eq:moment-criterion-tightness} holds, which concludes the proof of the weak convergence \eqref{FWLLNintensity}.

\smallskip

We now proceed to prove the weak convergence \eqref{eq:X_hat_convergence}. 
Using the definitions of $\bar{X}_\lambda$ and $\widehat{X}_\lambda$, 
we can rewrite Equation~\eqref{eq:CIR_ec} as 
\begin{equation}
    \widehat{X}_\lambda(t) = \widehat{X}_\lambda(0) + \sigma \underbrace{ \int_0^t \sqrt{\bar{X}_\lambda(s)} \dd{B(s)}}_{\coloneqq Y_\lambda(t)} - \kappa \int_0^t \widehat{X}_\lambda(s)\dd{s}. 
\label{eqn:intensitybyn1}
\end{equation}
By Theorem 4.1 of \cite{PangTalrejaWhitt07_ec}, we can view Equation \eqref{eqn:intensitybyn1} as a continuous mapping $h:  \mathbb{R}\times \mathcal{C}[0,T] \mapsto \mathcal{C}[0,T]$ which maps $(b,y)$ to $x$, where $x$ is the solution to the following integral equation: 
\begin{equation}
    x(t) = b + \sigma y(t) - \kappa \int_0^t x(s)\dd{s}. 
\label{eqn:intensitybyn2}
\end{equation}

Because $\bar{X}_\lambda \Rightarrow \BFone$ in $\mathcal{C}[0,T]$, 
we know from the continuous-mapping theorem and Lemma~\ref{lemma:weak-convergence}\ref{part:from-weak-to-prob} that 
$\sqrt{\bar{X}_\lambda} \to {\bf 1}$ in probability. 
Then, applying Lemma~\ref{App1:lem3} leads to 
$\| Y_\lambda-  B \|_T \to 0$ in probability.
Moreover, it follows from Condition~\ref{part:indepen-B}, Condition~\ref{part:converg-init}, and Lemma~\ref{lemma:weak-convergence}\ref{part:joint-weak-R-and-C} that $(\widehat{X}_\lambda(0), B)\Rightarrow (X_\dagger, B)$ in $\mathbb{R}\times \mathcal{C}[0,T]$. 
Thus, by Lemma~\ref{lemma:weak-convergence}\ref{part:joint-weak-additive}, we have  
\[\big( \widehat{X}_\lambda(0), Y_\lambda \big)\Rightarrow (X_\dagger, B) \quad \mbox{in $\mathbb{R}\times \mathcal{C}[0,T]$}.
\]

Lastly, applying the continuous-mapping theorem, 
we have 
\[
h(\widehat{X}_\lambda(0), Y_\lambda)\Rightarrow h(X_\dagger, B) \quad\mbox{in } \mathcal{C}[0,T].
\]
Note that $h(X_\dagger, B)$, as defined by Equation~\eqref{eqn:intensitybyn2}, is the solution to 
\[
U(t) = X_\dagger + \sigma B(t) - \kappa \int_0^t U(s)\dd{s} .
\]
This is precisely the OU process \eqref{OU}. Therefore, we conclude the proof of Theorem~\ref{FCLTofintensity}.
\Halmos
\endproof

\subsection{Proof of Theorem~\ref{thm:FCLTofA}}\label{EC13} 
We divide the process $\widehat{A}_\lambda$ into two components: 
$
\widehat{A}_\lambda= \widehat{A}_{\lambda,1}+\widehat{A}_{\lambda,2}$, 
where 
\begin{equation}
    \widehat{A}_{\lambda,1}(t)  = \lambda^{-\frac{\alpha+1}2}  \left[ A_\lambda(t) -  \int_0^t X_\lambda(s) \dd{s} \right] \quad \textrm{ and } \quad
    \widehat{A}_{\lambda,2}(t) =  \int_0^t \widehat{X}_\lambda(s) \dd{s}. \label{defA1A2}
\end{equation}
To simplify the presentation, we introduce the following notation. For any $\hat{x}\in\mathcal{C}[0,T]$, we denote the conditional process $\widehat{A}_{\lambda} | \widehat{X}_\lambda = \hat{x}$ as $\widehat{A}_{\lambda}^{\hat x}$. In other words, $\widehat{A}_{\lambda}^{\hat x}$ has the same distribution as $\widehat{A}_{\lambda}$ given the condition $\widehat{X}_\lambda = \hat{x}$.
Likewise, we define $\widehat{A}_{\lambda,1}^{\hat x}$ and $\widehat{A}_{\lambda,2}^{\hat x}$. 

We first consider the case $\alpha=0$. 
For any $\hat{x}_\lambda, \hat{x}\in \mathcal{C}[0,T]$ such that $\|\hat{x}_\lambda- \hat{x}\|_T \to 0$ as $\lambda\to\infty$, we observe that under the condition $\widehat{X}_\lambda=\hat{x}_\lambda$, the intensity process $X_\lambda=\lambda+\lambda^{1/2}\hat{x}_\lambda$ remains positive for all sufficiently large $\lambda$. 
Hence, for all sufficiently large $\lambda$, $A^{\hat{x}_\lambda}_\lambda$ is a non-homogeneous Poisson process. 
Additionally, 
\[\widehat{A}^{\hat{x}_\lambda}_{\lambda,1}(t) = \lambda^{-\frac{1}{2}}  \left[ A^{\hat{x}_\lambda}_\lambda(t) -  \int_0^t [\lambda + \lambda^{\frac{1}{2}}\hat{x}_\lambda(s)]  \dd{s} \right]  \]
is a square-integrable martingale with a predictable quadratic variation of $\int_0^t [1+ \lambda^{-\frac12} \hat{x}_\lambda(s)] \dd{s}$, according to Lemma~3.1 of \cite{PangTalrejaWhitt07_ec}. 
Notably, this predictable quadratic variation converges to $t$ as $\lambda\to\infty$. 
Invoking the martingale functional central limit theorem (see, for example, Theorem~2.1 of \citealt{Whitt07_ec}), we deduce that as $\lambda\to\infty$, 
\[
\widehat{A}^{\hat{x}_\lambda}_{\lambda,1} \Rightarrow \widetilde{B} \quad \text{in $\mathcal{D}[0,T]$},
\]
where $\widetilde{B}$ is an independent standard Brownian motion. 

Furthermore, given the condition $\widehat{X}_\lambda= \hat{x}_\lambda$, %
\begin{align*}
\left\|\widehat{A}_{\lambda,2}^{\hat{x}_\lambda} - \int_0^\cdot \hat{x}(s)\dd{s}\right\|_T 
={}& \left\|\int_0^\cdot \hat{x}_\lambda(s)\dd{s}  - \int_0^\cdot \hat{x}(s)\dd{s}\right\|_T 
 \\
\leq {}& \int_0^T \|\hat{x}_\lambda- \hat{x}\|_T \dd{s} \to 0, 
\end{align*}
as $\lambda\to \infty$.  
It then follows from  Lemma~\ref{lemma:weak-convergence}\ref{part:joint-weak-D-and-C} that as $\lambda\to\infty$,
\begin{equation}\label{eqn521}
\left(\widehat{A}^{\hat{x}_\lambda}_{\lambda,1},\widehat{A}^{\hat{x}_\lambda}_{\lambda,2} \right) \Rightarrow \left(\widetilde{B}, \int_0^\cdot \hat{x}(s)\dd{s}\right) \quad \mbox{in $\mathcal{D}[0,T]\times \mathcal{C}[0,T]$}.
\end{equation}
For any bounded and continuous functional $f: \mathcal{D}[0,T]\times\mathcal{C}[0,T] \to \mathbb{R}$, we define two bounded continuous functionals $g_\lambda: \mathcal{C}[0,T] \mapsto \mathbb{R}$ and $g: \mathcal{C}[0,T]\mapsto \mathbb{R}$ as follows:
\[
g_\lambda(\hat{x})\coloneqq \E\left[ f\big( \widehat{A}_{\lambda,1}, \widehat{A}_{\lambda,2} \big) \Big| \widehat{X}_\lambda= \hat{x} \right] \quad \text{and} \quad g(\hat{x}) \coloneqq \E\left[ f\left( \widetilde{B}, \int_0^\cdot \hat{x}(s)\dd{s} \right) \right]. 
\]
By the definition of weak convergence, Equation~\eqref{eqn521} yields
$g_\lambda(\hat{x}_\lambda) \to g(\hat{x})$ as $\lambda\to\infty$. 
It then follows from Theorem~\ref{FCLTofintensity} and Lemma~\ref{GCM} (i.e., the generalized continuous-mapping theorem) that $g_\lambda\big(\widehat{X}_\lambda\big) \Rightarrow g(U)$ as $\lambda\to\infty$. 
Furthermore, 
because $g_\lambda$ is bounded, we have $\E\big[g_\lambda\big(\widehat{X}_\lambda\big)\big] \to \E[g(U)]$ as $\lambda\to\infty$, 
which is equivalent to 
\begin{equation}\label{eq:joint-converg-A1A20}
\lim_{\lambda\to\infty}\E\left[ f\big( \widehat{A}_{\lambda,1}, \widehat{A}_{\lambda,2} \big) \right] %
=\E\left[ f\left( \widetilde{B}, \int_0^\cdot U(s) \dd{s} \right) \right],
\end{equation}
due to the tower property of conditional expectations. 
By the definition of weak convergence, Equation~\eqref{eq:joint-converg-A1A20} yields
\begin{equation}\label{eq:joint-converg-A1A2}
\left( \widehat{A}_{\lambda,1}, \widehat{A}_{\lambda,2} \right)  \Rightarrow  \left( \widetilde{B}, \int_0^\cdot U(s)\dd{s} \right) \quad \mbox{in $\mathcal{D}[0,T]\times \mathcal{C}[0,T]$},     
\end{equation}
as $\lambda\to\infty$. 
Lastly, we conclude from the continuous-mapping theorem that if $\alpha=0$, then as $\lambda\to\infty$, 
\[\widehat{A}_\lambda =\widehat{A}_{\lambda,1}+ \widehat{A}_{\lambda,2} \Rightarrow \int_0^\cdot U(s) \dd{s} + \widetilde{B}\quad \mbox{in } \mathcal{D}[0,T].
\]

We now proceed to the case $0<\alpha<1$. 
In this case, 
\[\widehat{A}^{\hat{x}_\lambda}_{\lambda,1} =  \lambda^{-\frac{\alpha+1}2}  \left[ A^{\hat{x}_\lambda}_\lambda(t) -  \int_0^t [\lambda + \lambda^{\frac{\alpha+1}{2}}\hat{x}_\lambda(s)]  \dd{s}  \right] \quad \mbox{ and } \quad 
    \widehat{A}^{\hat{x}_\lambda}_{\lambda,2}(t) =  \int_0^t \hat{x}_\lambda(s) \dd{s}.
\]
Recall that when $\alpha=0$, the weak convergence 
\eqref{eqn521} reads 
\[
\left(\lambda^{-\frac{1}{2}}  \left[ A^{\hat{x}_\lambda}_\lambda(t) -  \int_0^t [\lambda + \lambda^{\frac{1}{2}}\hat{x}_\lambda(s)]  \dd{s} \right] ,\int_0^t \hat{x}_\lambda(s) \dd{s} \right) \Rightarrow \left(\widetilde{B}, \int_0^\cdot \hat{x}(s)\dd{s}\right) \quad \mbox{in $\mathcal{D}[0,T]\times \mathcal{C}[0,T]$}.
\]
Therefore, when $0<\alpha<1$, 
\[
\left(\widehat{A}^{\hat{x}_\lambda}_{\lambda,1},\widehat{A}^{\hat{x}_\lambda}_{\lambda,2} \right) \Rightarrow \left(\BFzero, \int_0^\cdot \hat{x}(s)\dd{s}\right) \quad \mbox{in $\mathcal{D}[0,T]\times \mathcal{C}[0,T]$},
\]
where $\BFzero$ denotes the constant process in $\mathcal{C}[0,T]$ taking the value of $0$ for all $t\in[0,T]$. 
Then, using an argument similar to that for proving the convergence \eqref{eq:joint-converg-A1A2} in the case $\alpha=0$, we can show 
\[
\left( \widehat{A}_{\lambda,1}, \widehat{A}_{\lambda,2} \right)  \Rightarrow  \left( {\bf 0}, \int_0^\cdot U(s)\dd{s} \right) \quad \mbox{in $\mathcal{D}[0,T]\times\mathcal{C}[0,T]$},
\]
as $\lambda\to\infty$. 
Hence, it follows from the continuous-mapping theorem that if $\alpha\in(0,1)$, then 
\[
\widehat{A}_\lambda = \widehat{A}_{\lambda,1}+ \widehat{A}_{\lambda,2} \Rightarrow  \int_0^\cdot U(s) \dd{s} \quad\mbox{in } \mathcal{D}[0,T],
\] 
as $\lambda\to\infty$. 
This concludes the proof of Theorem~\ref{thm:FCLTofA}.
\Halmos
\endproof

\section{Proof of Theorem \ref{thm:QQ}}\label{sec:proof-scaled-num-in-sys}

With the weak convergence of the scaled arrival process $\widehat{A}_\lambda$, as stated in Theorem~\ref{thm:FCLTofA}, 
the weak convergence of the scaled number-in-system process $\widehat{Q}_\lambda$ in an infinite-server queueing system can be shown by verifying the conditions of the following result due to  Theorem 1 in Chapter 2.2 of \cite{Borovkov84_ec}. 

\begin{lemma}\label{lemma:Borovkov84}  
Consider an infinite-server queueing system indexed by $\lambda>0$. 
Let $A_\lambda(t)$ be the arrival process and let ${Q}_\lambda(t)$ be the number-in-system process. 
Let $F$ be the cumulative distribution function of the service times, and let $\bar F  \coloneqq 1-F$. 
Suppose the following conditions are satisfied: 
\begin{enumerate}[label=(\roman*)]
    \item There exists a random element $\widehat{A}_\infty$ in $\mathcal{D}[0,T]$ and     
    a positive, increasing sequence $\{B_\lambda:\lambda>0\}$ 
    with $B_\lambda\to \infty$ as $\lambda\to\infty$ such that 
    \[B_\lambda^{-1}[ A_\lambda(t) - \lambda t] \Rightarrow_{\mathcal{C}} \widehat{A}_\infty \quad\mbox{in } \mathcal{D}[0,T],\]    
    as $\lambda\to\infty$.
    \item The function $G(t)\coloneqq \int_0^t F(t-z)\dd{z}$ is H\"older continuous, that is, $|G(t_1)-G(t_2)|< c_1 |t_1-t_2|^{c_2}$ for some constants $c_1>0$ and $c_2>0$.
\end{enumerate}
Then, for $T>0$, there exists a random element $\widehat{Q}_\infty$ in $\mathcal{D}[0,T]$ such that 
\[ B_\lambda^{-1}\left[ Q_\lambda(t) -\lambda \int_0^t \bar{F}(t-z)\dd{z} \right] \Rightarrow_\mathcal{C} \widehat{Q}_\infty \quad \mbox{in } \mathcal{D}[0,T],\] 
as $\lambda\to\infty$, where 
\[
\widehat{Q}_\infty(t) = \xi(t)\mathbb{I}\{\alpha = 0\} + \int_0^t \bar{F}(t-z)\dd{\widehat{A}_\infty(z)}, 
\]
where $\xi$ is a centered Gaussian process, independent of $\widehat{A}_\infty$, having the covariance function 
\[
\E[\xi(t)\xi(t+u)] = \int_0^t F(t-z)\bar{F}(t+u-z)\dd{z}. 
\]
\end{lemma}

In Lemma~\ref{lemma:Borovkov84} , the $\mathcal{C}$-convergence, denoted as ``$\Rightarrow_{\mathcal{C}}$'', is defined as follows. 
Basically, it means a sequence of random elements  weakly converge in $\mathcal{D}[0,T]$ and the limit lies in $\mathcal{C}[0,T]$.
\smallskip

\begin{definition}\label{defofConv}
Fix $T>0$. 
For any $\lambda>0$, let $L_\lambda$ and $L_\dagger$ be random elements in the space $\mathcal{D}[0,T]$. 
Then, $L_\lambda$ is said to $\mathcal{C}$-converge to $L_\dagger$, denoted as $L_\lambda\Rightarrow_{\mathcal{C}} L_\dagger$, as $\lambda\to\infty$,  if 
\begin{enumerate}[label=(\roman*)]
    \item $L_\dagger\in\mathcal{C}[0,T]$ almost surely, and 
    \item $g(L_\lambda)\Rightarrow g(L_\dagger)$ as $\lambda\to\infty$, for all functionals $g:\mathcal{D}[0,T]\to \mathbb{R}$ which are continuous in $\mathcal{C}[0,T]$ with respect to the uniform metric.
\end{enumerate}
\end{definition}

\begin{lemma}\label{lemma:C-convergence}
Fix $T>0$. 
For any $\lambda>0$, let $L_\lambda$ and $L_\dagger$ be random elements in the space $\mathcal{D}[0,T]$. 
Then, $L_\lambda\Rightarrow_{\mathcal{C}} L_\dagger$ as $\lambda\to\infty$,  if and only if 
\begin{enumerate}[label=(\roman*)]
    \item $L_\dagger\in\mathcal{C}[0,T]$ almost surely, and 
    \item $L_\lambda \Rightarrow L_\dagger$ in $\mathcal{D}[0,T]$ as $\lambda\to\infty$.
\end{enumerate}
\end{lemma}

\proof{Proof.} 
For the ``if'' part: Fix an arbitrary functional $g:\mathcal{D}[0,T]\to \mathbb{R}$ which is continuous in $\mathcal{C}[0,T]$ with respect to the uniform metric.
Let $\mathrm{Disc}(g)$ be the set of discontinuity points of $g$ in $\mathcal{D}[0,T]$.
Then, 
$\mathrm{Disc}(g) \subseteq \mathcal{D}[0,T]\setminus \mathcal{C}[0,T]$.
Because $\pr(L_\dagger \in \mathcal{C}[0,1])=1$, 
we must have 
$\pr(L_\dagger \in \mathrm{Disc}(g)) = 0$. 
By the continuous-mapping theorem, 
$L_\lambda \Rightarrow L_\dagger$ implies that 
$g(L_\lambda) \Rightarrow g(L_\dagger)$, as $\lambda\to\infty$.

For the ``only if" part: 
Fix an arbitrary bounded and continuous functional $f:\mathcal{D}[0,T] \to\mathbb{R}$. 
By the definition of the weak convergence in $\mathcal{D}[0,T]$,
it suffices to show $\E[f(L_\lambda)] \to \E[f(L_\dagger)]$ as $\lambda\to\infty$.  

To this end, note that because $f$ is bounded, 
there exists a constant $M>0$ such that $|f|\leq M$. 
Therefore, using integration by parts, we have 
\begin{align*}
    \E[f(L_\lambda)]  = \int_0^M \pr( f(L_\lambda)> r) \dd{r} -  \int_{-M}^0 \pr(f(L_\lambda)\leq r) \dd{r}.
\end{align*}
It then follows from the bounded convergence theorem that 
\begin{align}
\lim_{\lambda\to\infty}\E[f(L_\lambda)]  = \int_0^M \lim_{\lambda\to\infty} \pr( f(L_\lambda)> r) \dd{r} -  \int_{-M}^0 \lim_{\lambda\to\infty} \pr(f(L_\lambda)\leq r) \dd{r}.
\label{eq:integ-by-parts}
\end{align}

Because $f$ is continuous in $\mathcal{D}[0,T]$, it must be continuous in $\mathcal{C}[0,T]$ with respect to the uniform metric. 
Therefore, by the definition of $\mathcal{C}$-convergence, 
$f(L_\lambda) \Rightarrow f(L_\dagger)$ as $\lambda\to\infty$, which implies 
\[
\lim_{\lambda\to\infty} \pr( f(L_\lambda)> r) = \pr( f(L_\dagger) > r)
\]
for any $r\in[-M,M]$. 
It then follows from Equation~\eqref{eq:integ-by-parts} that 
\begin{align*}
\lim_{\lambda\to\infty}\E[f(L_\lambda)] = 
\int_0^M \pr( f(L_\dagger)> r) \dd{r}  -  \int_{-M}^0 \pr(f(L_\dagger)\leq r) \dd{r} =  \E[f(L_\dagger)],
\end{align*}
which concludes the proof. 
\Halmos
\endproof
\smallskip

In light of Lemmas~\ref{lemma:Borovkov84} and \ref{lemma:C-convergence}, 
it is straightforward to prove Theorem~\ref{thm:QQ}. 
\smallskip

\proof{Proof of Theorem \ref{thm:QQ}.} 
Theorem \ref{thm:FCLTofA} states $\widehat{A}_\lambda \Rightarrow \widehat{A}_\infty$ in $\mathcal{D}[0,T]$ as $\lambda\to\infty$. 
Because $\widehat{A}_\infty \in \mathcal{C}[0,T]$ by definition, 
we know from Lemma~\ref{lemma:C-convergence}  that 
$\widehat{A}_\lambda \Rightarrow_{\mathcal{C}} \widehat{A}_\infty$. 
The proof is then concluded by applying Lemma~\ref{lemma:Borovkov84}. 
\Halmos
\endproof

\section{Asymptotic Normality of Self-normalized Arrival Counts} \label{sec:self-norm}

We present empirical evidence for the asymptotic normality of the ``unconventionally'' scaled arrival process. 
We first note that combining Theorems \ref{theo:CoD} and \ref{thm:FCLTofA} immediately implies  asymptotic normality of the ``self-normalized'' arrival processes. 

\begin{corollary}\label{coro}
If Assumptions~\ref{asp:Feller}--\ref{asp:alpha} hold, then 
for any $t>0$, we have that as $\lambda\to\infty$,
\[
 \frac{A_\lambda(t) - \E_\pi[A_\lambda(t)]}{\sqrt{\Var_\pi[A_\lambda(t)]}} \Rightarrow \mathsf{Normal}(0,1).
\]
\end{corollary}
\proof{Proof of Corollary~\ref{coro:Gauss}.} 
As per Equation~\eqref{Qinfty}, $Q_\infty(t)$ is a linear combination of three stochastic processes: $\int_0^t \bar{F}(t-z) U(z) \dd{z}$, $\int_0^t \bar{F}(t-z) \dd{\widetilde{B}(z)}$, and $\xi(t)$. 
Moreover, $U(t)$, $\widetilde{B}(z)$, and $\xi(t)$ are mutually independent by definition. 
Therefore, it suffices to show that all three are Gaussian processes.

Firstly, according to its definition in Equation~\eqref{OU}, 
$U(t)$ is an OU process with initial condition $U(0)=X_\dagger$.
Thus, if $X_\dagger$ is a constant or an independent normal random variable, $U(t)$ becomes a Gaussian process. 
Given that $\bar F$ is a deterministic function, following \citet[Lemma~6.2.3]{Oksendal03_ec}, we can assert that the Riemann integral $\int_0^t \bar{F}(t-z) U(z) \dd{z}$ constitutes a Gaussian process.

Secondly, since $\bar F$ is a deterministic function, the It\^{o} integral $\int_0^t \bar{F}(t-z) \dd{\widetilde{B}(z)}$ also forms a Gaussian process, as noted by \citet[Proposition~7.6]{Steel01_ec}.

Lastly, that $\xi(t)$ is a Gaussian process is a conclusion of Theorem~\ref{thm:QQ}. 
\Halmos 
\endproof

\begin{figure}[t]
\FIGURE{
    \includegraphics[width=0.33\textwidth]{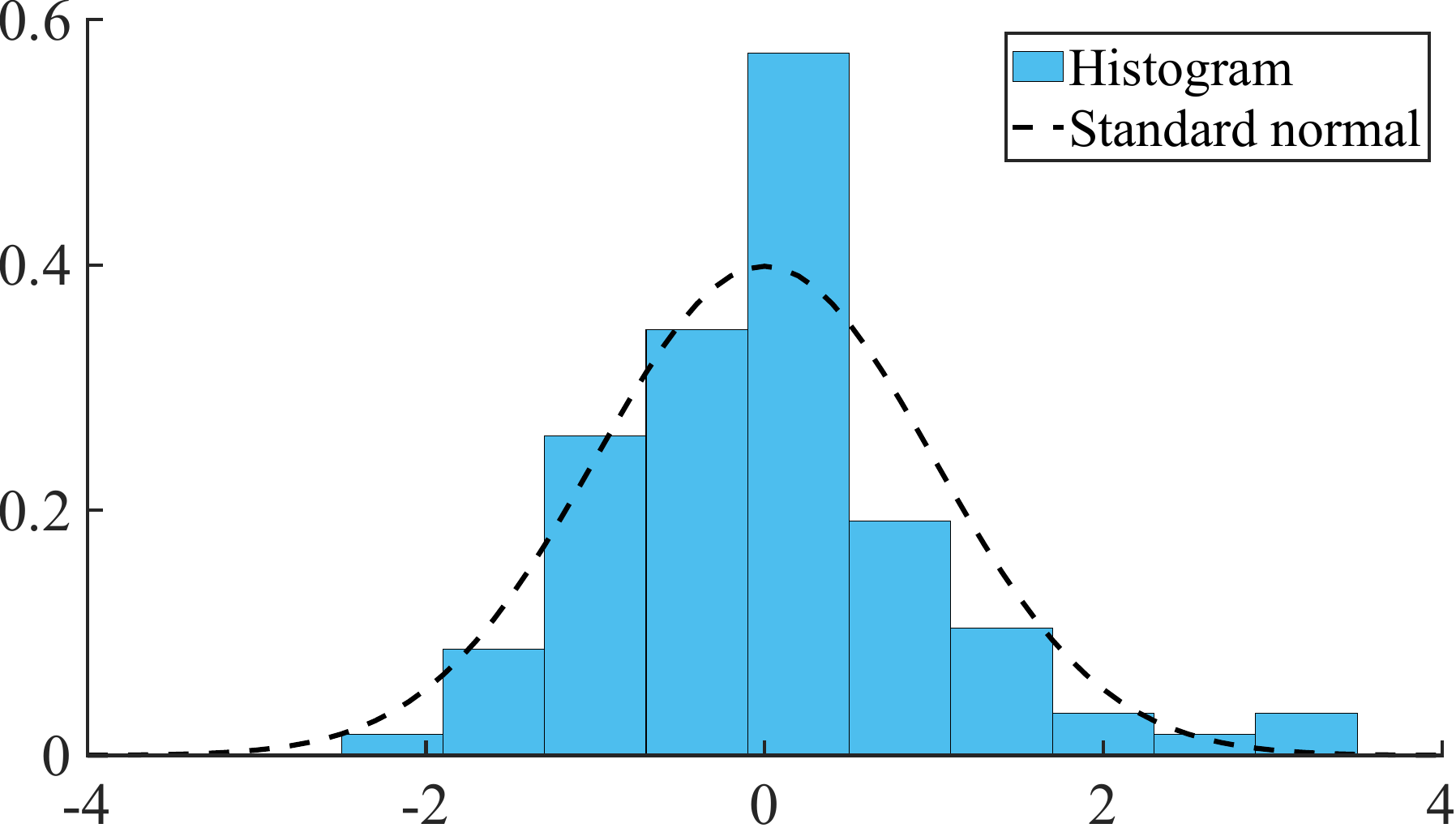}
	\includegraphics[width=0.33\textwidth]{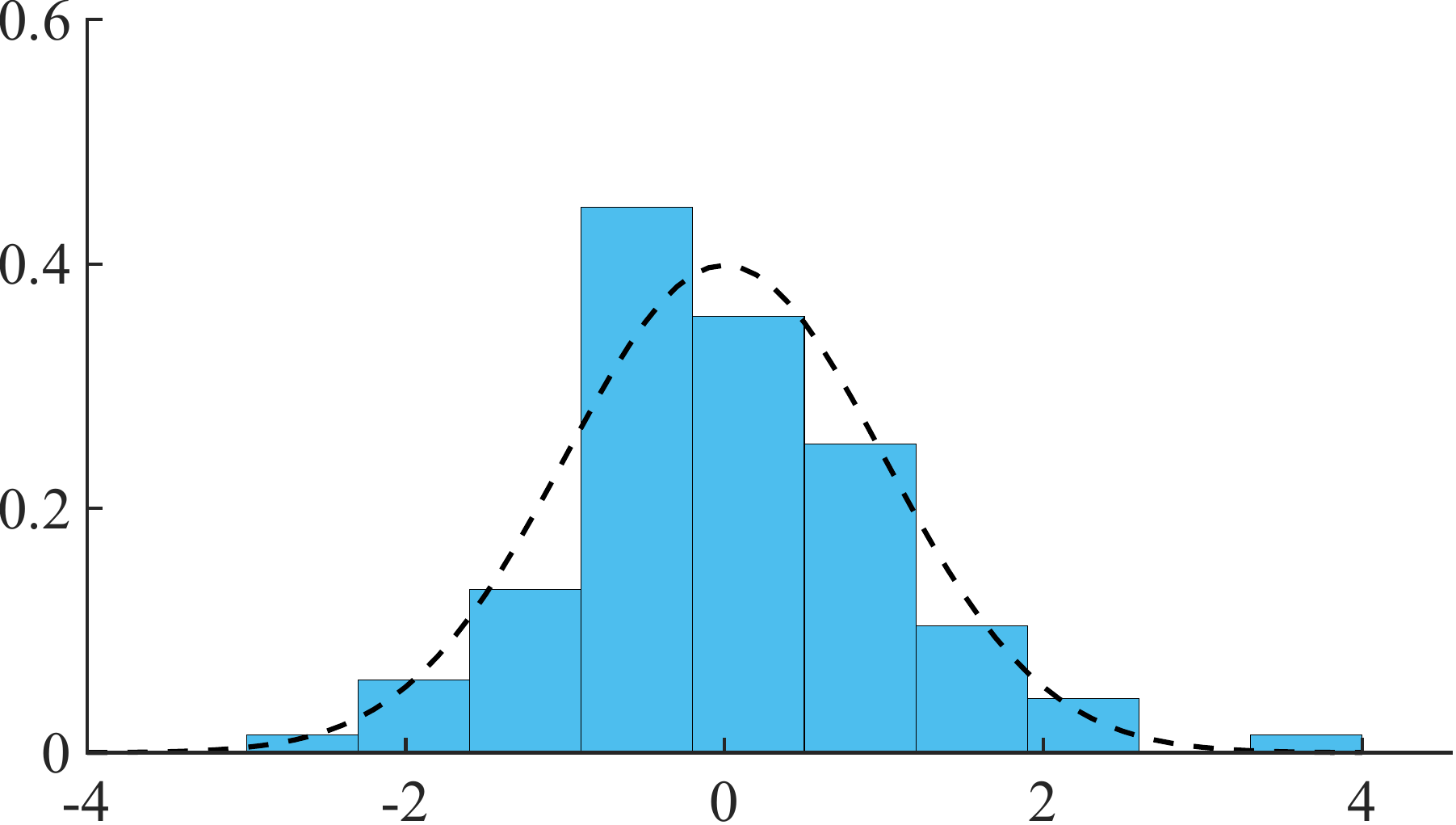}
	\includegraphics[width=0.33\textwidth]{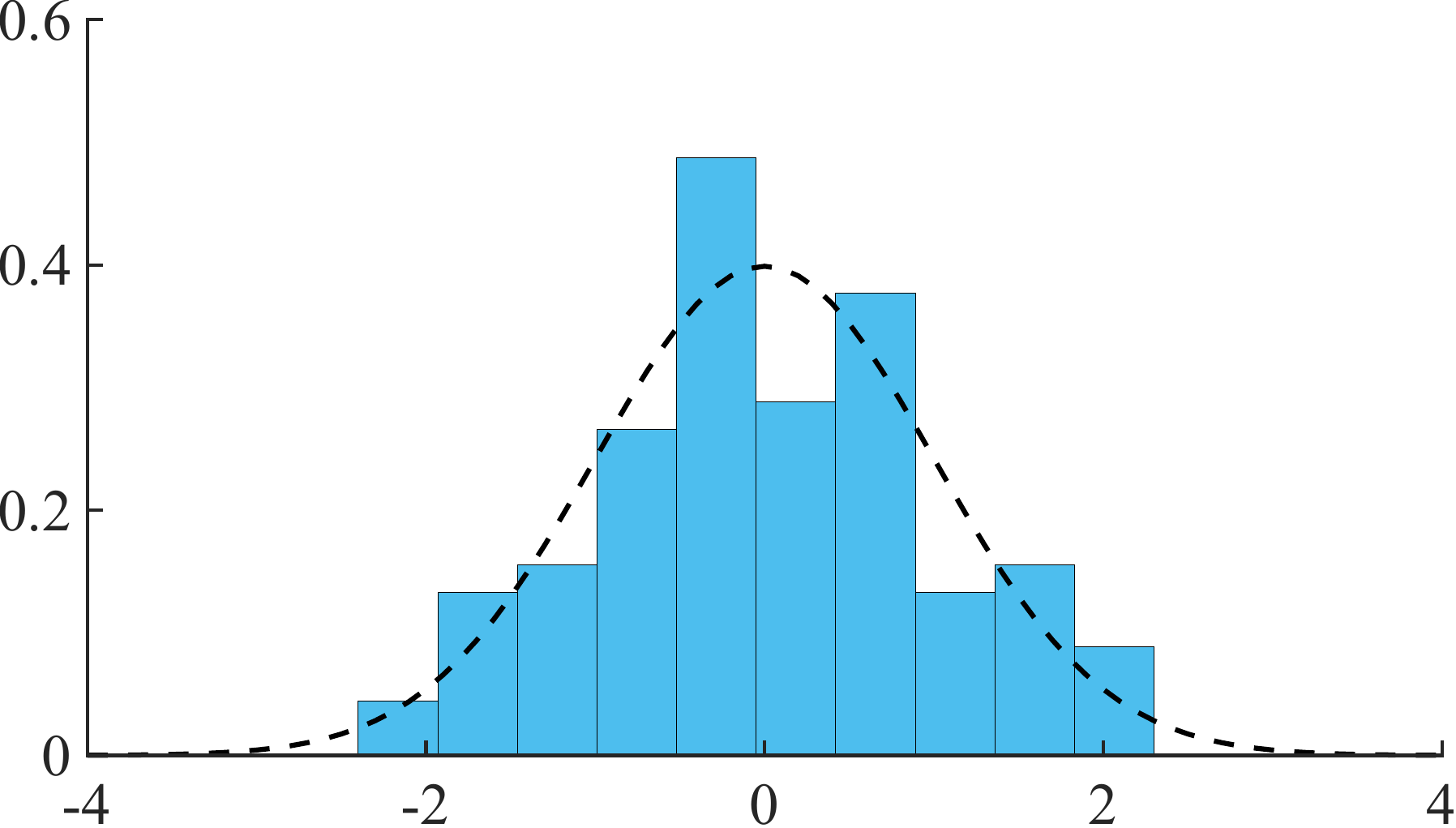}
}
{Distribution of Self-normalized Arrival Counts of the NYC 311 Call Center.\label{fig:kerneldensity}}
	{The arrival count $A_\lambda(t, t+\Delta]$ is observed across 96 days. The duration of the interval is $\Delta=\frac{1}{6}$ (hour). 
 Each of the three subfigures represents a distinct time of the day, with a corresponding arrival rate. Left: $t= 4$ AM and $\lambda=60$; Middle: $t=8$ AM and $\lambda=828$; Right: $t=12$ PM and $\lambda = 2052$. 
 }
\end{figure}

The asymptotic normality in Corollary~\ref{coro} is observed in the NYC 311 Call Center dataset. 
We assume the arrivals in each 10-minute time interval are stationary (similar to the treatment in Figure~\ref{fig:taylor}) and select three time intervals which begin at 4 AM, 8 AM, and 12 PM, respectively, and exhibit increasing mean numbers of arrivals. 
For each time interval, in Figure~\ref{fig:kerneldensity} we plot the distribution of the self-normalized arrival count across the 96 days in the dataset.

We stress here that one should not confuse the asymptotic normality shown in Figure~\ref{fig:kerneldensity}---which is a result of a ``heavy-traffic'' asymptotic regime---with the kind of asymptotic normality that stems from a ``large-sample'' asymptotic regime based on repeated observations of the same random variable (e.g., the arrival count in a given time interval). 
Heavy-traffic regimes are concerned with the behavior of the arrival process as the mean arrival rate $\lambda$ increases, while large-sample regimes are focused on the distribution of repeated observations of the same random variable.

The observation of asymptotic normality in the NYC 311 Call Center dataset supports the theoretical findings presented in this paper, indicating that the unconventionally scaled arrival process indeed converges to a Gaussian process in the heavy-traffic regime. 
This empirical evidence not only strengthens the theoretical underpinnings of our model for the arrival process, but also offers valuable insights for practical applications, such as staffing rules and other managerial decisions in call centers and similar service systems.

\section{Parameter Estimation for Non-stationary DSPP Model with Generalized CIR Intensity} \label{sec:MLE}

Let $N_{\Delta,1}, \ldots, N_{\Delta,k}$ denote the observed arrival counts in sequential time periods of length $\Delta$ under our non-stationary DSPP model with intensity defined in Equations~\eqref{eq:CIR_time-varying} and \eqref{eq:piecewise-constant-rate} in Section~\ref{sec:non-stat-DSPP}.

\subsection{Likelihood Approximation}

We extend the multivariate normal approximation of the likelihood of the arrival count vector $(N_{\Delta,1}, \ldots, N_{\Delta,k})$ from Section~\ref{sec:model-estimation} to the non-stationary setting.  

By definition, $N_{\Delta,i}$ is the arrival count from our stationary DSPP model with parameter $(\lambda_i, \alpha, \kappa, \sigma)$, 
for all $i=1,\ldots,k$. The heavy-traffic approximation of our stationary DSPP model in Theorem~\ref{thm:FCLTofA} implies that if $\lambda_i$ is sufficiently large, then
\begin{align*}
N_{\Delta,i} \stackrel{d}{\approx} {}& 
\lambda_i^{\frac{\alpha+1}{2}}\left(\int_{(i-1)\Delta}^{i\Delta} U(s) \dd{s} + \left(\widetilde{B}(i\Delta) - \widetilde{B}((i-1)\Delta) \right) -  \mathbb{I}\{\alpha = 0\} \right)+ \lambda_i \Delta,
\end{align*}
where $U$ is the OU process \eqref{OU} and $\widetilde{B}$ is a standard Brownian motion independent of $U$.

Thus, the joint distribution of $(N_{\Delta,1}, \ldots, N_{\Delta,k})$ can be approximated by a multivariate normal distribution with mean vector $\mathbf{m} = (\lambda_1\Delta,\ldots,\lambda_k\Delta)^\intercal$ and covariance matrix $\mathbf{\Sigma}$ with elements $\Sigma_{i,\ell}$ given by
\begin{align*}
\Sigma_{i,\ell}= \left\{
\begin{array}{ll}
     \displaystyle\lambda_i \Delta +  \frac{\sigma^2 \lambda_i^{\alpha+1} \Delta}{\kappa^2} \left( 1- \frac{1-e^{-\kappa \Delta}}{\kappa \Delta} \right),&  \mbox{if } i = \ell, \\[2ex] 
    \displaystyle (\lambda_i\lambda_\ell)^{\frac{\alpha+1}{2}}\frac{\sigma^2}{2\kappa^3} \left[ e^{-\kappa (\ell-1) \Delta} - e^{-\kappa \ell \Delta} \right] \left[ e^{\kappa i \Delta} - e^{\kappa (i-1) \Delta} + e^{-\kappa i \Delta} - e^{-\kappa (i-1) \Delta} \right], & \mbox{if } i < \ell, \\[2ex] 
    \displaystyle (\lambda_i\lambda_\ell)^{\frac{\alpha+1}{2}} \frac{\sigma^2}{2\kappa^3} \left[ e^{-\kappa (i-1) \Delta} - e^{-\kappa i \Delta} \right] \left[ e^{\kappa \ell \Delta} - e^{\kappa (\ell-1) \Delta} + e^{-\kappa \ell \Delta} - e^{-\kappa (\ell-1) \Delta} \right], & \mbox{if } i > \ell. 
\end{array}
\right.
\end{align*}  

Using this multivariate normal approximation of the likelihood, we can efficiently solve the MLE approximately and compute the corresponding AIC and BIC for model selection.

\subsection{Estimating $\lambda_i$ with SIPP}\label{sec:SIPP}

For large numbers of segments (e.g., $k=48$ in Section~\ref{sec:casestudy}), solving the MLE for models $\mathcal{M}_2$ through $\mathcal{M}_5$ becomes a challenging high-dimensional optimization problem. We address this through a two-step estimation approach: first estimating each $\lambda_i$ via SIPP using the sample mean of observed arrival counts in segment $i$ to estimate $\lambda_i \Delta$, then maximizing the log-likelihood function with respect to the remaining parameters such as $(\alpha, \kappa, \sigma)$ for $\mathcal{M}_5$ while holding these $\lambda_i$ estimates fixed. This reduces the problem to a more tractable low-dimensional optimization.

Suppose we observe $m$ cycles of a periodic arrival process (e.g., $m$ business days of call arrivals). Let $N_{\Delta, i, j}$ denote the arrival count in segment $i$ of length $\Delta$ during cycle $j$, for $i=1,\ldots,k$ and $j=1,\ldots,m$. The basic SIPP estimates $\lambda_i$ as
\[
\widehat{\lambda}_i = \frac{1}{m \Delta } \sum_{j=1}^m N_{\Delta, i, j}. 
\]
For the three SIPP variants---\textsf{SIPP Avg}, \textsf{SIPP Min}, and \textsf{SIPP Max}---we divide each segment into $h$ equal sub-intervals of length $\Delta/h$ (e.g., $h=3$ yielding 10-minute sub-intervals in Section~\ref{sec:casestudy}). Let $N_{\Delta, i, j}^{(l)}$ denote the arrival count in sub-interval $l$ of segment $i$ during cycle $j$. The mean arrival rate for each sub-interval is estimated as
\[
\widehat{\lambda}_{i}^{(l)} = \frac{h}{m \Delta} \sum_{j=1}^m N_{\Delta, i,j}^{(l)}. 
\]
The three variants are then defined as
\[
\widehat\lambda_i = 
\left\{\begin{array}{ll}
   \displaystyle \frac{1}{h}\sum_{h=1}^h \widehat{\lambda}_{i}^{(l)}, & \mbox{ for \textsf{SIPP Avg}}, \\[2ex]
   \displaystyle \min_{1\leq l \leq h}\widehat{\lambda}_{i}^{(l)},  & \mbox{ for \textsf{SIPP Min}}, \\[2ex] 
   \displaystyle \max_{1\leq l \leq h} \widehat{\lambda}_{i}^{(l)}, & \mbox{ for \textsf{SIPP Max}}.
\end{array}
\right.
\]

The \textsf{SIPP Mix} approach combines these three SIPP variants. 
For the NYC 311 Call Center data, it uses \textsf{SIPP Min} from 6 AM--9 AM when arrival rates increase rapidly, \textsf{SIPP Max} from 4 PM--8 PM when rates decrease rapidly, and \textsf{SIPP Avg} for all other periods.

\section{Analysis of Model $\mathcal{M}_3$}\label{StaticDSPP}

In this section, we derive a staffing rule under model $\mathcal{M}_3$ introduced in Section~\ref{sec:model-comparison}. 
To well pose the model, we fix a complete probability space $(\Omega, \mathscr{F}, \pr)$ and a right-continuous, complete information filtration $\mathbb{F}=\{\mathscr{F}_t: t\geq 0\}$, as in Section~\ref{themodel}. 
The arrival process $A^{\mathsf{S}}=\{A^{\mathsf{S}}(t):t\geq 0\}$ is modeled as a doubly stochastic Poisson process with intensity  
\begin{equation}\label{RVarrivalrate}
X^{\mathsf{S}} = \lambda + \lambda^{\frac{\alpha+1}{2}} Y, 
\end{equation}
where $0\leq \alpha \leq 1$ and $Y$ is a $\mathscr{F}_0$-measurable random variable with zero mean and finite variance $\sigma_Y^2$. 
The superscript $\mathsf{S}$ denotes ``static,'' indicating that the arrival rate is a random variable independent of time $t$ rather than a stochastic process. 

While this arrival model has been used in \cite{Maman09,BassambooRandhawaZeevi10,KocagaArmonyWard15} and \cite{HuChanDong25} for staffing service systems, these studies did not derive staffing rules through heavy-traffic analysis aimed at achieving a target delay probability, which we will do here.

\subsection{Asymptotic Normality}\label{sec:FCLTofA_ec}

We utilize the notation 
\[\widehat{A}^{\mathsf{S}}_\lambda(t) :=  \lambda^{-\frac{\alpha+1}2}[A^{\mathsf{S}}(t)-\lambda t]\] 
to denote the scaled version of the arrival process $A^{\mathsf{S}}(t)$. Notice that, we continue to use the superscript $(s)$ to denote that the intensity process follows Equation~\eqref{RVarrivalrate}. 

\begin{theorem}\label{thm:FCLTofA_ec}
Suppose $0\leq \alpha\leq 1$. Then, for any $T>0$, we have  
\begin{equation}\label{thm:FCLT:eq1_ec}
\widehat{A}^{\mathsf{S}}_\lambda \Rightarrow \widehat{A}^{\mathsf{S}}_\infty\quad \mbox{in $\mathcal{D}[0,T]$}
\end{equation}
as $\lambda\rightarrow\infty$,
where 
\begin{equation}\label{FCLTofA_ec}
\widehat{A}^{\mathsf{S}}_\infty(t) = Yt + \widetilde{B}(t) \mathbb{I}\{\alpha = 0\},
\end{equation}
and $\widetilde{B}$ is a standard Brownian motion.
\end{theorem}

\proof{Proof.}
By following a similar approach to the proof of Theorem~\ref{thm:FCLTofA} in Section~\ref{EC13}, we can establish the proof for Theorem~\ref{thm:FCLTofA_ec}. Notably, we can directly conclude that the scaled intensity process 
\[\widehat{X}^{\mathsf{S}}_\lambda(t):=\lambda^{-(\alpha+1)/2}[X^{\mathsf{S}}(t)-\lambda] = Y\] 
holds for all $\lambda$ and $t$, without the need to prove a weak convergence result as presented in Theorem~\ref{FCLTofintensity}. By substituting this clear conclusion into the proof of Theorem~\ref{thm:FCLTofA}, we can successfully establish the proof of Theorem~\ref{thm:FCLTofA_ec} directly and then obtain the result~\eqref{thm:FCLT:eq1_ec}--\eqref{FCLTofA_ec}.
\Halmos
\endproof

\subsection{Alpha Safety Rules}

We consider an infinite-server queue system, and utilize the notation 
\[\widehat{Q}^{\mathsf{S}}_\lambda(t) :=\lambda^{-\frac{\alpha+1}{2}}\left[Q^{\mathsf{S}}_\lambda(t)-\lambda\int_0^t\bar{F}(t-z)\dd z\right]
\]
to denote the scaled version of the number of customers process. Notice that, the superscript $(s)$ denotes that the arrival process follows $A^{\mathsf{S}}(t)$. 

\begin{theorem}\label{thm:QQ_ec}
Suppose that $0\leq \alpha\leq 1$, and the function $G(t)\coloneqq \int_0^t F(t-z)\dd{z}$ is H\"older continuous, that is, $|G(t_1)-G(t_2)|< c_1 |t_1-t_2|^{c_2}$ for some constants $c_1>0$ and $c_2>0$. Let $\widehat{A}^{\mathsf{S}}_\infty$ be the process defined in Equation~\eqref{FCLTofA_ec}, and let
$\xi$ be an independent, zero-mean Gaussian process with the following covariance function: 
\begin{equation*}\label{eq:xi-cov_ec}
\E[\xi(t)\xi(t+u)] = \int_0^t F(t-z)\bar{F}(t+u-z)\dd{z},    
\end{equation*}
for all $t, u\geq 0$.
Then, for any $T>0$, we have
\begin{equation*}\label{eqn:eq13_ec}
\widehat{Q}^{\mathsf{S}}_\lambda  \Rightarrow  \widehat{Q}^{\mathsf{S}}_\infty  \quad \mbox{in $\mathcal{D}[0,T]$} 
\end{equation*}
as $\lambda\rightarrow\infty$, where 
\begin{equation}\label{limitofQ_ec}
\widehat{Q}^{\mathsf{S}}_\infty(t) = \int_0^t \bar{F}(t-z)\dd{\widehat{A}^{\mathsf{S}}_\infty(z)} + \xi(t)\mathbb{I}\{\alpha = 0\}.  
\end{equation}
\end{theorem}

\proof{Proof} 
Theorem \ref{thm:FCLTofA_ec} states $\widehat{A}^{\mathsf{S}}_\lambda \Rightarrow \widehat{A}^{\mathsf{S}}_\infty$ in $\mathcal{D}[0,T]$ as $\lambda\to\infty$. 
Because $\widehat{A}^{\mathsf{S}}_\infty \in \mathcal{C}[0,T]$ by its definition, 
we know from Lemma~\ref{lemma:C-convergence}  that 
$\widehat{A}^{\mathsf{S}}_\lambda \Rightarrow_{\mathcal{C}} \widehat{A}^{\mathsf{S}}_\infty$. 
The proof is then concluded by applying Lemma~\ref{lemma:Borovkov84}.~\Halmos
\endproof
\smallskip

Applying Equation~\eqref{FCLTofA_ec} to Equation~\eqref{limitofQ_ec}, we have 
\begin{equation*}\label{Qinfty_ec}
\widehat{Q}^{\mathsf{S}}_\infty(t) = Y \int_0^t \bar{F}(t-z) \dd{z} + \biggl[\int_0^t \bar{F}(t-z) \dd{\widetilde{B}(z)} + \xi(t)\biggr] \mathbb{I}\{\alpha = 0\}.
\end{equation*}
Suppose $Y$ is a normal random variable, i.e., $Y\sim\mathsf{Normal}(0,\sigma_Y^2)$, then $\widehat{Q}^{\mathsf{S}}_\infty(t)$ is a Gaussian process. Following the discussion between Equations~\eqref{eq:Q-approx-dist-0}--\eqref{eqn:eq140}, we have 
\begin{equation}\label{eqn:eq140_ec}
Q^{\mathsf{S}}_\lambda(t) \stackrel{d}{\approx}  
\mathsf{Normal} \left( \frac{\lambda}{\mu}, \, \lambda^{\alpha+1} \biggl(  V_1^{\mathsf{S}}(\infty) + \frac{1}{\mu}   \mathbb{I}\{\alpha = 0\} \biggr) \right),
\end{equation}
for large values of $\lambda$ and $t$, where $V_1^{\mathsf{S}}(\infty) = \lim_{t\to\infty} \Var[Y \int_0^t\bar{F}(t-z)\dd z]$. The approximation~\eqref{eqn:eq140_ec} offers a staffing rule for large-scale service systems. To satisfy the target delay probability, i.e., $\pr(Q^{\mathsf{S}}_\lambda(t)>n)\approx \varepsilon$, we can configure the number of servers to be the $(1-\varepsilon)$-quantile of the normal distribution outlined in the approximation~\eqref{eqn:eq140_ec}.
More specifically, let $\Phi$ represent the CDF of the standard normal distribution and define $\beta = \Phi^{-1}(1-\varepsilon)$. 
We propose a staffing rule as:
\begin{equation}\label{eq:basic-alpha_ec}
n^{\mathsf{S}}  =   \frac{\lambda}{\mu} + \beta   \lambda^{\frac{\alpha +1}{2}} \sqrt{  V_1^{\mathsf{S}}(\infty) +  \frac{1}{\mu}   \mathbb{I}\{\alpha = 0\} }.
\end{equation}
We refer to this staffing rule $n^{\mathsf{S}}$ as the basic alpha safety rule under the setting of the random variable arrival rate~\eqref{RVarrivalrate}. Notice that, the rule $n^{\mathsf{S}}$ has the same form of the rule $n^\ast_{\mathsf{Basic}}$ given in Equation~\eqref{eq:basic-alpha} with the only difference being the term $V^{\mathsf{S}}_1(\infty)$ therein. %

\subsection{Model Estimation}

The parameter estimation method based on heavy-traffic approximation, developed in Appendix~\ref{sec:MLE} for our non-stationary DSPP model, can also be applied to $\mathcal{M}_3$ (which incorporates $\mathcal{M}_2$ as a special case) and its non-stationary version.

As established by Theorem~\ref{thm:FCLTofA_ec} and its proof in Appendix~\ref{sec:FCLTofA_ec}, the joint distribution of this arrival count vector $(N_{\Delta,1}, \ldots, N_{\Delta,k})$ can be approximated by a multivariate normal distribution with mean vector $\mathbf{m} = (\lambda_1\Delta,\ldots,\lambda_k\Delta)^\intercal$ and covariance matrix $\mathbf{\Sigma}$ with elements $\Sigma_{i,\ell}$ given by
\begin{align*}
\Sigma_{i,\ell}
=\left\{
\begin{array}{ll}
    \displaystyle  \lambda_i \Delta +  \lambda_i^{\alpha+1} \sigma_Y^2 \Delta^2, & \mbox{if } i = \ell, \\[2ex]
     \displaystyle (\lambda_i\lambda_\ell)^{\frac{\alpha+1}{2}} \sigma_Y^2 \Delta ^2, & \mbox{if } i \neq \ell. 
\end{array}
\right.
\end{align*}
Using this multivariate normal approximation, we can compute the likelihood function without dealing with high-dimensional numerical integration arising from unobservable intensities. This enables efficient computation of MLE estimates for parameters $(\alpha, \sigma_Y)$.

\section{Additional Numerical Experiments}\label{sec:additional-exp}

In this section, we present additional results omitted in Sections~\ref{sec:synthetic}--\ref{sec:casestudy} and 
conduct supplementary numerical experiments to validate the robustness of our proposed staffing rules.

\subsection{Additional Results for Section~\ref{sec:synthetic}}\label{sec:add-res-Sec6}

\begin{figure}[ht]
\FIGURE{    \includegraphics[width=0.5\textwidth]{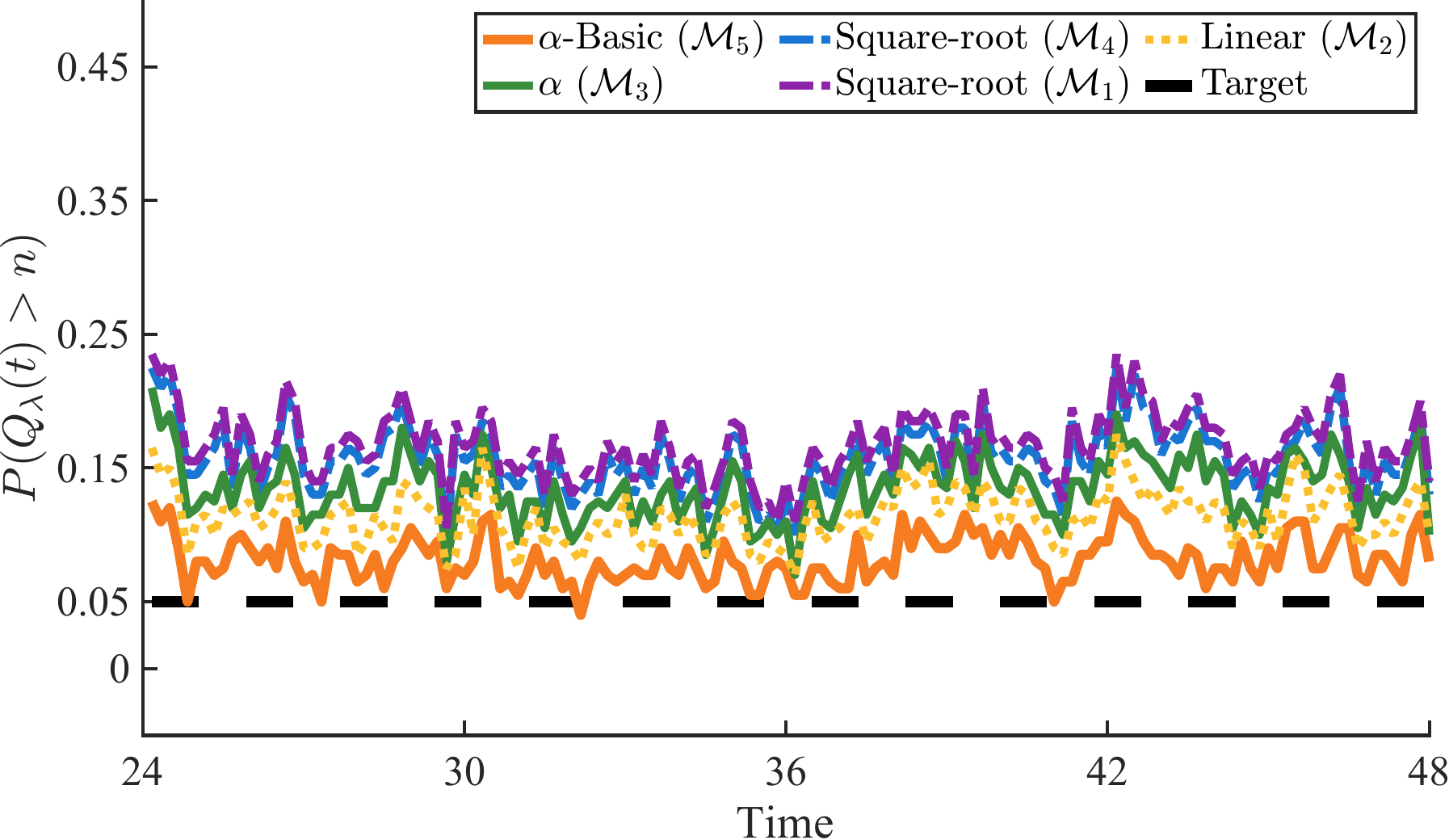}
	\includegraphics[width=0.5\textwidth]{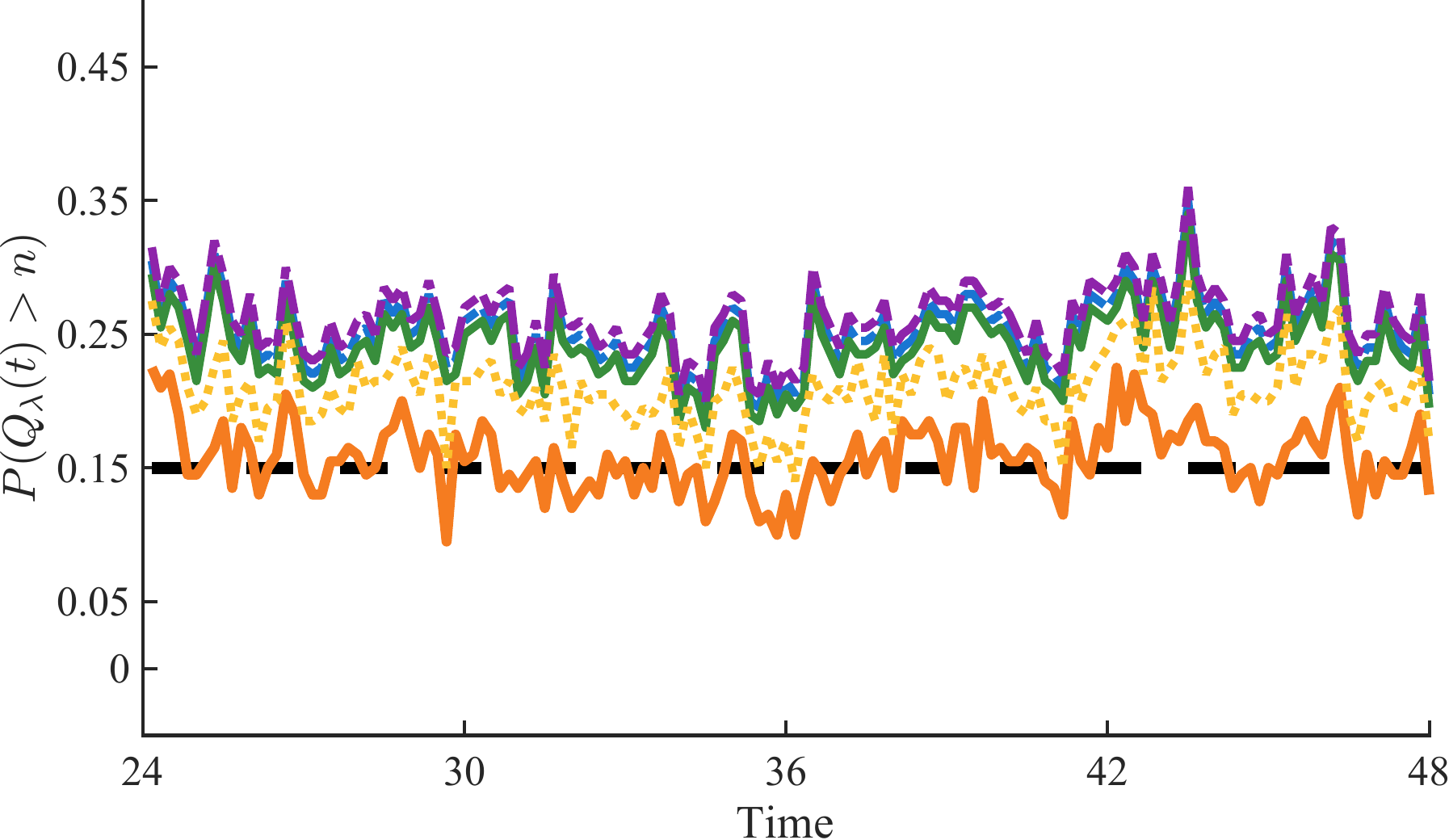}}
{Performance of Staffing Rules in Infinite-Server Systems with Arrival Model $\mathcal{M}_5$  ($\lambda=150$). \label{fig:sim-inf-small-lambda}}
{}

\FIGURE{    \includegraphics[width=0.5\textwidth]{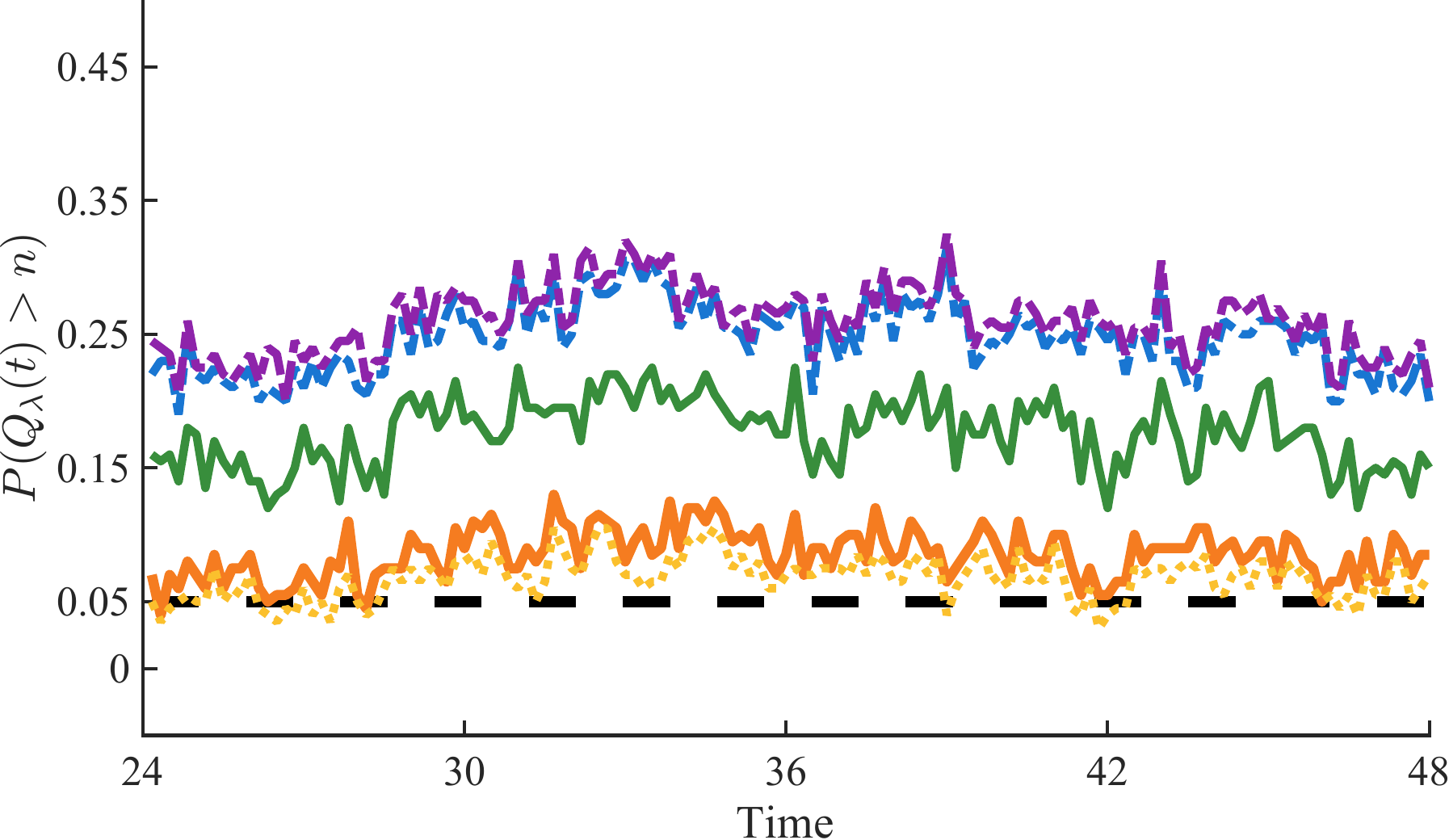}
\includegraphics[width=0.5\textwidth]{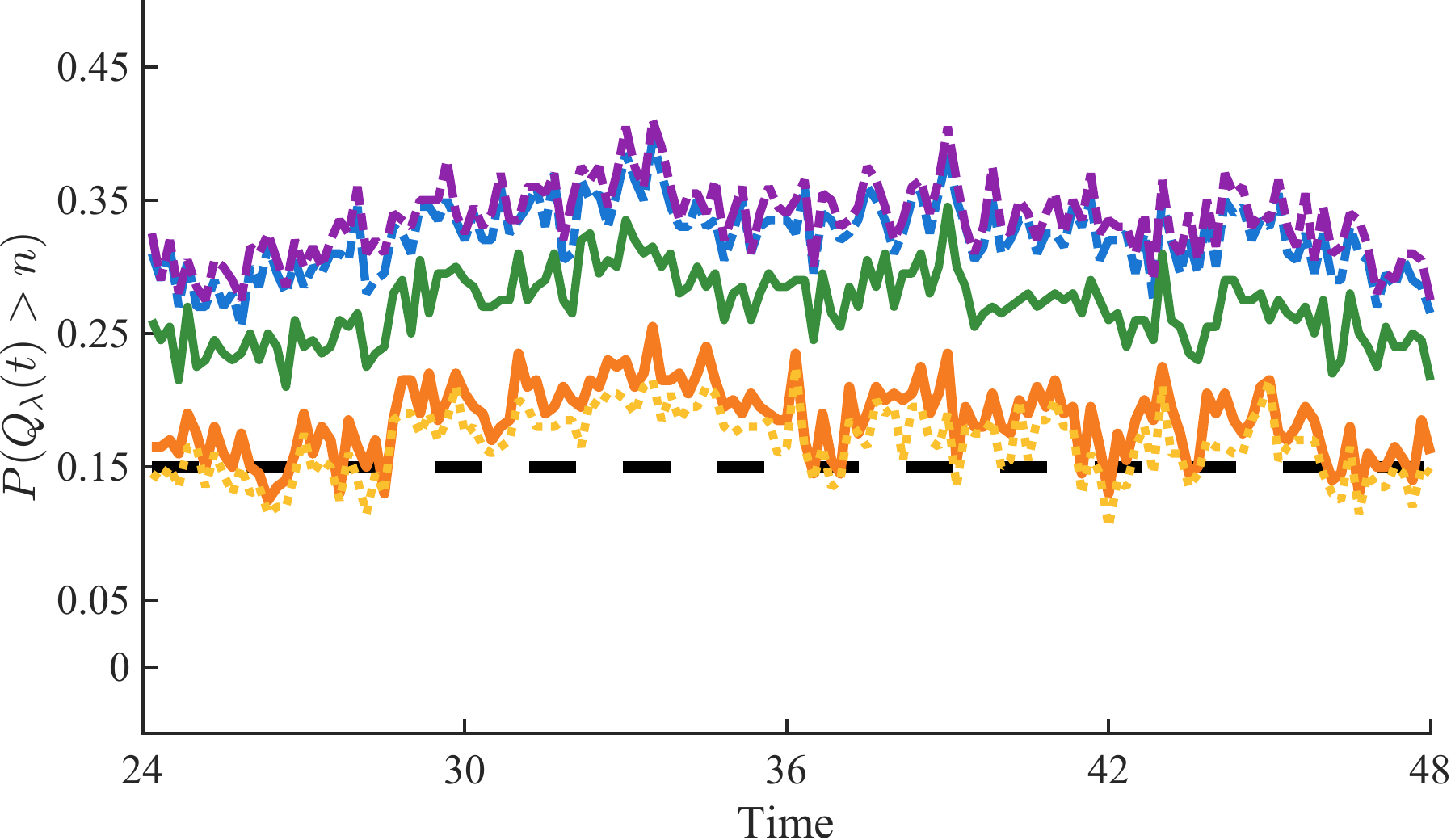}}
{Performance of Staffing Rules in Infinite-Server Systems with Arrival Model $\mathcal{M}_5$  ($\lambda=600$). \label{fig:sim-inf-med-lambda}}
{The probability that the number of customers in an infinite-server system exceeds $n$, where $n$ varies depending on the staffing rule. Results exclude the first 24 hours (warm-up period).
Left: $\varepsilon = 0.05$. Right: $\varepsilon = 0.15$.}
\end{figure}

\begin{figure}[ht]
\FIGURE{ 
    \includegraphics[width=0.5\textwidth]{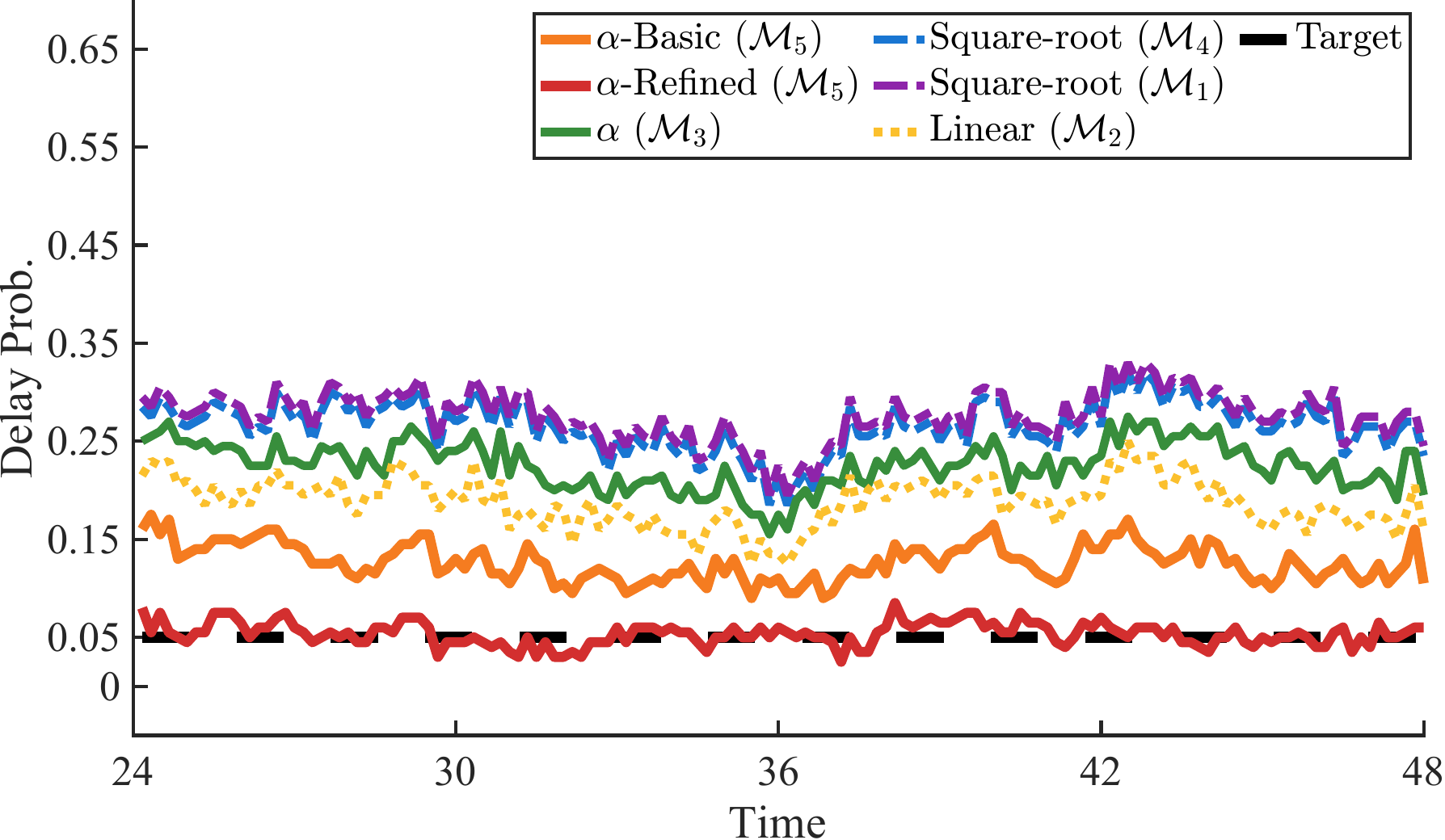}
	\includegraphics[width=0.5\textwidth]{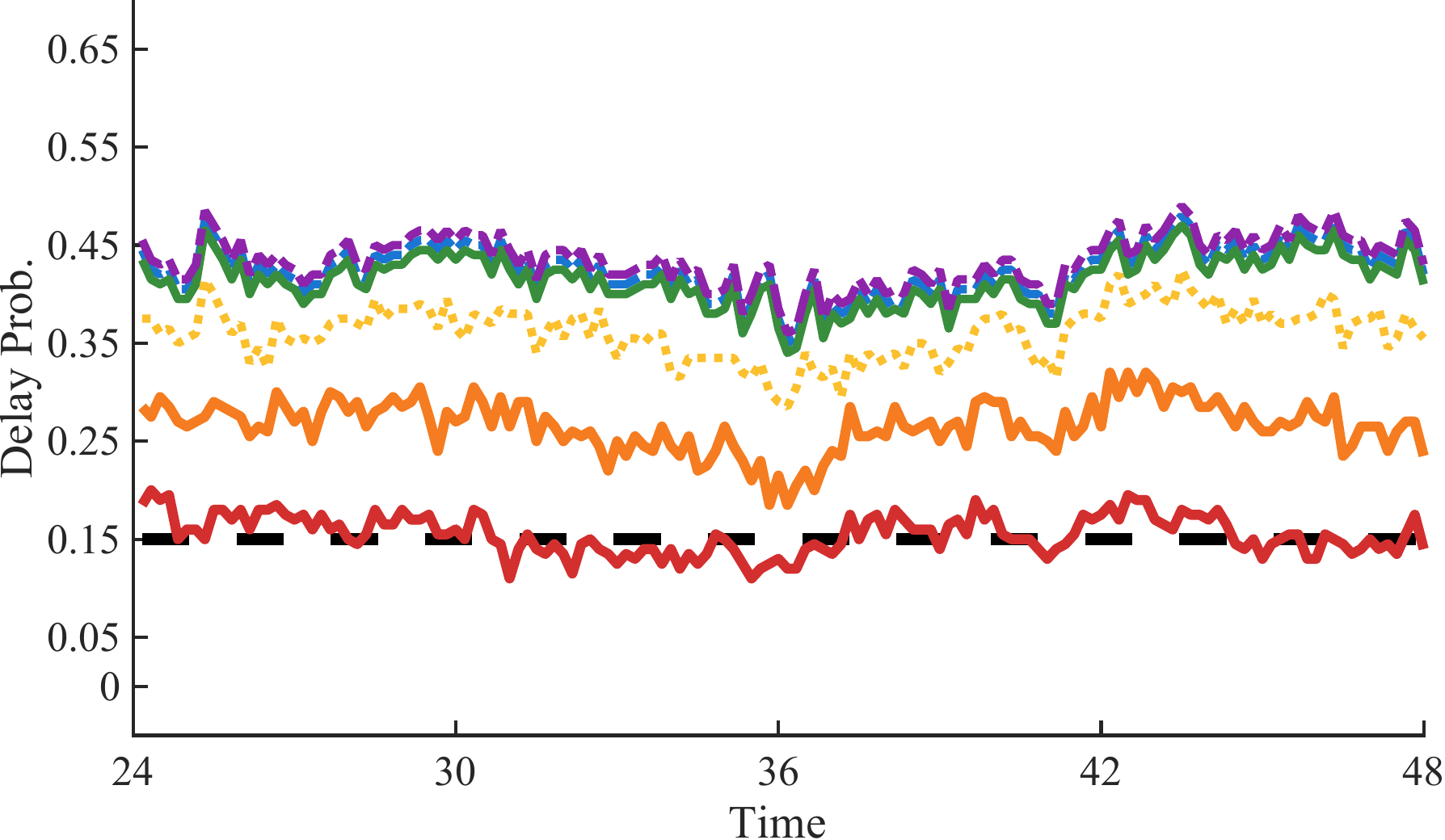}
	}
{Performance of Staffing Rules in Finite-server Systems with Arrival Model $\mathcal{M}_5$  ($\lambda=150$). \label{figsim120}}
{}

\FIGURE{
    \includegraphics[width=0.5\textwidth]{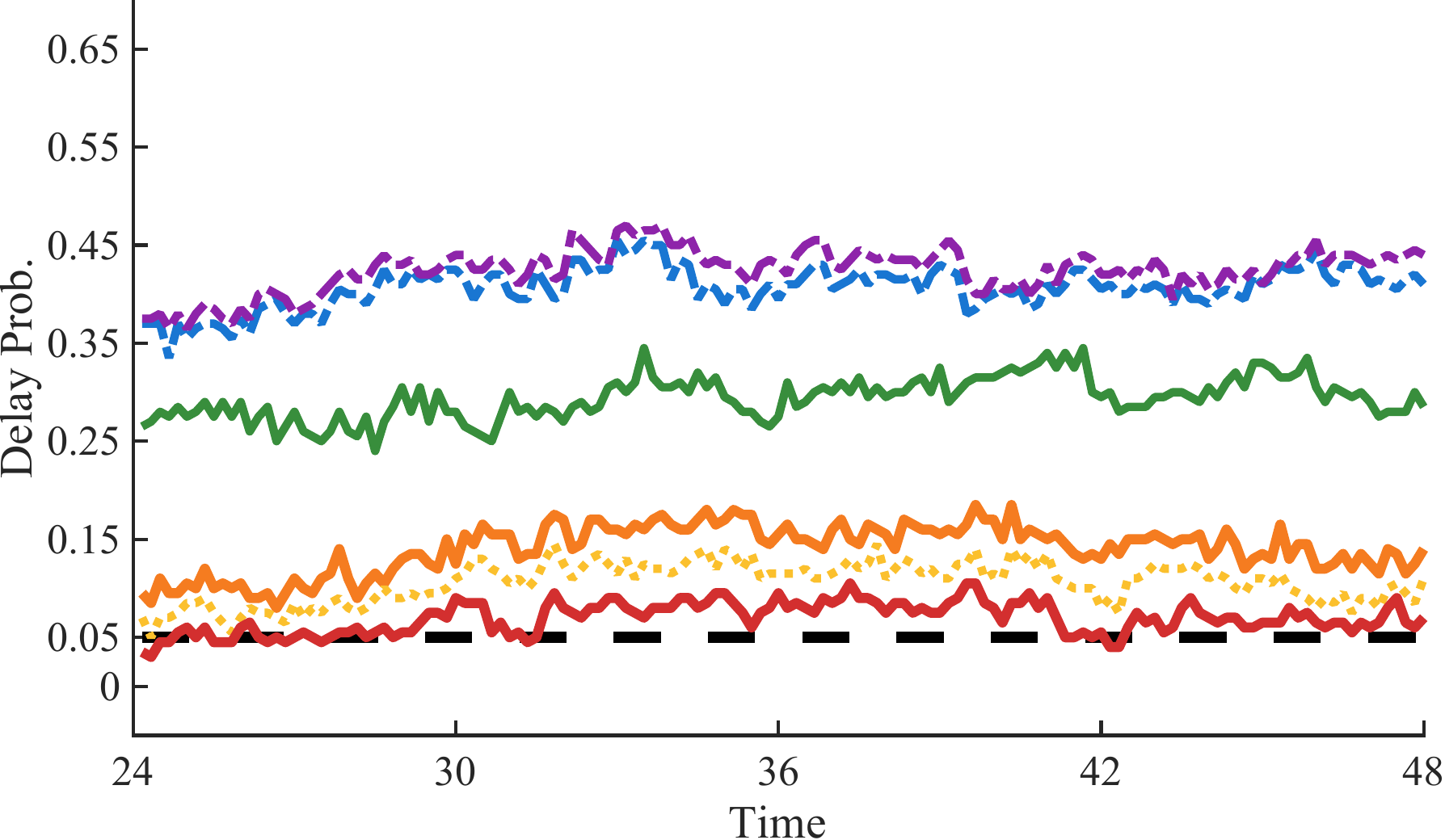}
	\includegraphics[width=0.5\textwidth]{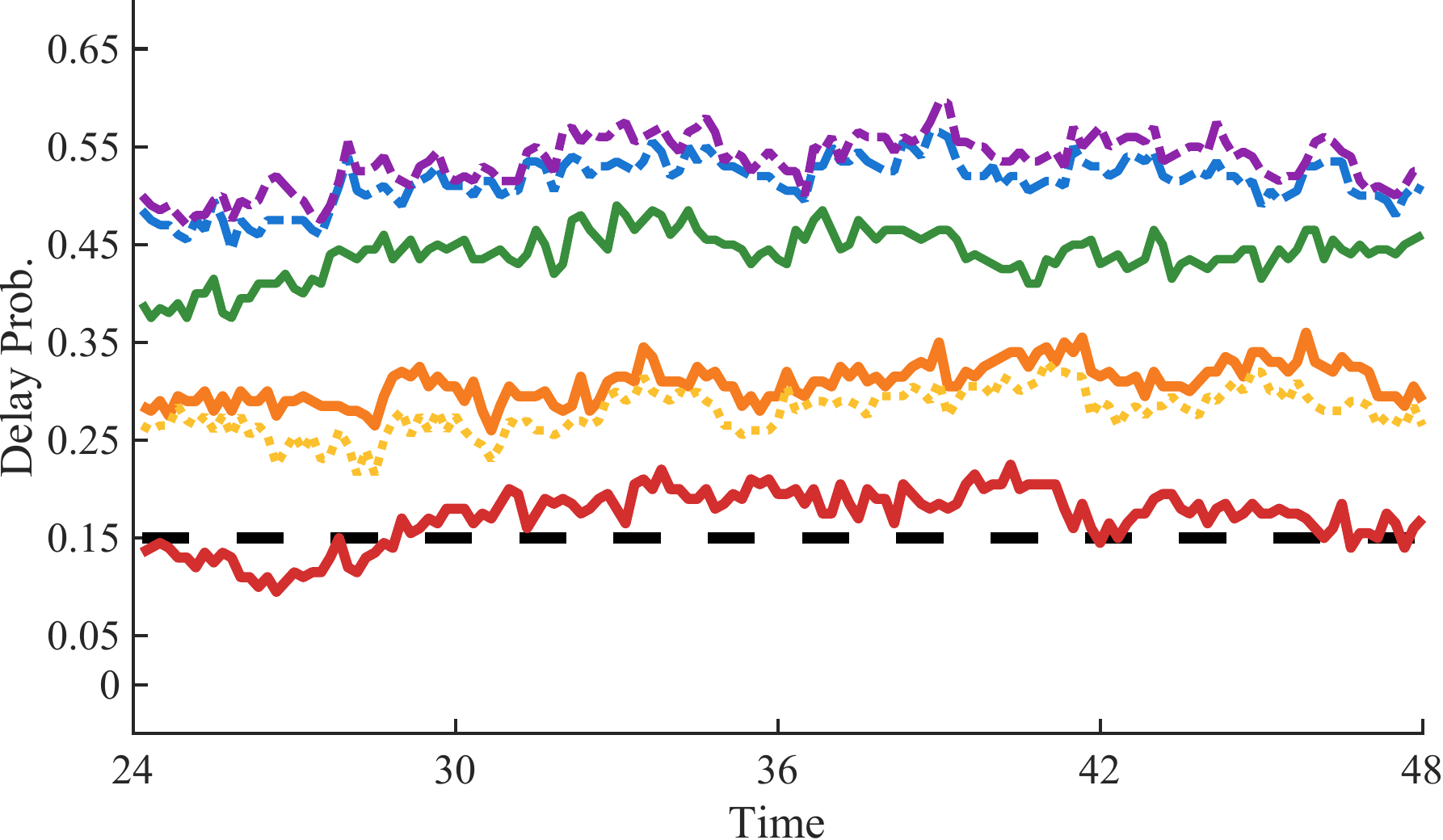}
	}
{Performance of Staffing Rules in Finite-server Systems with Arrival Model $\mathcal{M}_5$  ($\lambda=600$). \label{fig:sim-finite-server-med}}
{The delay probability as defined by 
$\pr( Q_\lambda^{(n)}(t) > n)$, where $n$ varies depending on the staffing rule. 
Results exclude the first 24 hours (warm-up period).
Left: $\varepsilon = 0.05$. Right: $\varepsilon = 0.15$. }
\end{figure}

Figure~\ref{fig:sim-inf-small-lambda} and Figure~\ref{fig:sim-inf-med-lambda} present the performance of staffing rules in infinite-server systems under our DSPP model for $\lambda=150$ and $600$, respectively. Figure~\ref{figsim120} and Figure~\ref{fig:sim-finite-server-med} present the corresponding results for finite-server systems.

\subsection{Model Misspecification} \label{sec:model-misspec}

In addition to the time-of-day effect and over-dispersion, 
service system arrivals typically exhibit temporal correlation across different intervals within a day. 
\cite{Avramidis04_ec} introduced three models to capture the complex correlation structure of arrival counts, 
primarily for simulation analysis of queueing systems rather than staffing rule development. 
Their Model 3, which offers superior flexibility and fitting performance, serves as our ground truth for comparing staffing rules under $\mathcal{M}_1$ through $\mathcal{M}_5$. This comparison demonstrates the robustness of our non-stationary DSPP model in Section~\ref{sec:non-stat-DSPP} and the refined alpha safety rule to model misspecification.

We refer to Model 3 of \cite{Avramidis04_ec} as the ADL model, based on the authors' last name initials. Let $Z_i$ denote the arrival count in segment $i$ for each $i=1,\ldots,k$, and let $Z = \sum_{i=1}^k Z_i$ represent the total arrivals over one period $T$. The ADL model specifies the joint distribution of $(Z_1,\ldots,Z_k)$ through a two-step process: first, it models the total count $Z$ using a gamma distribution; then, it introduces an independent $k$-dimensional random vector $\mathbf{p} \coloneqq (p_1,\ldots,p_k)$, where $\sum_{i=1}^k p_i = 1$, representing the proportion of $Z$ allocated to each $Z_i$. To generate $(Z_1,\ldots,Z_k)$, the model draws $\mathbf{p}$ from a Dirichlet distribution with parameters $(a_1,\ldots,a_k)$ and sets $Z_i = [p_i Z]$, where $[\cdot]$ denotes rounding. 
After obtaining the arrival counts, the model generates arrival times by uniformly distributing the $Z_i$ arrivals within each segment $[(i-1)\Delta, i\Delta)$.

\begin{figure}[t]
  \FIGURE{
    \includegraphics[width=0.5\textwidth]{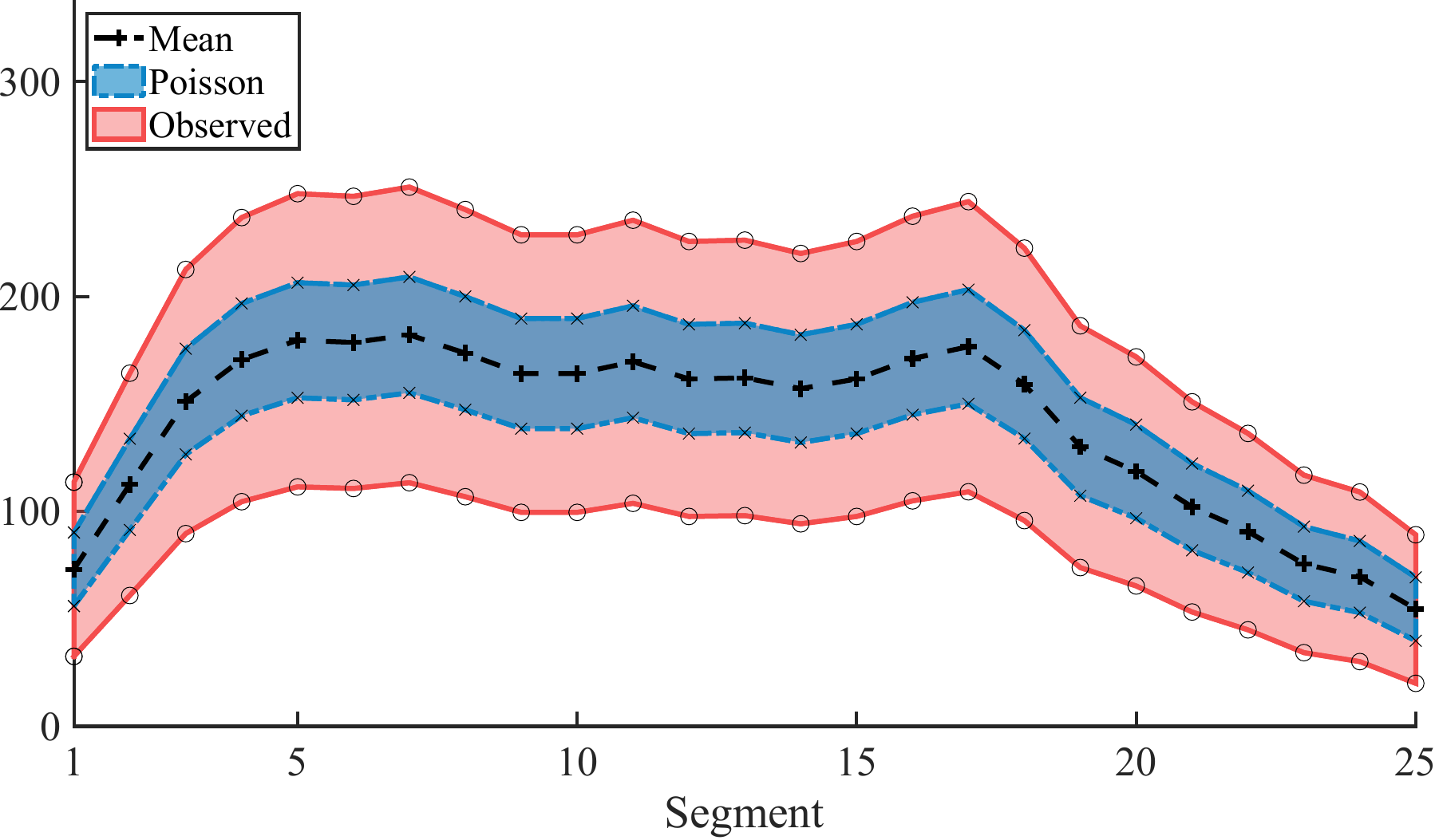} 
  }
  {Daily Arrival Pattern of the Fitted ADL Model. \label{fig:time-of-day-ADL}}
  {The black dashed line represents the mean arrival count of each 30-minute segment.  The half-width of the shaded areas represents two standard deviations. 
  The shaded area bounded by the red solid line represents the fitted ADL model, while the area bounded by the blue dash-dotted line corresponds the Poisson-distributed arrival counts. }
\end{figure}

\cite{Avramidis04_ec} fitted this model to arrival data from a Bell Canada call center operating 12.5 hours daily, dividing the day into $k=25$ segments of $\Delta=0.5$ hours each. Using their setup and estimated parameters, we present in Figure~\ref{fig:time-of-day-ADL} the means and standard deviations of arrival counts $(Z_1, \ldots, Z_k)$ under the fitted ADL model.

Using 1000 independent samples of the arrival count vector $(Z_1,\ldots,Z_k)$, we fit the non-stationary version of models $\mathcal{M}_1$ through $\mathcal{M}_5$ to this data and compute segment-specific staffing levels using each model's safety rule. To evaluate each safety rule's performance, we simulate $M=1000$ sample paths of customer arrival times from the ADL model, assuming an empty initial system and log-normal service times with mean and standard deviation both equal to $1/6$. The performance metric is the realized delay probability, estimated as
\begin{align}\label{eq:delay-prob-ADL}
\pr\big( \widetilde{Q}(t) > n(t) \big)   \approx  \frac{1}{M} \sum_{j=1}^{M} \mathbb{I}\big\{\widetilde{Q}_j( t ) > n(t) \big\}, \quad t\in[0, T],
\end{align}
where $\widetilde{Q}_j(t)$ is the $j$-th realization of the number-in-system process $\widetilde{Q}(t)$, $n(t)$ is the staffing level at time $t$ determined by each model's safety rule, and $T = 12.5 $ hours.

\begin{figure}[t]
  \FIGURE{
    \includegraphics[width=0.5\textwidth]{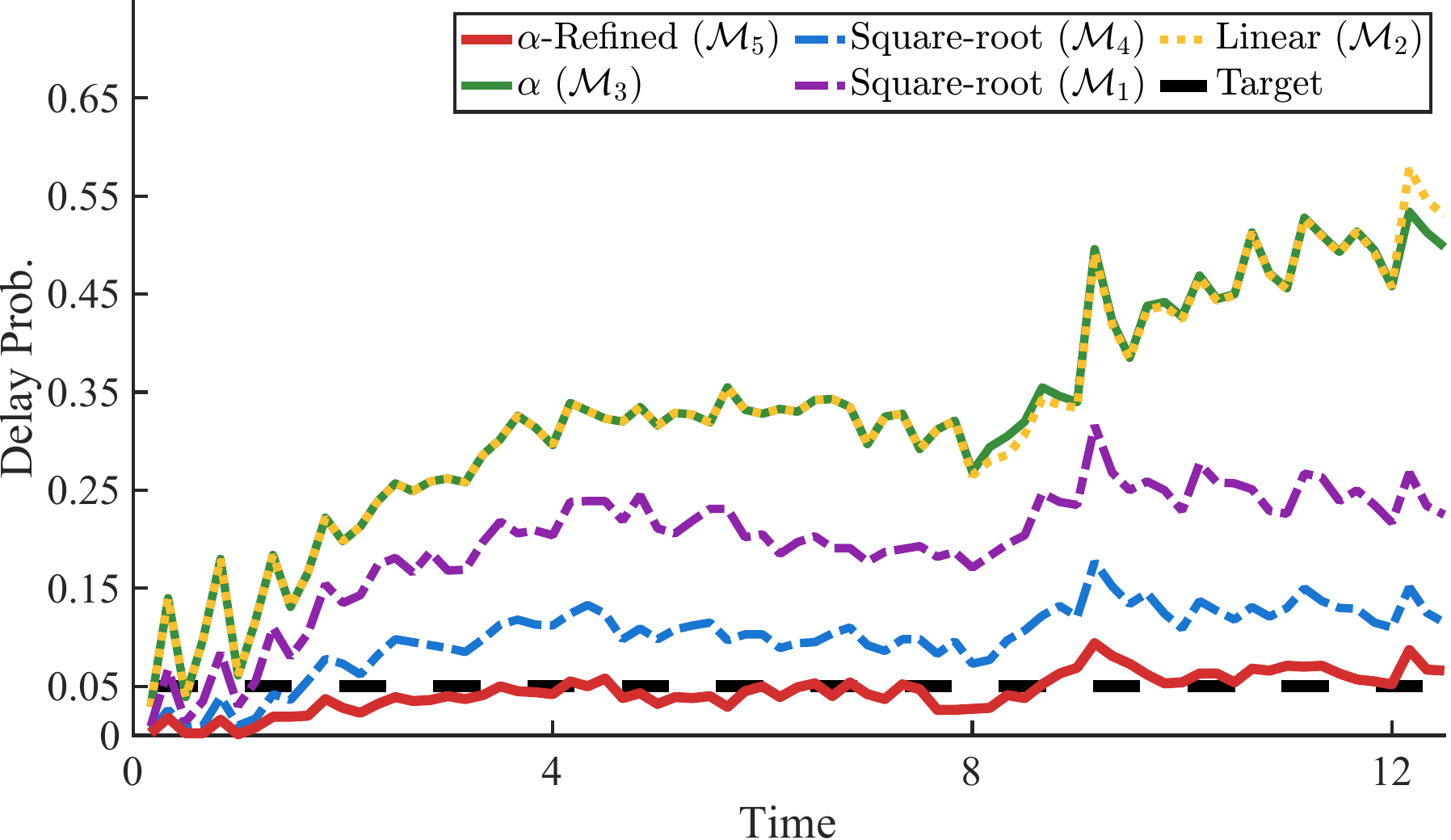} 
    \includegraphics[width=0.5\textwidth]{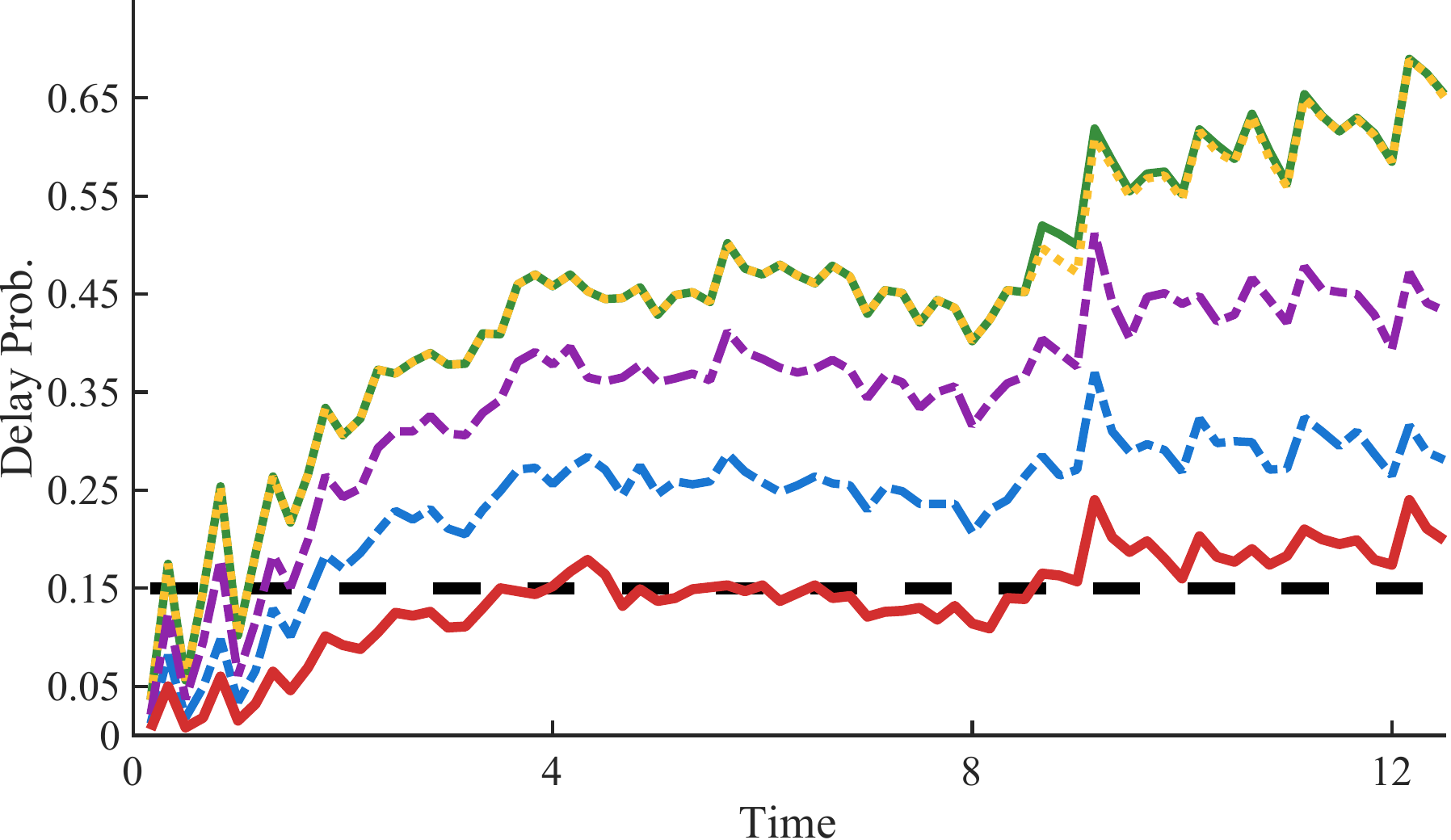} 
  }
  {Performance of Staffing Rules in Finite-server Systems with the ADL Model. \label{fig:Model-misspec}}
{Models $\mathcal{M}_1$ through $\mathcal{M}_5$ are the non-stationary version with piecewise-constant mean arrival rates. The delay probability is defined by Equation~\eqref{eq:delay-prob-ADL}. Left: $\varepsilon = 0.05$. Right: $\varepsilon = 0.15$.}
\end{figure}

Figure~\ref{fig:Model-misspec} shows that the refined alpha safety rule under our DSPP model $\mathcal{M}_5$ significantly outperforms other safety rules. Starting from zero (due to the empty initial system), its delay probability quickly approaches the target and stays around it for most of the simulated time period. This strong performance, given the ADL model as ground truth, demonstrates our DSPP model's robustness to model misspecification. 
In contrast, safety rules under models $\mathcal{M}_1$ through $\mathcal{M}_4$ show delay probabilities that substantially deviate from the target, primarily due to their more severe model misspecification relative to $\mathcal{M}_5$.

\subsection{Data Description of NYC 311 Call Center}
\label{sec:data-description}

The NYC 311 Call Center is a large-scale service system providing 24/7 public access to non-emergency government services and information in New York City. 
Its phone call arrival data is publicly available  on the NYC OpenData website (\url{https://portal.311.nyc.gov}), with each record containing seven attributes: a unique identifier (UNIQUE\_ID); the date (CALENDAR\_DATE) and time (TIME) the call was received; the associated agency (AGENCY\_NAME), such as Department of Finance or Metropolitan Transportation Authority; the inquiry topic (INQUIRY\_NAME) used to resolve the user's inquiry; a brief description of the topic (BRIEF\_DESCRIPTION); and the call resolution method (CALL\_RESOLUTION), such as ``Information Provided'' or ``Caller Hung Up.''

To illustrate, consider a call received on August 31, 2017, at 11:43:09 AM with ID 125728858. The inquiry topic was ``Status of a Birth Certificate Order,'' and the associated agency was the Department of Health and Mental Hygiene. The caller requested ``the status of a birth certificate order,'' and the call was resolved by ``Transfer to City Agency.''

\subsection{Alternative Staffing Intervals}

Call centers typically use short staffing intervals under one hour \citep{GreenKolesarWhitt07_ec,ShenHuang08_ec}. While our main analysis of the NYC 311 Call Center (Section~\ref{sec:casestudy}) used 30-minute intervals ($k=48$ segments each day for 24-hour operation), 
here we examine staffing rule performance with alternative 10-minute and 1-hour intervals, maintaining all other experimental specifications.

\begin{figure}[t]
\FIGURE{
    \includegraphics[width=0.5\textwidth]{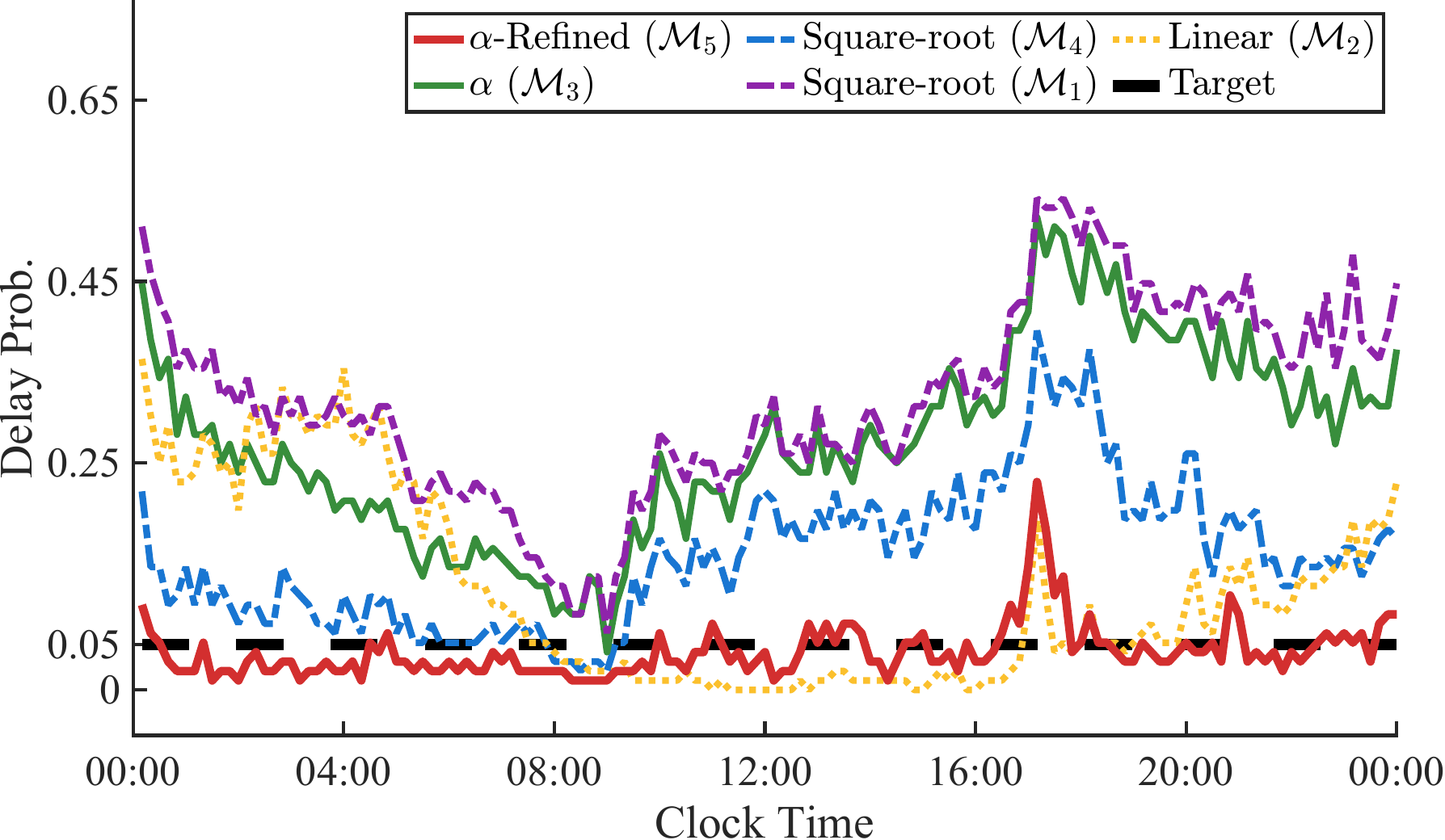}
	\includegraphics[width=0.5\textwidth]{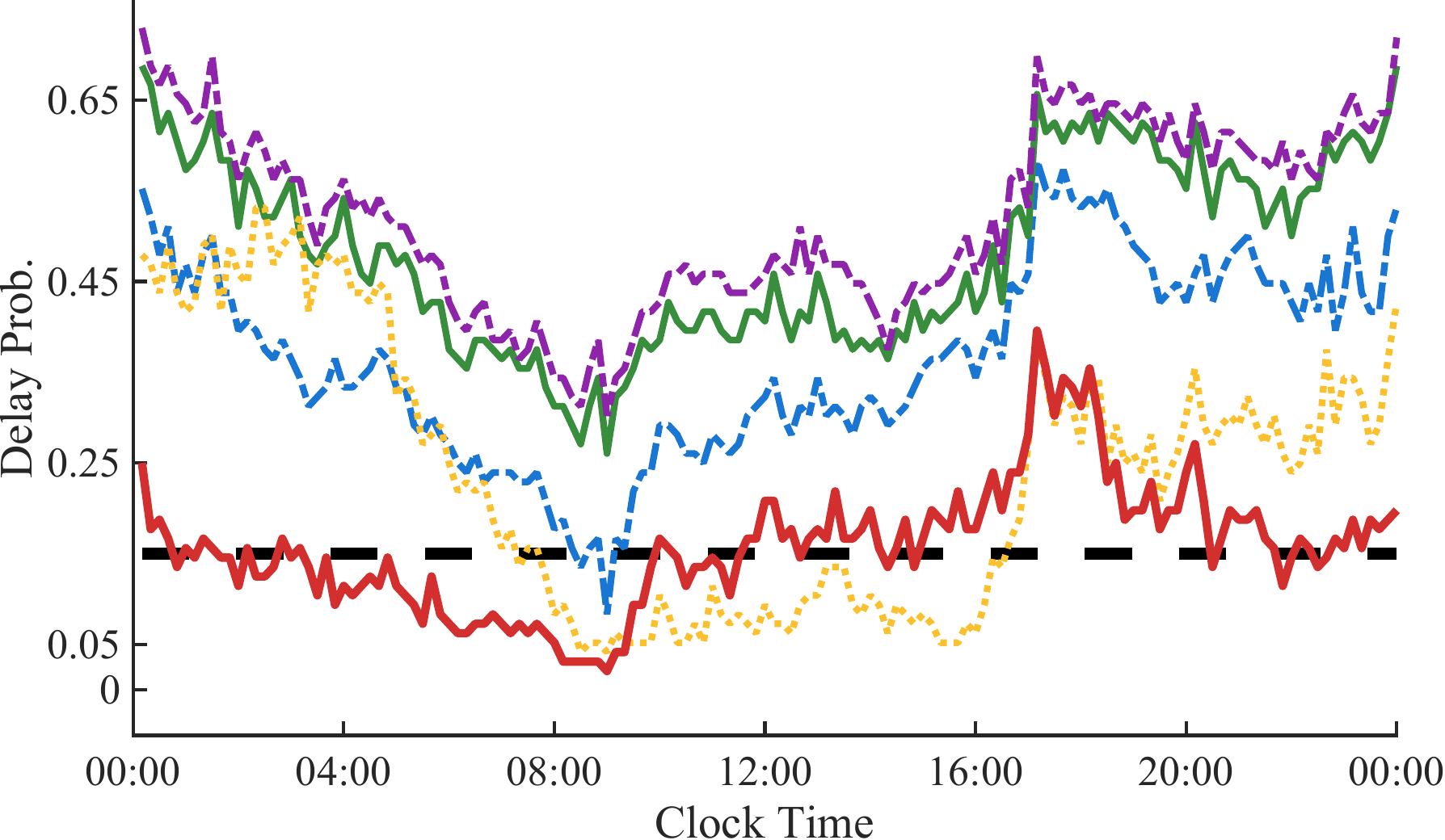}
	}
{Performance of Staffing Rules in Finite-server Systems with Real Arrival Data ($\Delta = 10$ Minutes). \label{figsim212}}
{}
\FIGURE{
    \includegraphics[width=0.5\textwidth]{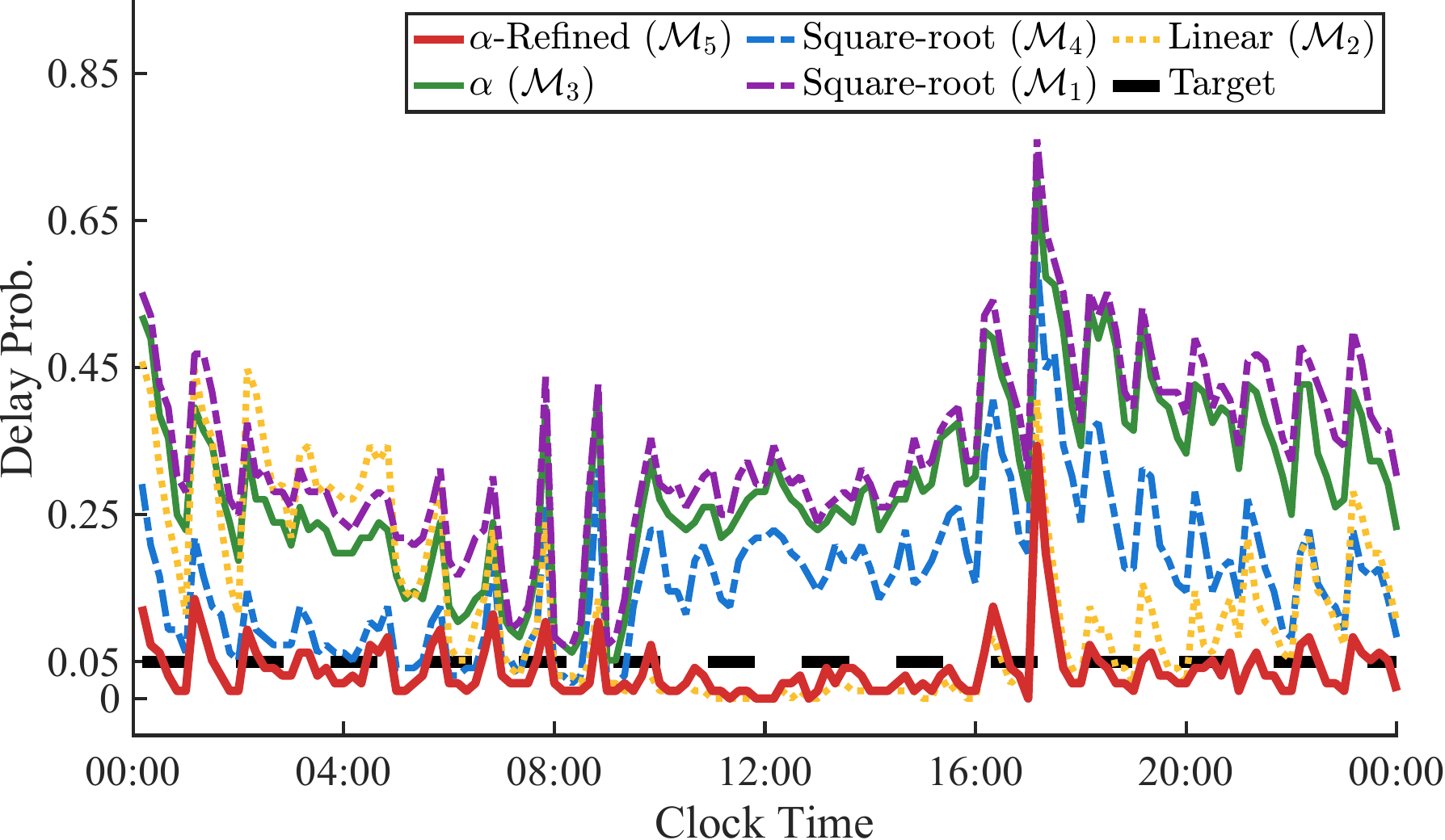}
	\includegraphics[width=0.5\textwidth]{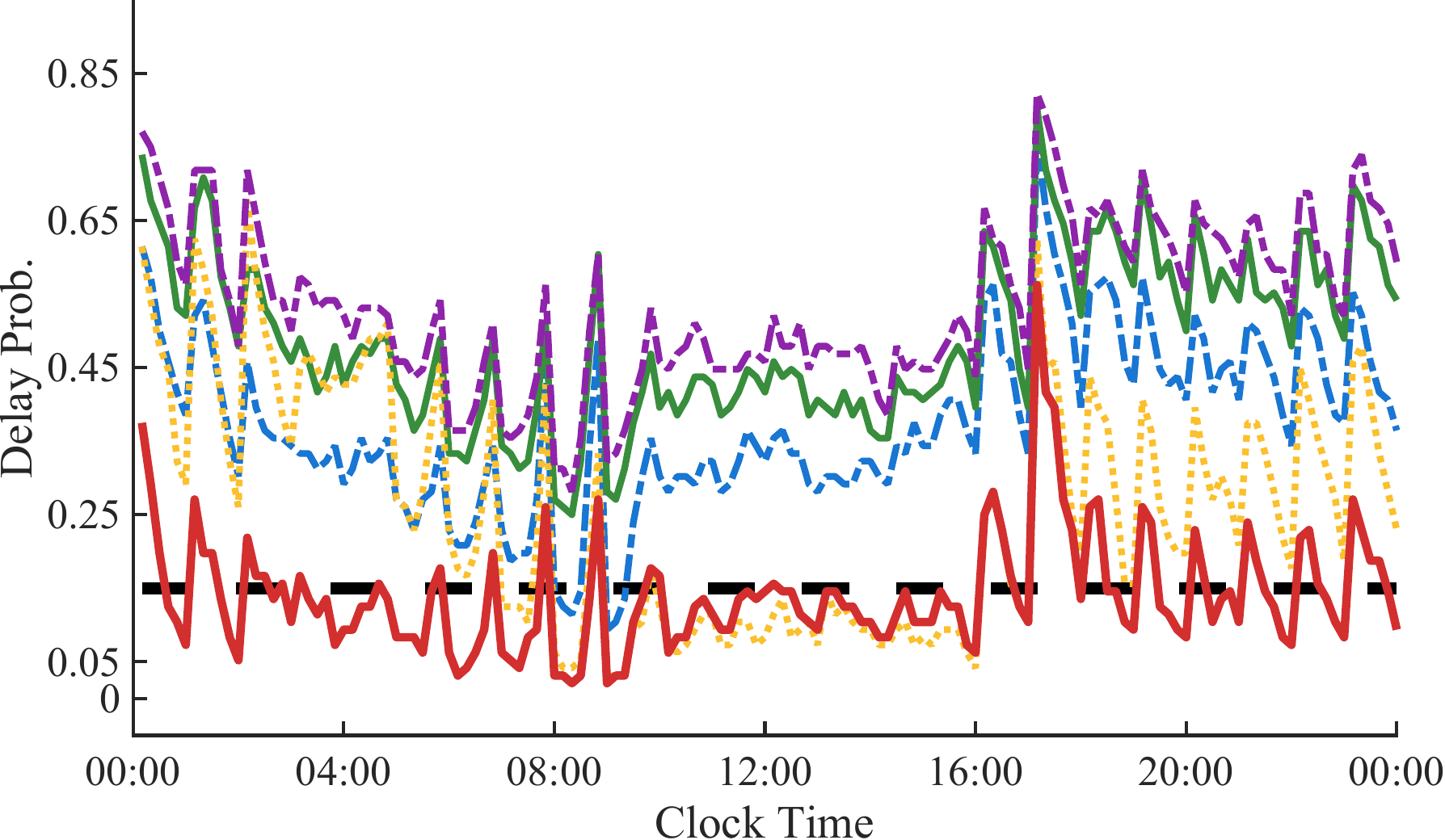}
	}
{Performance of Staffing Rules in Finite-server Systems with Real Arrival Data ($\Delta = 1$ Hour). \label{figsim2122}}
{Models $\mathcal{M}_1$ through $\mathcal{M}_5$ represent non-stationary processes with piecewise-constant mean arrival rates, with staffing levels computed for each 30-minute segment. The delay probability follows Equation~\eqref{delayprob00}, shown for $\varepsilon = 0.05$ (left) and $\varepsilon = 0.15$ (right). Model parameters are estimated using nighttime data (9 PM--8 AM) from January--May 2017, while delay probabilities are evaluated using full-day data from July--November 2017.  Each $\lambda_i$ is estimated using the \textsf{SIPP~Avg} approach (Appendix~\ref{sec:MLE}).}
\end{figure}

Results appear in Figure~\ref{figsim212} (10-minute intervals) and Figure~\ref{figsim2122} (1-hour intervals), both using the \textsf{SIPP Avg} estimator for the mean arrival rate process $\lambda(t)$. 
We do not use the \textsf{SIPP Mix} estimator as we did in Section~\ref{sec:casestudy} for the following reason. 
For 10-minute segments, \textsf{SIPP Min} and \textsf{SIPP Max} yield negligible differences from \textsf{SIPP Avg}, as the mean arrival rate changes little over such short periods. However, for 1-hour segments, using \textsf{SIPP Min} or \textsf{SIPP Max} would introduce substantial bias, as the mean arrival rate may change significantly within an hour. 

Consistent with the 30-minute interval results (Figure~\ref{figsim2120}), the refined alpha safety rule ($\mathcal{M}_5$) outperforms other staffing rules for both 10-minute and 1-hour intervals. All staffing rules perform better with shorter intervals. With 1-hour staffing intervals, even the refined alpha safety rule shows significant fluctuations, primarily because infrequent staffing adjustments cannot adequately respond to rapid changes in mean arrival rates.

\subsection{Gamma-distributed Service Times}
\label{app:sec:gamma}

The experiments in Sections~\ref{sec:synthetic} and~\ref{sec:casestudy} assume log-normally distributed service times. Here we extend our analysis to gamma-distributed service times, setting both mean and standard deviation to $1/6$ hours (equivalent to $\mu=6$ customers per hour), matching the parameters of the log-normal case.

\begin{figure}
\FIGURE{
    \includegraphics[width=0.5\textwidth]{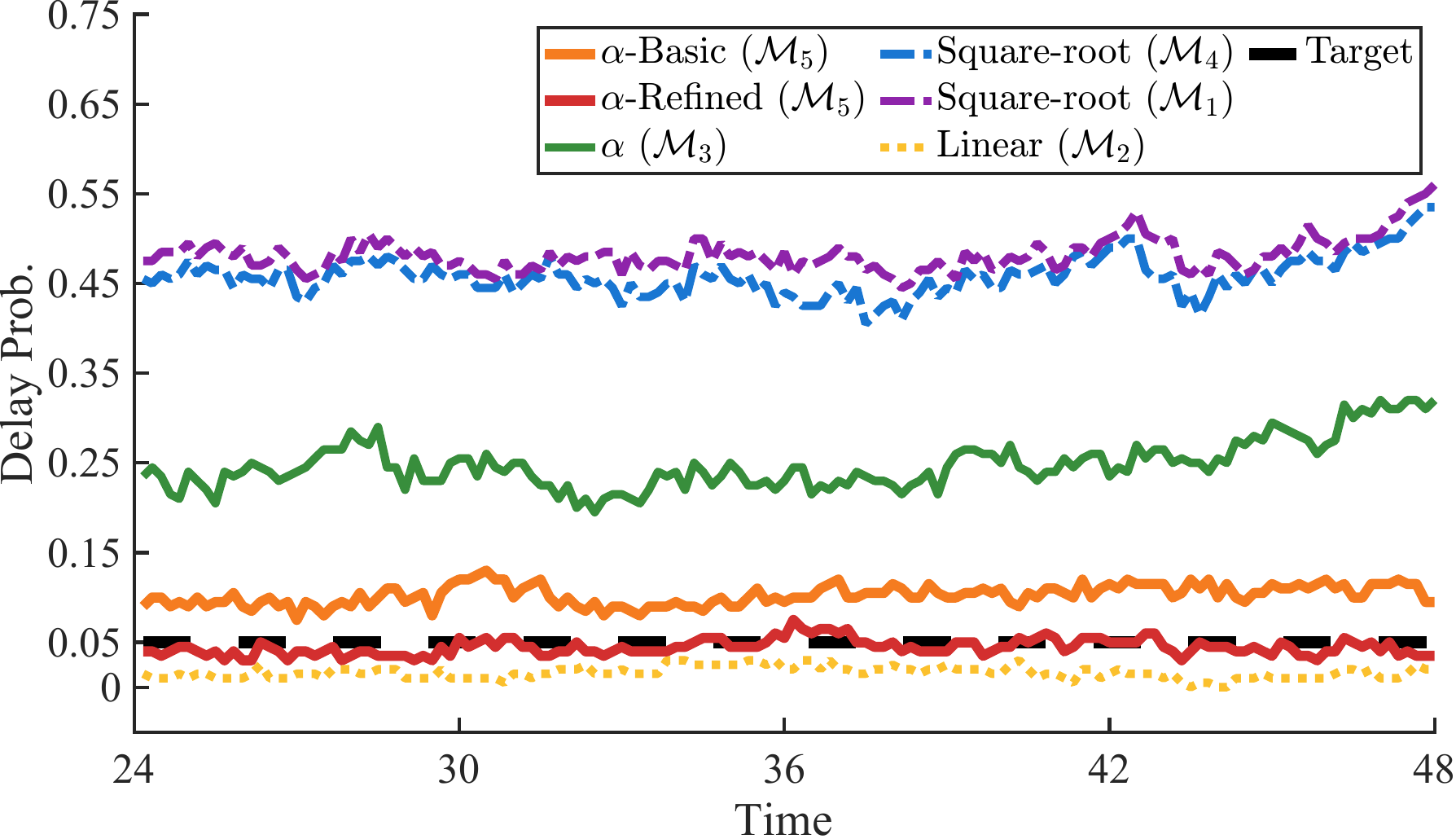}
	\includegraphics[width=0.5\textwidth]{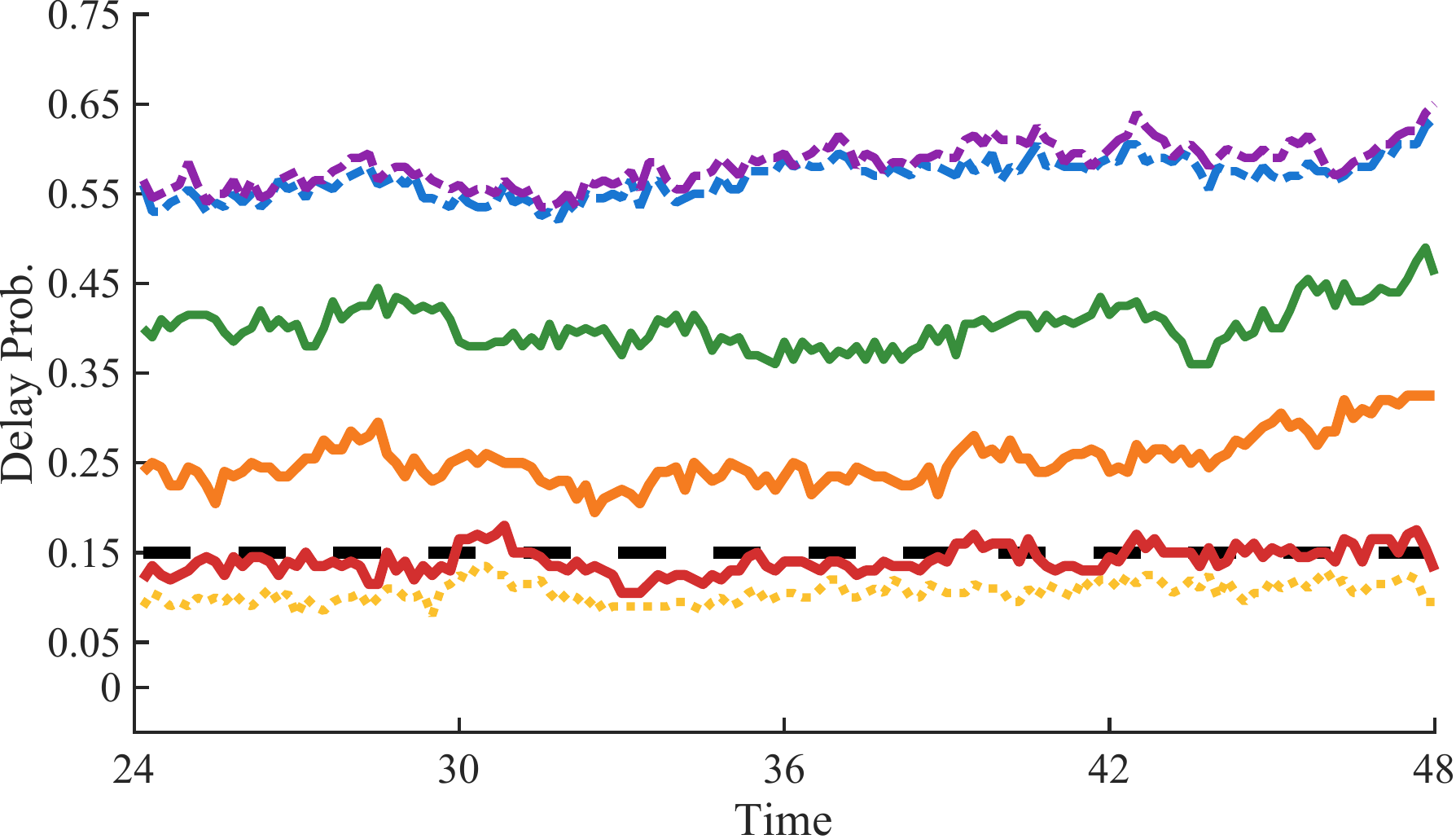}
	}
{Performance of Staffing Rules in Finite-server Systems with Arrival Model $\mathcal{M}_5$ and Gamma Service Times. \label{fig:Gamma+Finite+Sim}}
{The delay probability as defined by 
$\pr( Q_\lambda^{(n)}(t) > n)$, where $n$ varies depending on the staffing rule. 
The initial 24 hours serve as the warm-up period and thus, they are excluded. System scale: $\lambda=2400$. 
Left: $\varepsilon = 0.05$. Right: $\varepsilon = 0.15$.}
\end{figure}

\begin{figure}[t]
\FIGURE{
    \includegraphics[width=0.5\textwidth]{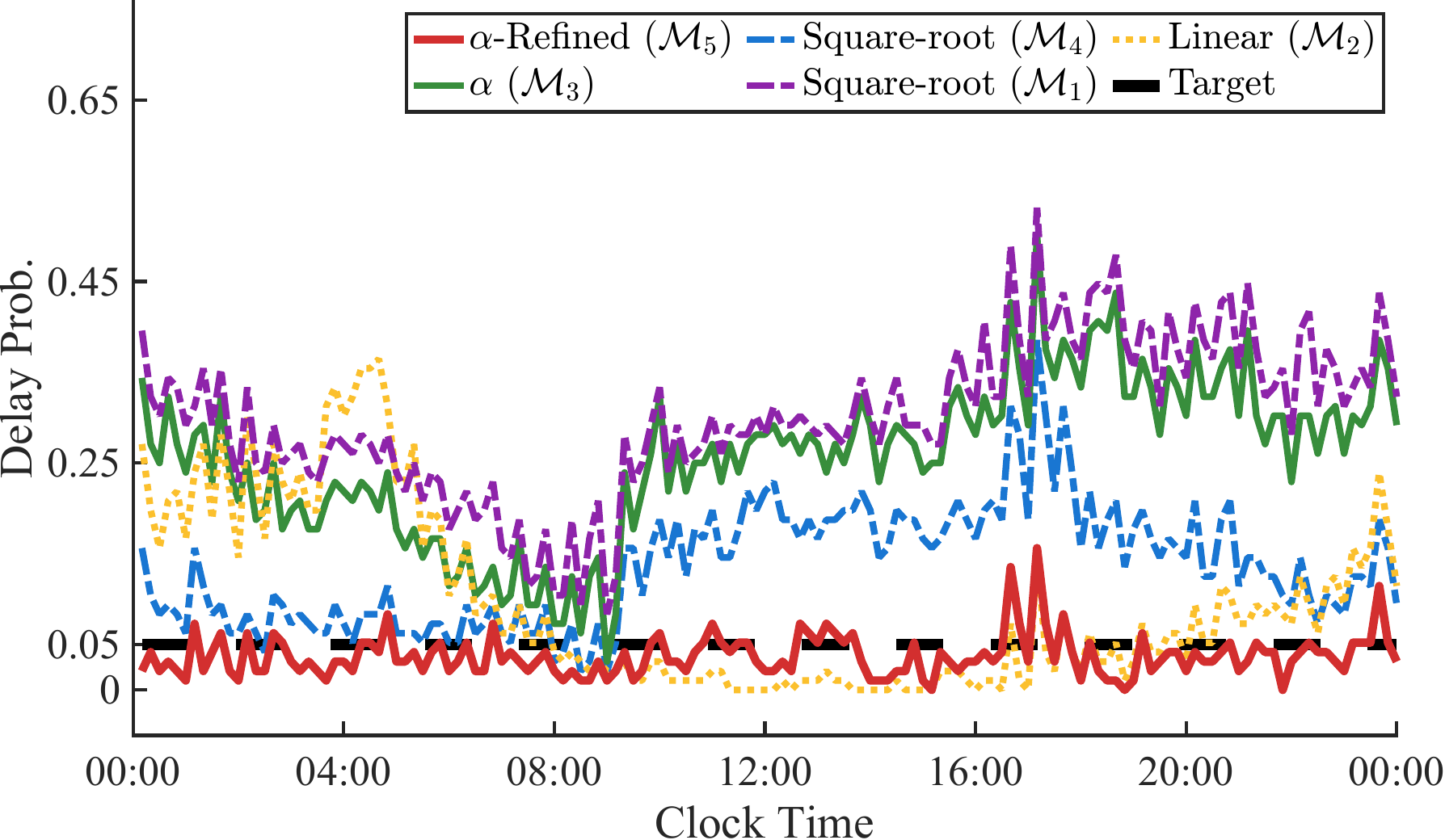}
	\includegraphics[width=0.5\textwidth]{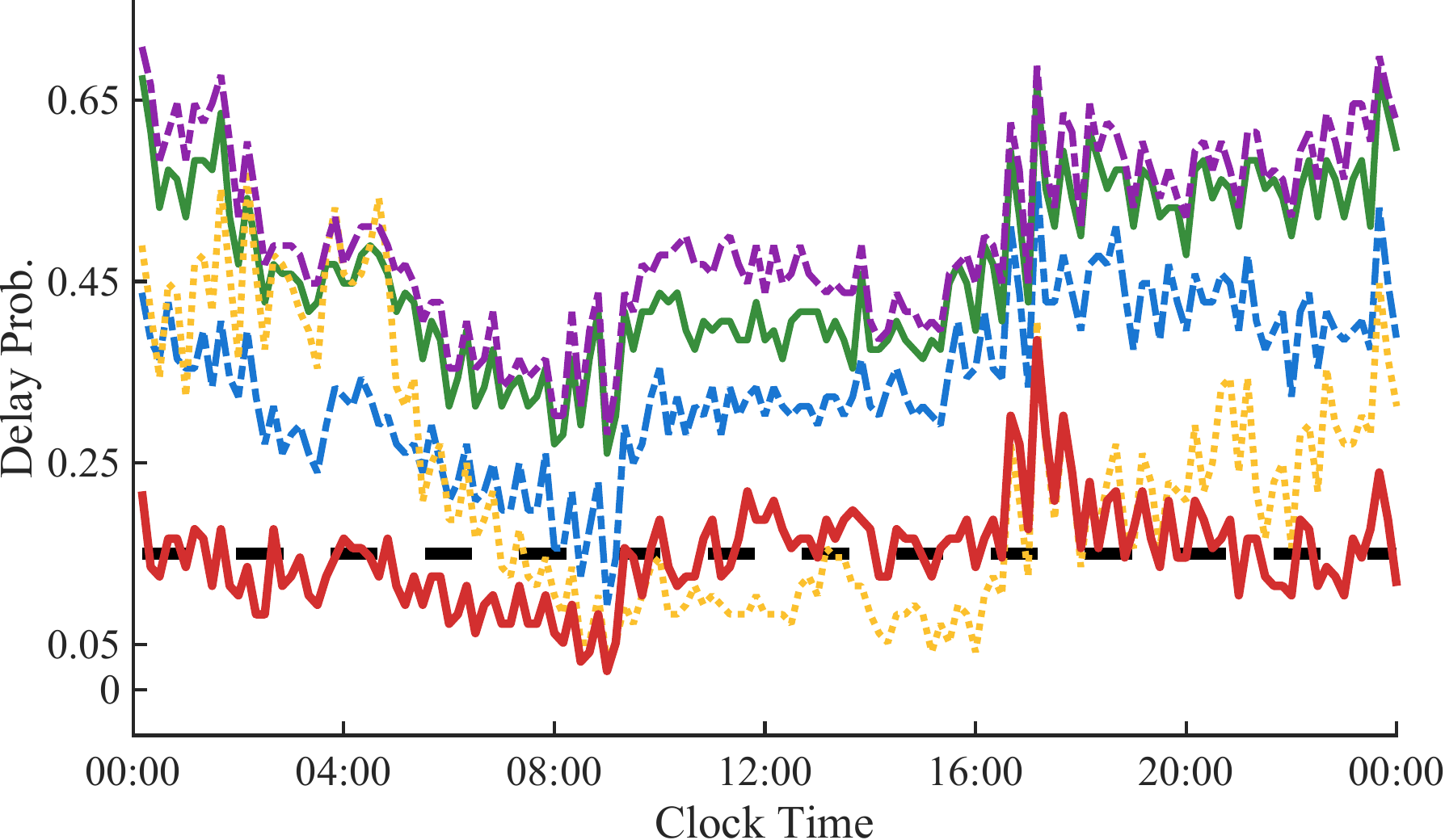}
	}
{Performance of Staffing Rules in Finite-server Systems with Real Arrival Data and Gamma Service Times (SIPP Avg). \label{fig:Gamma+Finite+Real+Avg}}
{}
\FIGURE{
    \includegraphics[width=0.5\textwidth]{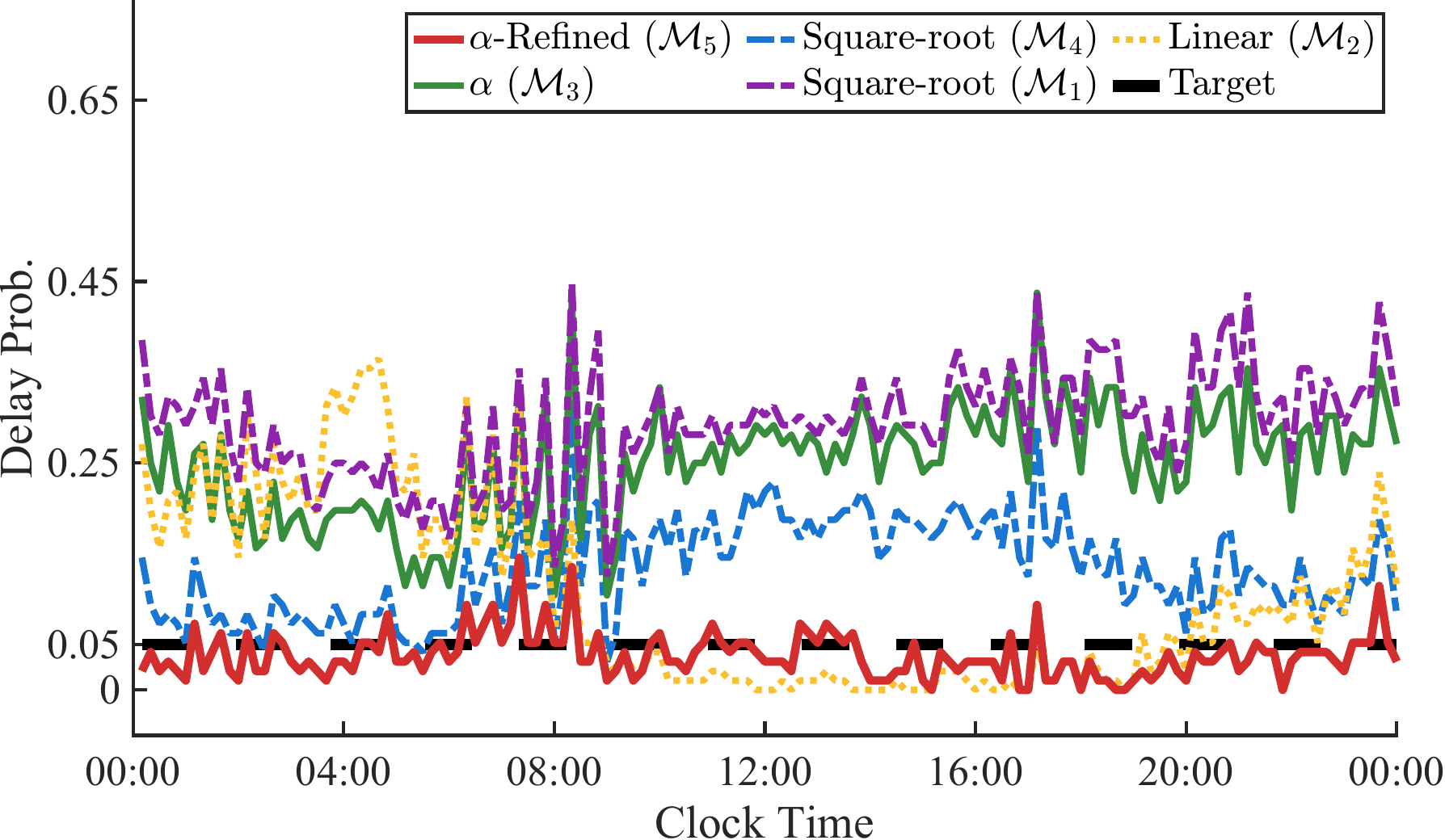}
	\includegraphics[width=0.5\textwidth]{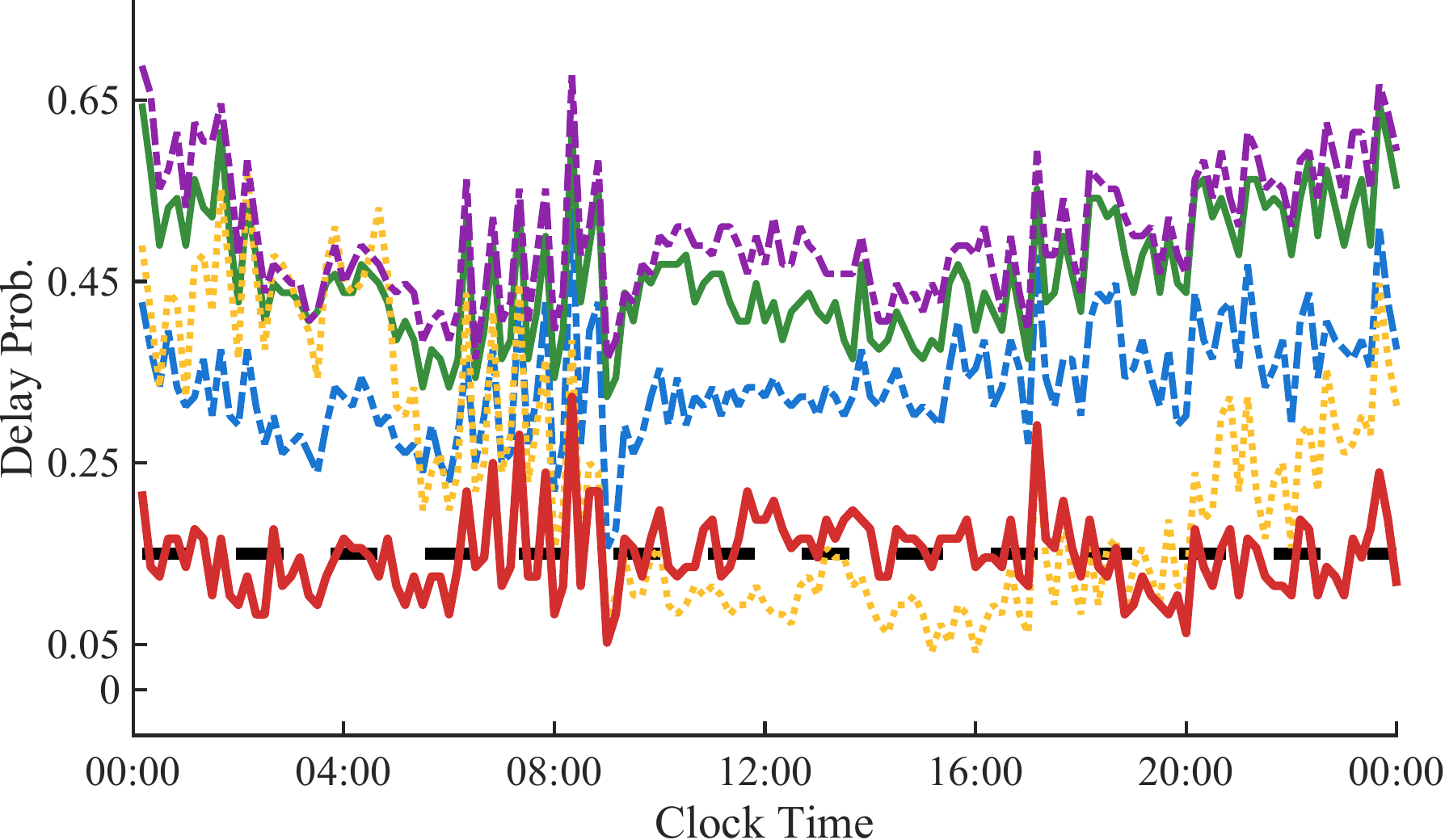}
	}
{Performance of Staffing Rules in Finite-server Systems with Real Arrival Data and Gamma Service Times (SIPP Mix). \label{fig:Gamma+Finite+Real+Mix}}
{Models $\mathcal{M}_1$ through $\mathcal{M}_5$ represent non-stationary processes with piecewise-constant mean arrival rates, with staffing levels computed for each 30-minute segment. The delay probability follows Equation~\eqref{delayprob00}, shown for $\varepsilon = 0.05$ (left) and $\varepsilon = 0.15$ (right). Model parameters are estimated using nighttime data (9 PM--8 AM) from January--May 2017, while delay probabilities are evaluated using full-day data from July--November 2017.  Each $\lambda_i$ is estimated using the \textsf{SIPP~Avg} and \textsf{SIPP~Mix} approach, respectively (Appendix~\ref{sec:MLE}).}
\end{figure}

We conduct two experiments: finite-server systems with stationary arrivals from our DSPP model $\mathcal{M}_5$ (Section~\ref{sec:fin-ser-sym}), and finite-server systems using NYC 311 Call Center data (Section~\ref{sec:casestudy}). Results appear in Figure~\ref{fig:Gamma+Finite+Sim} for the stationary arrivals and Figures~\ref{fig:Gamma+Finite+Real+Avg} and \ref{fig:Gamma+Finite+Real+Mix} for the call center data. These gamma-distribution results mirror the log-normal findings, with the refined alpha safety rule maintaining its superior performance.

\end{document}